\documentclass[12pt,reqno]{amsart}
 \author{Kai (Steve) Fan}
\author{Paul Pollack}
\address{Department of Mathematics\\ University of Georgia\\ Athens, GA 30602}
\email{Steve.Fan@uga.edu}
\email{pollack@uga.edu}
\usepackage{mathtools}
\usepackage{empheq}

\newcommand{\verteq}{\rotatebox{90}{$\,=$}}
\newcommand{\vertcong}{\rotatebox{-90}{$\cong\,$}}

\newcommand{\equalto}[2]{\underset{\scriptstyle\overset{\mkern4mu\verteq}{#2}}{#1}}
\newcommand{\conggto}[2]{\underset{\scriptstyle\overset{\mkern4mu\vertcong}{#2}}{#1}}

\title{The typical elasticity of a quadratic order}
\usepackage{amsmath,amssymb,amsthm,microtype,geometry,mathscinet,enumitem,url}
\usepackage{parskip}
\usepackage{stackengine}
\usepackage{colonequals}
\usepackage{xcolor}

\renewcommand\supset\supseteq
\renewcommand\epsilon\varepsilon
\renewcommand\phi\varphi
\usepackage[utf8]{inputenc}
\newcommand{\leg}[2]{\genfrac{(}{)}{}{}{#1}{#2}}

\makeatletter
\NewCommandCopy\@@pmod\pmod
\DeclareRobustCommand{\pmod}{\@ifstar\@pmods\@@pmod}
\def\@pmods#1{\mkern4mu({\operator@font mod}\mkern 6mu#1)}
\makeatother

\DeclareMathAlphabet{\curly}{U}{rsfs}{m}{n}
\newtheorem{thm}{Theorem}[section]

\newtheorem{prop}[thm]{Proposition}
\newtheorem{lem}[thm]{Lemma}

\newtheorem*{thmcheb}{GRH-dependent Chebotarev density theorem}

\theoremstyle{remark}

\DeclareMathOperator{\Rad}{Rad}
\DeclareMathOperator{\Exp}{Exp}
\DeclareMathOperator{\RadExp}{RadExp}
\DeclareMathOperator{\Nm}{Nm}
\DeclareMathOperator{\Tr}{Tr}
\DeclareMathOperator{\Dav}{Dav}
\DeclareMathOperator{\Var}{Var}
\DeclareMathOperator{\Cov}{Cov}
\DeclareMathOperator{\sgn}{sgn}
\usepackage{enumitem}

\begin{document}

\def\Cc{\curly{C}}
\def\Ll{\mathcal{L}}
\def\Ee{\curly{E}}
\def\Oo{\mathcal{O}}
\def\Aa{\mathcal{A}}
\def\I{\mathcal{I}}
\def\J{\mathcal{J}}
\def\N{\mathbb{N}}
\def\Q{\mathbb{Q}}
\def\Qq{\mathcal{Q}}
\def\Rr{\mathcal{R}}
\def\Cl{\mathrm{Cl}}
\def\radi{\mathrm{rad}}
\def\PrinCl{\mathrm{PrinCl}}
\def\PreCl{\mathrm{PreCl}}
\newcommand\rad{\mathrm{radexp}}
\def\E{\mathbb{E}}
\def\Z{\mathbb{Z}}
\def\F{\mathbb{F}}
\def\PP{\mathbb{P}}
\def\R{\mathbb{R}}
\def\C{\mathbb{C}}
\def\Gal{\mathrm{Gal}}
\def\Dd{\mathcal{D}}
\def\Pp{\mathcal{P}}
\def\Ss{\mathcal{S}}
\newcommand\Frob{\mathrm{Frob}}
\def\lcm{\mathop{\mathrm{lcm}}}
\renewcommand\subset\subseteq
\newcommand\Li{\mathrm{Li}}
\newcommand\altxrightarrow[2][0pt]{\mathrel{\ensurestackMath{\stackengine%
  {\dimexpr#1-7.5pt}{\xrightarrow{\phantom{#2}}}{\scriptstyle\!#2\,}%
  {O}{c}{F}{F}{S}}}}
\newcommand{\congto}{\altxrightarrow{\sim}}

\begin{abstract} For an atomic domain $D$, the \emph{elasticity} $\rho(D)$ of $D$ is defined as $\sup\{r/s: \pi_1\cdots \pi_r = \rho_1 \cdots \rho_s,~ \text{where each $\pi_i, \rho_j$ is irreducible}\}$; the elasticity provides a concrete measure of the failure of unique factorization in $D$. Fix a quadratic number field $K$ with discriminant $\Delta_K$, and for each positive integer $f$, let $\Oo_f = \Z + f\Oo_K$ denote the order of conductor $f$ in $K$. Results of Halter-Koch imply that $\Oo_f$ has finite elasticity precisely when $f$ is \emph{split-free}, meaning not divisible by any rational prime $p$ with $(\Delta_K/p)=1$. When $K$ is imaginary, we show that for almost all split-free $f$,
\[ \rho(\Oo_f) = f/(\log{f})^{\frac{1}{2}\log\log\log{f} + \frac{1}{2}C_K+o(1)}, \]
for a constant $C_K$ depending on $K$. When $K$ is real, we prove under the assumption of the Generalized Riemann Hypothesis that \[ \rho(\Oo_f)= (\log{f})^{\frac12 +o(1)} \]
for almost all split-free $f$. Underlying these estimates are new statistical theorems about class groups of orders in quadratic fields, whose proofs borrow ideas from investigations of Erd\H{o}s, Hooley, Li, Pomerance, Schmutz, and others into the multiplicative groups $(\Z/m\Z)^\times$. One novelty of the argument is the development of a weighted version of the Tur\'{a}n--Kubilius inequality to handle a variety of sums over split-free integers.
\end{abstract}
\subjclass[2020]{Primary 11R27; Secondary 11N25, 11N37, 11R11, 11R65, 13A05}

\maketitle

\section{Introduction}
Let $K$ be a number field. It was realized in the 19th century that the ring of integers $\Oo_K$ of $K$ may fail to be a unique factorization domain. In such cases, Dedekind showed that unique factorization is restored by looking at ideals of $\Oo_K$ rather than elements. From this ideal-theoretic perspective, the failure of elementwise unique factorization is due to the principal ideals failing to ``fill out'' the space of all ideals, with the deficit captured concretely by the class group $\Cl(\Oo_K)$ of $K$. The finiteness of the class group can be interpreted as stating that $\Oo_K$ is always ``a finite distance away'' from possessing unique factorization.

The term ``distance'' is used informally here, but it is reasonable to wonder whether it can be made precise. Is there a nonnegative real number we can assign to a domain $D$ measuring the failure of unique factorization? If $D = \Oo_K$, one natural candidate is the class number $h_K$ of $K$. But there is an alternative and in some ways more appealing candidate for such a measure, foreshadowed in work of Carlitz \cite{carlitz60} and discussed explicitly by Valenza \cite{valenza90},\footnote{While Valenza's paper did not appear until 1990, it was submitted in 1980!} Steffan \cite{steffan86}, and Narkiewicz \cite{narkiewicz95}.
Let $D$ be an atomic domain, meaning an integral domain in which each nonzero nonunit possesses a factorization into irreducible elements. Following Valenza, the \emph{elasticity} $\rho(D)$ of $D$ is defined by
\[ \rho(D) = \sup\left\{\frac{r}{s}: \pi_1\cdots \pi_r = \rho_1\cdots \rho_s, \text{ all $\pi_i, \rho_j$ irreducible}\right\}. \]
That is, we consider all coincidences between a product of $r$ irreducibles and a product of $s$ irreducibles, and we take the supremum of the ratios $\frac{r}{s}$. For example, $\rho(D)=1$ if and only if the \emph{length} of a factorization into irreducibles (length meaning the number of factors in the product) is uniquely determined by the element being factored. In this case, $D$ is called a \emph{half-factorial domain} or HFD.  
Both half-factoriality and elasticity have been extensively investigated from the viewpoint of commutative algebra.  For an account of these theories as they stood at the end of the 20th century, see the surveys \cite{anderson97} and \cite{CC00}.

A beautiful theorem of Valenza, Steffan, and Narkiewicz (op.~cit.) explicitly determines the elasticity of $\rho(\Oo_K)$ for every number field $K$. Let $G$ be a finite abelian group, written multiplicatively with identity element $1$. Call a finite sequence $g_1,\dots,g_n$ of elements of $G$  a \emph{$1$-product sequence} if $g_1\cdots g_n = 1$. Say $g_1,\dots,g_n$ is \emph{$1$-product-free} if no nonempty subsequence of $g_1,\dots,g_n$ is a $1$-product sequence. The \emph{Davenport constant $\Dav G$} of $G$ is the largest integer $D$ for which there is a $1$-product sequence $g_1,\dots,g_D$ with no nonempty, proper $1$-product subsequence. It is straightforward to show that $\Dav G$ is finite (in fact, that $\Dav G \le \#G$) and that $(\Dav G)-1$ is the length of the longest $1$-product-free sequence of elements of $G$. (See, for instance, \cite[Lemma 1.4.9]{GHK06}.) Valenza, Narkiewicz, and Steffan's elegant result asserts that 
\[  \rho(\Oo_K) = \max\left\{1, \frac{1}{2} \Dav \Cl(\Oo_K)\right\}. \]
A particularly attractive corollary is that $\Oo_K$ is half-factorial precisely when $h_K \le 2$. This special case of the Valenza--Steffan--Narkiewicz theorem was already shown by Carlitz in 1960 \cite{carlitz60}. 

It is of great interest to understand the extent to which unique factorization in $\Oo_K$ fails as $K$ varies across a family of number fields. For the collection of quadratic fields $K$, such questions can be seen as implicit in Gauss's \emph{Disquisitiones} (couched in the language of binary quadratic forms). Two hundred years post-\emph{Disquisitiones}, such problems continue to elude us. To mention one concrete instantiation of our ignorance, though it is universally believed that there are infinitely many real quadratic fields of class number $1$, we cannot yet \underline{dis}prove that the class numbers of real quadratic fields tend to infinity with the discriminant. 

In this paper we consider a natural-seeming orthogonal family of problems. Specifically, we fix a quadratic field $K$ with discriminant $\Delta$ and examine the factorization behavior as we range over all orders inside $K$. Here an \emph{order} in $K$ is a subring of $\Oo_K$ strictly containing $\Z$.\footnote{For not-necessarily-quadratic number fields $K$, one needs a more complicated definition: An \emph{order in $K$} is a subring of $\Oo_K$ containing a $\Q$-basis for $K$. } Orders in a quadratic field are parametrized by the positive integers in a natural way: For each $f\in \N$, there is a unique order $\Oo_f$ with index $[\Oo_K:\Oo_f] = f$, namely $\Oo_f = \Z + f\Oo_K$; that is,
\[ \Oo_f = \{\alpha \in \Oo_K:  \alpha \equiv a\pmod*{f\Oo_K}\text{ for some $a \in \Z$}\}. \]
Furthermore, every order $\Oo$ is one of the $\Oo_f$, with the maximal order $\Oo_K = \Oo_1$. The integer $f$ is referred to as the \emph{conductor} of $\Oo_f$. 

In the late 1970s, Zaks \cite{zaks76, zaks80} observed that $\Z[\sqrt{-3}]$ --- the order of conductor $2$ in $\Q(\sqrt{-3})$ --- is a half-factorial domain. Later, Halter-Koch \cite{HK83} and Coykendall \cite{coykendall01} wrote down necessary and sufficient algebraic conditions for an order in a quadratic field to be half-factorial. In the same paper \cite{coykendall01} of Coykendall, it was shown that $\Z[\sqrt{-3}]$ is the only half-factorial nonmaximal order in an imaginary quadratic field. The situation for real quadratic fields is more complicated. In \cite{coykendall01} it is conjectured that (a) if one varies both the real quadratic field $K$ and the conductor $f$, then $\Oo_f$ is a half-factorial domain infinitely often, and (b) (stronger) fixing $K=\Q(\sqrt{2})$ and varying only $f$, one still finds infinitely many half-factorial orders.

In \cite{PP23}, 
Pollack
proved form (a) of Coykendall's conjecture. The stronger form (b) is shown subject to the assumption of the Generalized Riemann Hypothesis (GRH).\footnote{Throughout the paper, GRH refers to the Riemann Hypothesis for Dedekind zeta functions.} Both arguments rely on analogues of methods initially introduced to study Artin's primitive conjecture. Some results of \cite{PP23} are extended in \cite{pollackm}.  For instance, it is proved in \cite{pollackm} --- still assuming GRH --- that every element of $\mathcal{E}:= \{1, \frac32, 2, \frac52, 3, \frac72, \dots\} \cup \{\infty\}$ is the elasticity of infinitely many orders of $\Q(\sqrt{2})$. As follows from Lemma \ref{lem:elasticityirreducible} below, an order in a quadratic field always has an elasticity belonging to $\mathcal{E}$, so that $\Z[\sqrt{2}]$ is extremal in a natural sense.

The orders constructed in \cite{PP23} and \cite{pollackm} are quite special. For instance, each of their conductors is composed of at most two distinct primes. In this paper we investigate elasticities corresponding to general conductors $f$.  

Some care is needed to decide what ``general'' should mean. As already mentioned, the elasticity of a quadratic order may be infinite. In \cite[Corollary 4]{HK95}, Halter-Koch shows that if $\Oo$ is an order in a number field $K$, then $\rho(\Oo) < \infty$ precisely when every nonzero prime ideal of $\Oo$ lies below a unique prime ideal of $\Oo_K$. Placed in the setting of quadratic fields $K$, Halter-Koch's theorem says that $\rho(\Oo_f)<\infty$ precisely when $f$ is \emph{split-free}, meaning not divisible by any prime $p$ that splits (completely) in $K$ (cf. the proof of \cite[Theorem 1.1]{pollackm}).

The Selberg--Delange method \cite[Theorems II.5.2, II.5.4]{Ten15}, \cite[Theorem 13.2]{Kou19} (or a mean value theorem of the type discussed in \cite{MC99}) shows that the split-free numbers have a counting function asymptotic to a certain constant multiple of $x/\sqrt{\log x}$. In particular, the split-free numbers $f$ make up a set of density zero. Thus, if one fixes a quadratic field $K$, then $\rho(\Oo_f)=\infty$ for asymptotically 100\% of conductors $f$. This suggests that the proper object of study is not the arithmetic function $\rho(\Oo_f)$, but the corresponding function restricted to split-free $f$. Our main theorems address the almost-everywhere behavior of this restricted function. 

First, we deal with imaginary quadratic orders. Let 
\begin{equation}\label{eq:Dd}
\Dd\colonequals\{-4,\pm8\}\cup\left\{(-1)^{\frac{p-1}{2}}p\colon p>2\text{~is prime}\right\}.   
\end{equation}
Every discriminant of a quadratic field has a unique expression as a product of elements of $\Dd$. 

\begin{thm}\label{thm:imaginary} If $K$ is a fixed imaginary quadratic field with discriminant $\Delta$, then for almost all split-free $f$,
\begin{equation}\label{eq:imagestimate} 
\rho(\Oo_f) = f/(\log{f})^{\frac{1}{2}\log_3{f} + \frac{1}{2}C_K+O((\log_4f)^3/\log_3f)},
\end{equation}
where $\log_k$ denotes the $k$th iterate of the natural logarithm, and
\[C_K\colonequals\sum_{p>2}\frac{\log p}{(p-1)^2}-1-\frac{\sgn(\Delta)1_{\Dd}(\Delta)|\Delta|\log\Rad(|\Delta|)}{\phi(|\Delta|)^2}.\]
\end{thm}
``Almost all'', here and below, means that the estimate holds for all but $o(x/\sqrt{\log{x}})$ split-free numbers $f\le x$, as $x\to\infty$. (When an $O$ appears in the estimate, as is the case here, we are asserting  all of this happens for some choice of implied constant.) In what follows, in place of ``almost all'', we sometimes use the term ``typically'' or the phrases ``almost always'' or  ``for a typical $f$''. It is to be understood that $f$ is restricted to split-free numbers.

For real quadratic orders, we are again able to determine the typical elasticity to within a factor of $(\log{f})^{o(1)}$, but this time we need to assume GRH. The answer is quite a bit smaller than in the imaginary case.

\begin{thm}[conditional on GRH]\label{thm:real} If $K$ is a fixed real quadratic field, then for almost all split-free $f$,
\[ \rho(\Oo_f) = (\log{f})^{\frac{1}{2}+O(1/\log_4{f})}. \]
\end{thm}

In the forthcoming paper \cite{FPextremal}, the authors study the extremal orders of $\rho(\Oo_f)$ restricted to split-free $f$ in both real and imaginary quadratic fields. In addition, it would also seem natural to investigate its average order. Parts of this problem seem attackable by combining ideas of Erd\H{o}s--Pomerance--Schmutz \cite{EPS} with the technology developed in the present paper. We hope to revisit this question on a future occasion.


We now turn to discussing the main ideas needed to establish Theorems \ref{thm:imaginary} and \ref{thm:real}. In \S2, these results are shown to follow from two key propositions, whose proofs occupy \S\S\ref{sec:weightedTK}--\ref{sec:completionkeythm2}.




\section{Theorems 1 and 2: The big picture}\label{sec:bigpicture} For each positive integer $f$, let $I_K(f)$ denote the group of fractional ideals of $K$ generated by integral ideals comaximal with $f\Oo_K$. We let $P_{K,\Z}(f)$ denote the subgroup of $I_K(f)$ generated by principal ideals $\alpha \Oo_K$, where $\alpha\equiv a\pmod{f\Oo_K}$ for some rational integer $a$ with $\gcd(a,f)=1$. Then the \emph{class group $\Cl(\Oo_f)$ of the order $\Oo_f$} is defined to be the quotient $I_K(f)/P_{K,\Z}(f)$. When $f=1$, we will write $I_K$ instead of $I_K(1)$ and $P_K$ instead of $P_{K,\Z}(1)$. Then $I_K$ (resp.~$P_K$) is the group of all fractional (resp. principal fractional) ideals, and the quotient $I_K/P_K$ is the usual class group $\Cl(\Oo_K)$ of the number field $K$.

\subsection{Elasticity in terms of Davenport constants} 
Theorems \ref{thm:imaginary} and \ref{thm:real} are consequences of statistical theorems we establish for the associated class groups. Our starting point is the following lemma, which relates the elasticty $\rho(\Oo_f)$ to the Davenport constant $\Dav \Cl(\Oo_f)$.
\begin{lem}\label{lem:elasticitydavenport} For each split-free $f$,
\[ \frac{1}{2} \Dav \Cl(\Oo_f) \le \rho(\Oo_f) \le \max\left\{\frac{1}{2}\Dav \Cl(\Oo_f) + \frac{3}{2}\Omega(f), 1\right\}.\]
\end{lem}
It will emerge from our later arguments that for any fixed quadratic field $K$, the quantity $\Omega(f)$ is typically of a smaller order than $\Dav \Cl(\Oo_f)$. Hence, $\rho(\Oo_f) \approx \frac{1}{2} \Dav\Cl(\Oo_f)$ most of the time. 

The proof of Lemma \ref{lem:elasticitydavenport} requires some preparation. Lemma \ref{lem:elasticityirreducible} first appeared as \cite[Lemma 2.2]{pollackm}. Proposition \ref{prop:weber} is a classical result of Weber \cite{weber82}; a modern reference is Corollary 2.11.16 on p.\ 159 of \cite{GHK06}. For the statement of the next result, recall that $\Omega(f)$ denotes the number of prime divisors of $f$, counted with multiplicity. For example, $\Omega(30) = \Omega(-27)=3$.

\begin{lem}\label{lem:elasticityirreducible} For each natural number $f$,
\[ \rho(\Oo_f) =\frac{1}{2} \sup \{\Omega(\Nm_{K/\Q} \pi): \pi\text{ an irreducible element of $\Oo_f$}\}. \]
\end{lem}

\begin{prop}\label{prop:weber} Let $f \in \N$. Every class in $\Cl(\Oo_f)$ is represented by infinitely many prime ideals of $\Oo_K$.
\end{prop}

Our argument for Lemma \ref{lem:elasticitydavenport} also depends on the following simple but useful observation: Suppose $\alpha, \beta \in \Oo_f$ and that $\alpha$ divides $\beta$ in the ring $\Oo_K$. Choose $a \in \Z$ with $\alpha \equiv a\pmod{f\Oo_K}$, and suppose that $\gcd(a,f)=1$. (This coprimality condition is equivalent to the comaximality of $\alpha\Oo_K$ and $f\Oo_K$.) Then $\alpha$ divides $\beta$ in $\Oo_f$. To see this, say that $\beta\equiv b\pmod{f\Oo_K}$. Write $\beta = \alpha\gamma$, where $\gamma \in \Oo_K$. Then $b \equiv a\gamma \pmod{f\Oo_K}$. Multiplying both sides by an integer $\bar{a}$ that inverts $a$ mod $f$ shows that $\gamma \equiv b\bar{a} \pmod{f\Oo_K}$, so that $\gamma \in \Oo_f$. 

\begin{proof}[Proof of Lemma \ref{lem:elasticitydavenport}] For notational convenience, let $D = \Dav \Cl(\Oo_f)$. 

We begin with the lower bound half of Lemma \ref{lem:elasticitydavenport}. Using Proposition \ref{prop:weber}, choose prime ideals $P_1, P_2, \dots, P_D$, comaximal with $f\Oo_K$, such that the product $P_1\cdots P_D$ represents the identity in $\Cl(\Oo_f)$ and  no nonempty, proper subproduct of $P_1, \dots, P_D$ represents the identity. Write $P_1 \cdots P_D = \pi \Oo_K$, where $\pi \in \Oo_f$. Then $\pi$ is irreducible in $\Oo_f$. Otherwise, $\pi = \alpha \beta$ for $\alpha,\beta$ nonunits in $\Oo_f$ and, after rearranging the $P_i$, we have $\alpha \Oo_K = P_1 \cdots P_d$ for some $1\le d < D$. But then $P_1 \cdots P_d$ represents the identity in $\Cl(\Oo_f)$, contrary to the choice of $D$. Invoking Lemma \ref{lem:elasticitydavenport},
\[ \rho(\Oo_f) \ge \frac{1}{2} \Omega(\Nm_{K/\Q} \pi) = \frac{1}{2} \Omega\left(\prod_{i=1}^{D} \Nm(P_i)\right) \ge \frac{1}{2} D. \] 

For the upper bound, we let $\pi$ be an arbitrary irreducible element of $\Oo_f$ and proceed to bound $\Omega(\Nm\pi)$. Write \begin{equation}\label{eq:piidealfact} \pi\Oo_K = \prod_{i=1}^{g} P_i \prod_{j=1}^{h} Q_j^{e_j}, \end{equation}
where the $P_i$ are prime ideals of $\Oo_K$ comaximal with $f\Oo_K$ and $Q_1,\dots, Q_h$ are prime ideals containing $f\Oo_K$. 

Suppose to start with that some $P_i$ has degree $2$. Then $P_i = p_i\Oo_K$ for an inert prime $p_i$ not dividing $f$. Since $p_i$ divides $\pi$ in $\Oo_K$, the observation preceding the proof shows that $p_i$ divides $\pi$ in $\Oo_f$. Hence, $\pi$ is a unit multiple of $p_i$, we have $\Nm_{K/\Q}(\pi) = p_i^2$, and
\begin{equation}\label{eq:normpi1} \frac{1}{2} \Omega(\Nm_{K/\Q}{\pi}) = 1.\end{equation}
For the rest of the proof, we suppose that each $P_i$ has degree $1$.  We proceed to bound $g$ and the exponents $e_1,\dots,e_h$ appearing in \eqref{eq:piidealfact}. 

We must have $g\le D$: If $g > D$, select a subsequence of $P_1,\dots, P_D$ multiplying to $\alpha \Oo_K$ for some $\alpha \in \Oo_f$. Then $\alpha$ divides $\pi$ in $\Oo_K$ and $\alpha\Oo_K$ and $f\Oo_K$ are comaximal. Hence, $\alpha$ divides $\pi$ in $\Oo_f$. Since $\pi$ is irreducible, $\pi$ is a unit multiple of $\alpha$, contradicting that $\frac{\pi}{\alpha}$ is contained in $P_g$.

Now let $Q$ be one of the $Q_j$, and let $e$ be the corresponding $e_j$. Let $q$ be the rational prime lying below $Q$. Since $f$ is split-free, $Q$ is the \emph{unique} prime of $\Oo_K$ lying above $q$. 

We consider first the case that $q$ is inert, so that $q \Oo_K = Q$. If $e \ge v_{q}(f)+1$, we argue that the equation \begin{equation}\label{eq:prospectivefact} \pi = q \cdot \frac{\pi}{q}. \end{equation}
exhibits a nontrivial factorization of $\pi$ over $\Oo_f$, contradicting the irreducibility of $\pi$. Both right-hand factors in \eqref{eq:prospectivefact} are nonunits in $\Oo_K$, so it suffices to show that $\frac{\pi}{q}\in \Oo_f$. Write $f = q^v f'$, where $v=v_q(f)$ and $q\nmid f'$. Since $\pi \in \Oo_f$, there is a rational integer $u$ with $\pi \equiv u\pmod{f\Oo_K}$. Then $q\frac{\pi}{q}=\pi \equiv u\pmod{f'}$. Hence, if $\bar{q}\in \Z$ inverts $q$ mod $f$, we have $\frac{\pi}{q} \equiv u \bar{q}\pmod{f'\Oo_K}$. It follows that $\frac{\pi}{q} \in \Oo_{f'}$. If $e\ge v+1$, we also have $\frac{\pi}{q} \equiv 0\pmod{q^v \Oo_K}$, so that $\frac{\pi}{q} \in \Oo_{q^v}$. Therefore, $\frac{\pi}{q} \in \Oo_{f'} \cap \Oo_{q^v} = \Oo_{f}$. 

If $q$ is ramified, then $q \Oo_K = Q^2$. In this case, an entirely analogous argument to that of the last paragraph shows that \eqref{eq:prospectivefact} contradicts the irreducibility of $\pi$ if $e \ge 2(v_{q}(f)+1)$. 

Collecting our bounds,
\begin{align*} \Omega(N\pi)&= \sum_{i=1}^{g} \Omega(\Nm{P_i}) + \sum_{\substack{1\le j \le h\\ Q_j\text{ inert}}} e_j \cdot \Omega(\Nm{Q_j}) + \sum_{\substack{1\le j \le h\\ Q_j\text{ ramified}}} e_j \cdot \Omega(\Nm{Q_j})\\ 
&\le g + \sum_{\substack{1\le j \le h\\ Q_j\text{ inert}}} 2v_{q_j}(f) + \sum_{\substack{1\le j \le h\\ Q_j\text{ ramified}}} (2v_{q_j}(f)+1) \\
&\le \Dav \Oo_f+ 3 \sum_{1\le j \le h} v_{q_j}(f) \\
&= \Dav \Oo_f + 3\Omega(f).
\end{align*}
(Here and below, $v_{q}$ is the usual $q$-adic valuation.) Hence,
\begin{equation}\label{eq:normpi2} \frac{1}{2}\Omega(\Nm_{K/\Q} \pi) \le \frac{1}{2}\Dav(\Oo_f) + \frac{3}{2}\Omega(f). 
\end{equation}

Lemma \ref{lem:elasticitydavenport} follows from \eqref{eq:normpi1}, \eqref{eq:normpi2}, and Lemma \ref{lem:elasticityirreducible}.
\end{proof}

\subsection{The principal subgroup $\PrinCl(\Oo_f)$ as a proxy for $\Cl(\Oo_f)$} 
For our purposes, rather than directly analyze $\Cl(\Oo_f)$ it is more convenient to work with the group
\[ \PrinCl(\Oo_f):= (\Oo_K/f\Oo_K)^{\times}/\langle \text{images of integers prime to $f$, units of $\Oo_K$}\rangle.\] We can, and will, identify $\PrinCl(\Oo_f)$ with the subgroup $(I_K(f)\cap P_K)/P_{K,\Z}(f)$ of $I_K(f)/P_{K,\Z}(f) = \Cl(\Oo_f)$, with the identification prescribed by the exact sequence
\[ (\Z/f\Z)^{\times} \times \Oo_K^{\times} \stackrel{\mu}{\longrightarrow} (\Oo_K/f\Oo_K)^{\times} \stackrel{\iota}{\longrightarrow} (I_K(f) \cap P_K)/P_{K,\Z}(f) \longrightarrow 1;\]
here $\mu((a\bmod{f\Z},\eta)) := a\eta \bmod{f\Oo_K}$ and $\iota(\alpha\bmod{f\Oo_K}) := [\alpha \Oo_K]$. We think of $\PrinCl(\Oo_f)$ as the ``principal subgroup'' of the class group of $\Oo_f$.

$\PrinCl(\Oo_f)$ fits into the exact sequence
\[ 1 \longrightarrow \conggto{(I_K(f)\cap P_K)/P_{K,\Z}(f)}{\PrinCl(\Oo_f)} \longrightarrow \equalto{I_K(f)/P_{K,\Z}(f)}{\Cl(\Oo_f)} \longrightarrow \equalto{I_K/P_K}{\Cl(\Oo_K)} \longrightarrow 1. \]
(The maps here are the obvious ones. To show exactness at the last position, one needs that every ideal class in $\Oo_K$ has a representative comaximal with $f\Oo_K$; this is immediate from Proposition \ref{prop:weber}, although easier arguments are also possible.) Hence, $\PrinCl(\Oo_f)$ has index $h_K:=\#\Cl(\Oo_K)$ when viewed as a subgroup of $\Cl(\Oo_f)$. Now we are always working with a fixed quadratic field $K$. The following easy lemma will guarantee that, at the level of precision we are aiming for, there is no harm in working with $\PrinCl(\Oo_f)$ in place of $\Cl(\Oo_f)$.

\begin{lem}\label{lem:easydav}  Let $G$ be a finite abelian group, and let $H$ be a subgroup of $G$. Then 
\begin{equation}\label{eq:easydav} \max\{\Dav(H), \Dav(G/H)\} \le \Dav(G) \le \Dav(H) \Dav(G/H). \end{equation}
\end{lem}

As a consequence of Lemma \ref{lem:easydav},
\begin{equation}\label{eq:reductiontoPrinCl} \Dav \PrinCl(\Oo_f) \le \Dav{\Cl(\Oo_f)} \le h_K \Dav \PrinCl(\Oo_f). \end{equation}

\begin{proof}[Proof of Lemma \ref{lem:easydav}] The first inequality in \eqref{eq:easydav} is clear, as a $1$-product-free sequence in $H$ or in $G/H$ corresponds to a $1$-product-free sequence in $G$ of the same length. Turning to the second inequality, we argue that a sequence $g_1, \dots, g_{\Dav(H) \Dav(G/H)}$ in $G$ is never $1$-product-free. Split the sequence into $\Dav(H)$ blocks of $\Dav(G/H)$ terms. Each block contains a nonempty subsequence multiplying to the identity in $G/H$ --- that is, a subsequence whose product belongs to $H$. List these products as $h_1, \dots, h_{\Dav(H)}$. Then some subsequence of $h_1, \dots, h_{\Dav(H)}$ multiplies to $1$. Rewriting each $h_j$ as a product of terms from a certain block of the $g_i$ yields a $1$-product subsequence of $g_1, \dots, g_{\Dav(H) \Dav(G/H)}$.
\end{proof}

\subsection{The pre-class group and two key propositions}\label{sec:preclass} We will study $\PrinCl(\Oo_f)$ by viewing it as a quotient of a yet-simpler object, termed the \emph{pre-class group}, defined by
\[ \PreCl(\Oo_f) = (\Oo_K/f\Oo_K)^{\times}/\langle \text{images of integers prime to $f$}\rangle. \]
Comparing the definitions of $\PrinCl(\Oo_f)$ and $\PreCl(\Oo_f)$, we see that the former is obtained from the latter upon quotienting by the images of units of $\Oo_K$. Our arguments ultimately depend on the groups $\PreCl(\Oo_f)$ being close cousins of the more familiar multiplicative groups $(\Z/m\Z)^{\times}$. 

To set up the analogy, we begin by computing the order of $\PreCl(\Oo_f)$ for split-free $f$. By the Chinese Remainder Theorem, $\PreCl(\Oo_f) \cong \prod_{p^k \parallel f} \PreCl(\Oo_{p^k})$. If $p$ is inert, then 
\[ \#\PreCl(\Oo_{p^k}) = \frac{(\#\Oo_K/(p\Oo_K)^k)^{\times}}{\# (\Z/p^k \Z)^{\times}}= \frac{\Nm((p\Oo_K)^k) - \Nm((p\Oo_K)^{k-1})}{p^k-p^{k-1}} = p^k + p^{k-1}. \]
If $p$ is ramified, with $p\Oo_K = P^2$, then 
\[ \#\PreCl(\Oo_{p^k}) = \frac{(\#\Oo_K/P^{2k})^{\times}}{\# (\Z/p^k \Z)^{\times}}= \frac{\Nm(P^{2k}) - \Nm(P^{2k-1})}{p^k-p^{k-1}} = p^k. \]
Hence, $\#\PreCl(\Oo_f) = \psi(f)$, where $\psi$ is the multiplicative function whose values at prime powers are given by $\psi(p^k) = p^k(1-\frac1p\leg{\Delta}{p})$. This is of course reminiscent of Euler's classical formula $\phi(m) = \prod_{p^k\parallel m} p^k(1-\frac{1}{p})$ for the order of $(\Z/m\Z)^{\times}$.

An important invariant of the group $(\Z/m\Z)^{\times}$ is its exponent, denoted $\lambda(m)$. While typically referred to as \emph{Carmichael's lambda-function}, the study of $\lambda(m)$ goes back to Gauss. In fact, already in the \emph{Disquisitiones}, one can read the result that $\lambda(m)$ is the least common multiple of the numbers $\phi(p^k)$ for the prime powers $p^k\parallel m$, with the caveat that if $p=2$ and $k\ge 3$ one should replace $\phi(2^k)$ with $\frac{1}{2}\phi(2^k)$. 

To prove Theorems \ref{thm:imaginary} and \ref{thm:real}, we will need to understand how the exponents of $\Cl(\Oo_f)$ are distributed. We begin by establishing a (partial) analogue of Gauss's formula. We get this going with the observation that 
$\Exp \PreCl(\Oo_f) = \lcm\{\Exp \PreCl \Oo_{p^k}: p^k \parallel f\}$. The next lemma, due essentially to Halter-Koch, determines almost all of the exponents $\Exp \PreCl(\Oo_{p^k})$.

\begin{lem}\label{lem:HK} If $p > 3$, and $p$ is inert or ramified in $K$, then $\PreCl(\Oo_{p^k})$ is cyclic.
\end{lem}

\begin{proof} Write $K = \Q(\sqrt{D})$, where $D$ is a squarefree integer.

Suppose first that $p$ is inert in $K$. By results of Halter-Koch, summarized in the table on p.\ 77 of \cite{HK72}, $(\Oo_K/p^k \Oo_K)^{\times}$ is generated by (the images of) $w, 1+p$, and $1+p\sqrt{D}$, of respective orders $p^2-1$, $p^{k-1}$, and $p^{k-1}$. Here $w$ (denoted $w_s$ in \cite{HK72}) is an element of $\Oo_K$ comaximal with $p\Oo_K$ whose precise definition does not concern us here (see \cite[p.\ 75]{HK72}). Since $\PreCl(\Oo_{p^k})$ is a quotient of $(\Oo_K/p^k \Oo_K)^{\times}$ in which $1+p$ becomes trivial, $\PreCl(\Oo_{p^k})$ is generated by $w$ and $1+p\sqrt{D}$, of respective orders dividing $p^2-1$ and $p^{k-1}$. But $p^2-1$ and $p^{k-1}$ are relatively prime. Hence, $w(1+p\sqrt{D})$ generates $\PreCl(\Oo_{p^k})$. 

The ramified case is similar. Here  $(\Oo_K/p^k\Oo_K)^{\times}$ is generated by three elements $w$, $1+p$, and $1+\sqrt{D}$ (same table in \cite{HK72}). In this case, $w$ can be chosen as a rational integer that generates $(\Z/p\Z)^{\times}$. So both $w$ and $1+p$ become trivial in the quotient defining $\PreCl(\Oo_{p^k})$, implying that $\PreCl(\Oo_{p^k})$ is generated by the image of $1+\sqrt{D}$.
\end{proof}

From Lemma \ref{lem:HK} and the preceding remarks, if we define $L(f):= \lcm\{\psi(p^k): p^k \parallel f\}$ and $L'(f):= \lcm\{\psi(p^k): p^k \parallel f, \,p> 3\}$, then
\begin{equation}\label{LLprimeDivides} L'(f) \mid \Exp \PreCl(\Oo_f) \mid L(f) \end{equation}
for all split-free numbers $f$. The typical behavior of $\Exp \PreCl(\Oo_f)$ for split-free $f$ is now determined by our next result, which is the analogue of a theorem proved for Carmichael's $\lambda$-function by Erd\H{o}s, Pomerance, and Schmutz \cite[Theorem 2]{EPS}.

\begin{prop}\label{prop:Lforder} For almost all split-free $f$,
\[ L(f) = f/(\log{f})^{\frac{1}{2}\log_3{f} + \frac{1}{2}C_K+O((\log_4f)^3/\log_3f)}.\]
The same estimate holds with $L'(f)$ replacing $L(f)$.
\end{prop}

As we explain in the next subsection, our Theorem \ref{thm:imaginary} (typical elasticity for imaginary quadratic orders) can be quickly deduced from Proposition \ref{prop:Lforder}. To prove Theorem \ref{thm:real} (concerning real quadratic orders) we must work a bit harder. 

Recall that when $K$ is a real quadratic field, we are using $\epsilon$ to denote the fundamental unit of $K$, normalized so that $\epsilon > 1$. Then $\PrinCl(\Oo_f) \cong \PreCl(\Oo_f)/\langle \text{image of $\epsilon$}\rangle$. We let $\ell(f)$ denote the order of $\epsilon$ viewed inside $\PreCl(\Oo_f)$, so that 
\[ \ell(f) = \frac{\#\PreCl(\Oo_f)}{\#\PrinCl(\Oo_f)}.\]
Concretely, $\ell(f)$ can be described as the least positive integer $\ell$ for which $\epsilon^{\ell} \in \Oo_f$.

The next proposition collects various estimates needed in our proof of Theorem \ref{thm:real}.

\begin{prop}[conditional on GRH]\label{prop:keythm2} Let $K$ be a real quadratic field. For almost all split-free $f$,
\begin{enumerate}[label={\rm (\alph*)}]
\item $\RadExp \PrinCl(\Oo_f)= (\log f)^{\frac12 + O(1/\log_4 f)}$,
\item $\Exp \PrinCl(\Oo_f) = (\RadExp \PrinCl(\Oo_f)) \cdot (\log f)^{O(1/\log_4 f)}$,
\item $L(f)/\ell(f) = (\log f)^{O(1/\log_4 f)}$.
\end{enumerate}
\end{prop}
As with Proposition \ref{prop:Lforder}, the estimates of Proposition \ref{prop:keythm2} take their inspiration from research into the multiplicative groups $(\Z/m\Z)^{\times}$. Specifically, we adapt ideas introduced by Pollack
in \cite{pollack21} to show that $\frac{\phi(m)}{\lambda(m)}$ typically has about $\frac{\log_2 m}{\log_3 m}$ distinct prime factors. We also draw heavily on methods introduced by Li and Pomerance \cite{LP03}, who do a statistical comparison of $\lambda(m)$ with $\ell_a(m)$, the order of $a$ (a fixed integer) in $(\Z/m\Z)^{\times}$. 

\subsection{Proofs of Theorems \ref{thm:imaginary} and \ref{thm:real}, modulo Propositions \ref{prop:Lforder} and \ref{prop:keythm2}} We now describe how to complete the proofs of Theorems \ref{thm:imaginary} and \ref{thm:real}, taking for granted Propositions \ref{prop:Lforder} and \ref{prop:keythm2}. In addition to the tools already introduced, we rely on known results relating the Davenport constant and the exponent.

\begin{prop}\label{prop:davexp} For every finite abelian group $G$,
\[ 1\le \frac{\Dav G}{\Exp G} \le 1 + \log \frac{\#G}{\Exp G}. \]
\end{prop}

The lower bound in Proposition \ref{prop:davexp} is trivial: If $g$ has order $\Exp G$, then $g, g, \dots, g$ (repeated $(\Exp G)-1$ times) is $1$-product-free. The nontrivial and elegant upper bound is due to van Emde Boas and Kruyswijk \cite{EBK69}. A simplified proof can be found in \cite{AGP94}.

Now let $K$ be an imaginary quadratic field. The image $H$ (say) of $\Oo_K^{\times}$ in $\PreCl(\Oo_f)$ has size at most $\#\Oo_K^{\times} \le 6$. Identifying $\PrinCl(\Oo_f)$ with $\PreCl(\Oo_f)/H$, Lemma \ref{lem:easydav} implies that 
\[ \frac{\Dav \PreCl(\Oo_f)}{\Dav(H)} \le \Dav \PrinCl(\Oo_f) \le \Dav \PreCl(\Oo_f). \]
It is well-known and simple-to-show that the Davenport constant of a group is bounded by the size of the group. Hence, $\Dav H \le 6$, and $\Dav \PrinCl(\Oo_f)$ is within a factor of $6$ of $\Dav \PreCl(\Oo_f)$. 

Applying Proposition \ref{prop:davexp}, 
\begin{equation}\label{eq:imagdeduction} 1 \le \frac{\Dav \PreCl(\Oo_f)}{\Exp \PreCl(\Oo_f)} \le 1+ \log \frac{\psi(f)}{\Exp \PreCl(\Oo_f)}.\end{equation}
Notice that $\psi(f) \le f \prod_{p\mid f}(1+\frac{1}{p}) \ll f \log\log{(3f)}$, for all split-free $f$. Furthermore, from \eqref{LLprimeDivides} and Proposition \ref{prop:Lforder}, we typically have
\begin{equation}\label{eq:expprecltypical} \Exp \PreCl(\Oo_f) = f/(\log{f})^{\frac{1}{2}\log_3{f} + \frac{1}{2}C_K+O((\log_4f)^3/\log_3f))} \end{equation} Hence, $\psi(f)/\Exp \PreCl(\Oo_f)$ is typically at most $(\log{f})^{\log_3{f}}$, and $1 + \log(\psi(f)/\Exp \PreCl(\Oo_f)) \ll \log_2{f} \log_3{f}$. Referring back to \eqref{eq:imagdeduction},
\[ \Dav \PreCl(\Oo_f) = (\Exp \PreCl(\Oo_f)) \exp(O(\log_3{f})), \]
typically. Plugging in the estimate \eqref{eq:expprecltypical},
\[ \Dav \PreCl(\Oo_f) = f/(\log{f})^{\frac{1}{2}\log_3{f} + \frac{1}{2}C_K+O((\log_4f)^3/\log_3f))}, \]
for almost all split-free $f$. Our discussion in the last paragraph permits replacing $\PreCl$ with $\PrinCl$, and  \eqref{eq:reductiontoPrinCl} allows us to replace $\PrinCl$ with $\Cl$.  Theorem \ref{thm:imaginary} then follows from Lemma \ref{lem:elasticitydavenport}, in view of the trivial upper bound $\Omega(f) \le \frac{\log{f}}{\log{2}}$.

Turning to Theorem \ref{thm:real}, suppose that $K$ is real quadratic. By Proposition \ref{prop:keythm2}(a,b), 
\begin{equation}\label{eq:expprinclreal} \Exp \PrinCl(\Oo_f) = (\log{f})^{\frac{1}{2} + O(1/\log_4{f})}\end{equation}
for almost all split-free $f$. Furthermore, 
\[ \frac{\#\PrinCl(\Oo_f)}{\Exp \PrinCl(\Oo_f)} \le \#\PrinCl(\Oo_f) = \frac{\psi(f)}{\ell(f)}= \frac{\psi(f)}{L(f)} \frac{L(f)}{\ell(f)} \le (\log{f})^{\log_3{f}}, \]
typically. For the final inequality we used Proposition \ref{prop:Lforder} and Proposition \ref{prop:keythm2}(c). Invoking Proposition \ref{prop:davexp},
\[ \Dav \PrinCl(\Oo_f) = (\Exp \PrinCl(\Oo_f)) \exp(O(\log_3{f})), \]
almost always.
Substituting in the estimate of \eqref{eq:expprinclreal}, we see that typically
\[ \Dav \PrinCl(\Oo_f) = (\log{f})^{\frac{1}{2} + O(1/\log_4{f})}. \]
Owing to \eqref{eq:reductiontoPrinCl}, the same estimate holds with $\PrinCl$ replaced by $\Cl$. Theorem \ref{thm:real} now follows from Lemma \ref{lem:elasticitydavenport}, since the count of all positive integers $f\le x$ with $\Omega(f) > 10\log\log{f}$ (say) is $\ll x/\log{x} = o(S_{\alpha}(x))$. (This last estimate follows from well-known results on the distribution of numbers with many prime factors; see for instance \cite[Exercise 08]{HT88} or  \cite[Lemma 13]{LP07}. To apply the estimates as stated there, treat separately the cases when $f\le \sqrt{x}$ and $\sqrt{x} < f \le x$.) 

\subsection*{What lies ahead} The rest of the paper is organized as follows. In \S\ref{sec:weightedTK}, we establish a variant of the Tur\'{a}n--Kubilius inequality result for additive functions appearing with multiplicative weights. In \S\ref{sec:Lforder} we prove the key Proposition \ref{prop:Lforder}. There our weighted Tur\'{a}n--Kubilius inequality, with weight function $1_{\text{split-free}}$, plays an important role. Section \ref{sec:radicalpsiL} is something of a waypoint. There we show that $\Rad \frac{\psi(f)}{L(f)}$ is typically of size $\approx (\log{f})^{1/2}$. It will turn out that $\Rad \frac{\psi(f)}{L(f)}$ is a reasonable approximation to $\Exp \PreCl(\Oo_f)$ (typically, and assuming GRH). To connect those two quantities requires us to relate $L(f)$ and $\ell(f)$, for typical $f$. In \S\ref{sec:algebraicnonsense} we set up the algebraic framework needed to carry out this comparison. In \S\ref{sec:completionkeythm2} we present the details; this work, supplemented by various `anatomical' arguments, allows us to complete the proof of the key Proposition \ref{prop:keythm2}.

\section{A weighted Tur\'{a}n--Kubilius inequality}\label{sec:weightedTK}
An important tool for studying the normal order of an additive function is the Tur\'{a}n--Kubilius inequality. Let $f\colon\N\to\C$ be an additive function. It is often reasonable to think of the mean value of $f$ as a good candidate for its normal order. Since the mean value of $f$ over $[1,x]$ is
\[\frac{1}{x}\sum_{n\le x}f(n)=\frac{1}{x}\sum_{p^k\le x}f(p^k)\sum_{\substack{n\le x\\ p^k\parallel n}}1=\frac{1}{x}\sum_{p^k\le x}f(p^k)\left(\left\lfloor\frac{x}{p^k}\right\rfloor-\left\lfloor\frac{x}{p^{k+1}}\right\rfloor\right)=A_f(x)+O\left(\frac{1}{x}\sum_{p^k\le x}|f(p^k)|\right),\]
where
\[A_f(x)\colonequals\sum_{p^k\le x}\frac{f(p^k)}{p^k}\left(1-\frac{1}{p}\right),\]
in many situations we may think of $A_f(x)$ as an approximation to the mean value of $f$ over $[1,x]$. The Tur\'{a}n--Kubilius inequality furnishes an upper bound for the mean square of $f(n)-A_f(x)$, which may be thought of as the ``variance" of $f$ over $[1,x]$. In its simplest form, the Tur\'{a}n--Kubilius inequality asserts that
\[\frac{1}{x}\sum_{n\le x}|f(n)-A_f(x)|^2\ll B_f(x),\]
where 
\[B_f(x)\colonequals\sum_{p^k\le x}\frac{|f(p^k)|^2}{p^k}.\]
An immediate corollary of this inequality is that if $B_f(x)=o(|A_f(x)|^2)$ as $x\to\infty$, then for every fixed $\epsilon>0$, we have $|f(n)-A_f(x)|<\epsilon|A_f(x)|$ for all but $o(x)$ values of $n\in\N\cap[1,x]$. Moreover, if one can show that $A_f(n)$ is close to $A_f(x)$ on average, then $A_f(n)$ serves as a normal order of $f(n)$.

For our applications, we will need a version of the Tur\'{a}n--Kubilius inequality in which all the means involved are taken over the set of split-free positive integers. From a probabilistic point of view, we think of $n\in\N\cap[1,x]$ as a discrete random variable with probability distribution given by $\mathbb{P}(n=m)=1_{\text{split-free}}(m)/\sum_{k\le x}1_{\text{split-free}}(k)$ for all $m\in\N\cap[1,x]$, rather than as a discrete random variable with uniform distribution given by $\mathbb{P}(n=m)=1/\lfloor x\rfloor$ for all $m\in\N\cap[1,x]$. More generally, one may replace $1_{\text{split-free}}$ by an arbitrary nonnegative multiplicative function $\alpha$. In this direction, we prove the following weighted version of the Tur\'{a}n--Kubilius inequality inspired by Remark 3.11.1 in 
Fan's PhD thesis on weighted Erd\H{o}s--Kac theorems \cite{SFTh23}.

\begin{thm}\label{thm:WTK}
Let $f\colon\N\to\C$ be an additive function, and let  $\alpha\colon\N\to\R_{\ge0}$ be a multiplicative function with partial sums
\[S_{\alpha}(x)\colonequals\sum_{n\le x}\alpha(n).\]
Suppose that there exist $c_{\alpha},\delta>0$, $\sigma\ge0$ and $\kappa\in\R$, such that
\begin{equation}\label{eq:S_alpha}
\sum_{\substack{n\le x\\(n,a)=1}}\alpha(n)=c_{\alpha}x^{\sigma}(\log 3x)^{\kappa-1}\left(F_{\alpha}(a)^{-1}+O\left(\frac{1}{\log\log 3x}\right)\right)    
\end{equation}
uniformly for all $x\ge1$ and all squarefree $a\in\N\cap[1,x^{\delta}]$ with at most two prime factors, where
\[F_{\alpha}(a)\colonequals\prod_{p\mid a}\sum_{k\ge0}\frac{\alpha(p^k)}{p^{k\sigma}}<\infty.\]
Furthermore, suppose that
\begin{equation}\label{eq:Chebyshevalpha}
\sum_{p^k\le x}\frac{\alpha(p^k)}{p^{k\sigma}}\log p^k\ll \log x
\end{equation}
for all $x\ge2$. Then we have
\begin{equation}\label{eq:WEK}
S_{\alpha}(x)^{-1}\sum_{n\le x}\alpha(n)\left|f(n)-A_{\alpha,f}(x)\right|^2\ll B_{\alpha,f}(x)    
\end{equation}
for all $x\ge1$, where
\begin{align*}
A_{\alpha,f}(x)&\colonequals\sum_{p^k\le x}\alpha(p^k)F_{\alpha}(p)^{-1}\frac{f(p^k)}{p^{k\sigma}},\\
B_{\alpha,f}(x)&\colonequals\sum_{p^k\le x}\alpha(p^k)\frac{|f(p^k)|^2}{p^{k\sigma}}\left(1-\frac{\log p^k}{\log 3x}\right)^{\min(\kappa-1,0)}.
\end{align*}
The implied constant in \eqref{eq:WEK} depends at most on the parameters $\delta,\kappa$ and the implied constants in \eqref{eq:S_alpha} and \eqref{eq:Chebyshevalpha}.
\end{thm}
\begin{proof}
We may suppose $x\ge2$, since \eqref{eq:WEK} holds trivially when $x\in[1,2)$. For technical reasons, we first prove \eqref{eq:WEK} for the additive function $f_{\eta}$ defined by $f_{\eta}(p^k)=f(p^k)1_{p^k\le x^{\eta}}$, where $\eta=\delta/(2+2\delta)$. For convenience, we adopt the shorthand notation $1_{p^k}(n)\colonequals 1_{p^k\parallel n}(n)$ and the notation
\[\E_{\le x}^{\alpha}[g]\colonequals S_{\alpha}(x)^{-1}\sum_{n\le x}\alpha(n)g(n)\]
for any arithmetic function $g\colon\N\to\C$. We start by computing the expectation of $f_{\eta}$. For $p^k\le x^{\eta}$, we have by \eqref{eq:S_alpha} that
\begin{align*}
 \E_{\le x}^{\alpha}[1_{p^k}]&=S_{\alpha}(x)^{-1}\sum_{\substack{n\le x\\p^k\parallel n}}\alpha(n)\\
 &=S_{\alpha}(x)^{-1}\alpha(p^k)\sum_{\substack{n\le x/p^k\\(n,p)=1}}\alpha(n)\\
 &=\frac{\alpha(p^k)}{p^{k\sigma}}\left(1-\frac{\log p^k}{\log x}\right)^{\kappa-1}\left(F_{\alpha}(p)^{-1}+O\left(\frac{1}{\log\log 3x}\right)\right)\\
 &=\frac{\alpha(p^k)}{p^{k\sigma}}\left(F_{\alpha}(p)^{-1}+O\left(\frac{\log p^k}{\log x}+\frac{1}{\log\log 3x}\right)\right).
\end{align*}
Hence, the expectation of $f_{\eta}$ is
\[\E_{\le x}^{\alpha}[f_{\eta}]=\sum_{p^k\le x^{\eta}}f(p^k)\E_{\le x}^{\alpha}[1_{p^k}]=A_{\alpha,f_{\eta}}(x)+O\left(\sum_{p^k\le x^{\eta}}\alpha(p^k)\frac{|f(p^k)|}{p^{k\sigma}}\left(\frac{\log p^k}{\log x}+\frac{1}{\log\log 3x}\right)\right).\]
By Cauchy--Schwarz and \eqref{eq:Chebyshevalpha}, the error term on the previous line is
\[\ll \left(\sum_{p^k\le x^{\eta}}\alpha(p^k)\frac{|f(p^k)|^2}{p^{k\sigma}}\right)^{1/2}\left(\sum_{p^k\le x^{\eta}}\frac{\alpha(p^k)}{p^{k\sigma}}\left(\left(\frac{\log p^k}{\log x}\right)^2+\frac{1}{(\log\log 3x)^2}\right)\right)^{1/2}\ll \sqrt{B_{\alpha,f_{\eta}}(x)},\]
whence
\begin{equation}\label{eq:E[f_eta]}
 \E_{\le x}^{\alpha}[f_{\eta}]=A_{\alpha,f_{\eta}}(x)+O\left(\sqrt{B_{\alpha,f_{\eta}}(x)}\right).   
\end{equation}
Next, we compute the variance of $f_{\eta}$. It is clear that
\[\Var_{\le x}^{\alpha}[1_{p^k}]\colonequals \E_{\le x}^{\alpha}\left[\left|1_{p^k}-\E_{\le x}^{\alpha}[1_{p^k}]\right|^2\right]=\E_{\le x}^{\alpha}[1_{p^k}^2]-\left(\E_{\le x}^{\alpha}[1_{p^k}]\right)^2\ll\frac{\alpha(p^k)}{p^{k\sigma}}.\]
For $p^k,q^{l}\le x^{\eta}$ with $p\ne q$, we have by \eqref{eq:S_alpha} that
\begin{align*}
\E_{\le x}^{\alpha}[1_{p^k}1_{q^l}]&=S_{\alpha}(x)^{-1}\alpha(p^k)\alpha(q^l)\sum_{\substack{n\le x/p^kq^l\\(n,pq)=1}}\alpha(n)\\ 
&=\frac{\alpha(p^k)\alpha(q^l)}{p^{k\sigma}q^{l\sigma}}\left(1-\frac{\log p^kq^l}{\log x}\right)^{\kappa-1}\left(F_{\alpha}(pq)^{-1}+O\left(\frac{1}{\log\log 3x}\right)\right).
\end{align*}
It follows that
\begin{align*}
\Cov_{\le x}^{\alpha}[1_{p^k},1_{q^l}]&\colonequals \E_{\le x}^{\alpha}[1_{p^k}1_{q^l}]-\E_{\le x}^{\alpha}[1_{p^k}]\E_{\le x}^{\alpha}[1_{q^l}]\\ 
&=\frac{\alpha(p^k)\alpha(q^l)}{p^{k\sigma}q^{l\sigma}}\left(1-\frac{\log p^kq^l}{\log x}\right)^{\kappa-1}\left(F_{\alpha}(pq)^{-1}+O\left(\frac{1}{\log\log 3x}\right)\right)\\
&\quad-\frac{\alpha(p^k)\alpha(q^l)}{p^{k\sigma}q^{l\sigma}}\left(1-\frac{\log p^k}{\log x}\right)^{\kappa-1}\left(1-\frac{\log q^l}{\log x}\right)^{\kappa-1}\left(F_{\alpha}(pq)^{-1}+O\left(\frac{1}{\log\log 3x}\right)\right)\\
&\ll\frac{\alpha(p^k)\alpha(q^l)}{p^{k\sigma}q^{l\sigma}}\left(\frac{\log p^k}{\log x}\cdot\frac{\log q^l}{\log x}+\frac{1}{\log\log 3x}\right),
\end{align*}
where we have used the inequality $|a^{\kappa-1}-b^{\kappa-1}|\ll |a-b|$ for any $a,b\in[(1-\eta)^2,1]$, which is a direct consequence of the mean value theorem in calculus. Thus, the variance of $f_{\eta}$ is
\begin{align*}
\Var_{\le x}^{\alpha}[f_{\eta}]=\E_{\le x}^{\alpha}\left[\left|f_{\eta}-\E_{\le x}^{\alpha}[f_{\eta}]\right|^2\right]&=\sum_{p^k\le x^{\eta}}|f(p^k)|^2\Var_{\le x}^{\alpha}[1_{p^k}]+\sum_{\substack{p^k,q^l\le x^{\eta}\\p\ne q}}f(p^k)\overline{f(q^l)}\Cov_{\le x}^{\alpha}[1_{p^k},1_{q^l}]\\
&\ll B_{\alpha,f_{\eta}}(x)+\left(\sum_{p^k\le x^{\eta}}\alpha(p^k)\frac{|f(p^k)|}{p^{k\sigma}}\left(\frac{\log p^k}{\log x}+\frac{1}{\sqrt{\log\log 3x}}\right)\right)^2\\
&\ll B_{\alpha,f_{\eta}}(x)
\end{align*}
by Cauchy--Schwarz and \eqref{eq:Chebyshevalpha} as before. Combining this with \eqref{eq:E[f_eta]}, we find that
\[\E_{\le x}^{\alpha}\left[|f_{\eta}-A_{\alpha,f_{\eta}}(x)|^2\right]\le 2\Var_{\le x}^{\alpha}[f_{\eta}]+2\left|\E_{\le x}^{\alpha}[f_{\eta}]-A_{\alpha,f_{\eta}}(x)\right|^2\ll B_{\alpha,f_{\eta}}(x),\]
which is exactly \eqref{eq:WEK} with $f_{\eta}$ in place of $f$.

Now it is an easy matter to deduce \eqref{eq:WEK} for a general additive function $f$. Since
\[|f-A_{\alpha,f}(x)|^2\le 3\left(|f-f_{\eta}|^2+|f_{\eta}-A_{\alpha,f_{\eta}}(x)|^2+|A_{\alpha,f_{\eta}}(x)-A_{\alpha,f}(x)|^2\right),\]
applying $\E_{\le x}^{\alpha}$ to both sides yields
\[\E_{\le x}^{\alpha}\left[|f-A_{\alpha,f}(x)|^2\right]\ll \E_{\le x}^{\alpha}\left[|f-f_{\eta}|^2\right]+B_{\alpha,f_{\eta}}(x)+|A_{\alpha,f_{\eta}}(x)-A_{\alpha,f}(x)|^2.\]
Clearly, $B_{\alpha,f_{\eta}}(x)\ll B_{\alpha,f}(x)$. To estimate $\E_{\le x}^{\alpha}\left[|f-f_{\eta}|^2\right]$, we observe that for any $n\in\N\cap[1,x]$, the number of prime powers $p^k>x^{\eta}$ exactly dividing $n$ is at most $1/\eta$. Hence, we have by \eqref{eq:S_alpha} that
\begin{align*}
\E_{\le x}^{\alpha}\left[|f-f_{\eta}|^2\right]&=S_{\alpha}(x)^{-1}\sum_{n\le x}\alpha(n)\left|\sum_{\substack{p^k\parallel n,\,p^k>x^{\eta}}}f(p^k)\right|^2\\
&\ll S_{\alpha}(x)^{-1}\sum_{n\le x}\alpha(n)\sum_{\substack{p^k\parallel n,\,p^k>x^{\eta}}}|f(p^k)|^2\\
&\le S_{\alpha}(x)^{-1}\sum_{p^k\le x}|f(p^k)|^2\sum_{\substack{n\le x\\ p^k\parallel n}}\alpha(n)\\
&\le S_{\alpha}(x)^{-1}\sum_{p^k\le x}\alpha(p^k)|f(p^k)|^2\sum_{n\le x/p^k}\alpha(n)\\
&\ll B_{\alpha,f}(x).
\end{align*}
Finally, Cauchy--Schwarz together with \eqref{eq:Chebyshevalpha} yields the upper bound
\begin{align*}
|A_{\alpha,f_{\eta}}(x)-A_{\alpha,f}(x)|^2&=\left|\sum_{x^{\eta}<p^k\le x}\alpha(p^k)F_{\alpha}(p)^{-1}\frac{f(p^k)}{p^{k\sigma}}\right|^2\\
&\ll B_{\alpha,f}(x)\sum_{x^{\eta}<p^k\le x}\frac{\alpha(p^k)}{p^{k\sigma}}\\
&\ll \frac{B_{\alpha,f}(x)}{\log x}\sum_{p^k\le x}\frac{\alpha(p^k)}{p^{k\sigma}}\log p^k\\
&\ll B_{\alpha,f}(x).
\end{align*}
Putting everything together, we obtain $\E_{\le x}^{\alpha}\left[|f-A_{\alpha,f}(x)|^2\right]\ll B_{\alpha,f}(x)$ as desired.
\end{proof}

Examining the proof, one sees readily that Theorem \ref{thm:WTK} remains valid if $F_{\alpha}$ is replaced by any positive multiplicative function whose restriction on primes is bounded away from 0. Nevertheless, the definition of $F_{\alpha}$ given in Theorem \ref{thm:WTK} is intuitive and follows heuristically from \eqref{eq:S_alpha}. With this definition, \eqref{eq:S_alpha} holds true for a myriad of multiplicative functions $\alpha$; see \cite[Lemma 3.2]{SF23}, for instance. In the sequel, we shall always take $\alpha$ to be the characteristic function of the set of split-free integers. For this particular choice of $\alpha$ we have 
\begin{equation}\label{eq:sfcount}
 S_{\alpha}(x)=c_{\alpha}x(\log 3x)^{-1/2}\left(1+O\left(\frac{1}{\log3x}\right)\right)   
\end{equation}
and
\[\sum_{p^k\le x}\frac{\alpha(p^k)}{p^{k\sigma}}\log p^k=\frac{1}{2}\log x+O(1)\]
for all $x\ge1$, where
\begin{align*}
c_{\alpha}&=\frac{1}{\sqrt{\pi}}\prod_{p}\left(1-\frac{1}{p}\right)^{1/2}\sum_{k\ge0}\frac{\alpha(p^k)}{p^{k}}\\
&=\sqrt{\frac{1}{\pi L(1,\chi)}\cdot\frac{|\Delta|}{\phi(|\Delta|)}}\prod_{p\text{~inert}}\left(1-\frac{1}{p^2}\right)^{-1/2},
\end{align*} 
with $\chi\colonequals(\Delta/\cdot)$. We remind the reader that the value $L(1,\chi)$ can be expressed in terms of the arithmetic invariants of $K$ via 
Dirichlet's class number formula \cite[Ch. 26]{Rib01}:
\[h_K=\begin{cases}
 \displaystyle\frac{\sqrt{\Delta}}{2\log\epsilon}L(1,\chi),\quad&\text{if~}K\text{~is real},\\[10pt]
 \displaystyle\frac{w\sqrt{-\Delta}}{2\pi}L(1,\chi),\quad&\text{if~}K\text{~is imaginary},
\end{cases}\]
where $\epsilon>1$ is the normalized fundamental unit of $K$, and
\[w=\begin{cases}
 2,\quad&\text{if~}\Delta<-4,\\
 4,\quad&\text{if~}\Delta=-4,\\
 6,\quad&\text{if~}\Delta=-3,
\end{cases}\]
is the number of roots of unity in $K$.
Moreover, it can be shown \cite[Lemma 3.2]{SF23} that 
\[\sum_{\substack{n\le x\\(n,a)=1}}\alpha(n)=c_{\alpha}x(\log 3x)^{-1/2}\left(F_{\alpha}(a)^{-1}+O\left(\frac{1}{\log3x}\right)\right)\]
uniformly for all $x\ge1$ and all $a\in\N\cap[1,x]$, where $F_{\alpha}(a)$ is given by
\[F_{\alpha}(a)=\prod_{p\mid a}\sum_{k\ge0}\frac{\alpha(p^k)}{p^{k}}=\prod_{\substack{p\mid a\\\chi(p)\ne1}}\left(1-\frac{1}{p}\right)^{-1}.\]
Hence, $\alpha$ satisfies the conditions of Theorem \ref{thm:WTK}. We will apply our weighted Tur\'{a}n--Kubilius inequality to prove Proposition \ref{prop:Lforder}.

\section{The typical size of $L(f)$: Proof of Proposition \ref{prop:Lforder}}\label{sec:Lforder}
Our proof of Proposition \ref{prop:Lforder} builds on that of \cite[Theorem 2]{EPS} on the normal order of Carmichael's $\lambda$-function. Due to the introduction of the weight $\alpha$ used to capture only the split-free integers $f$, some special care needs to be taken of the estimation of various weighted sums in our proof. For this reason, our argument is more delicate than its counterpart in \cite{EPS}. We shall only prove Proposition \ref{prop:Lforder} for $L(f)$, as one sees readily from the proof that primes in any bounded interval which divide a typical $f$ contribute a negligible amount and can thus be left out. 

To begin with, we write
\begin{align*}
\log\psi(f)&=\sum_{q}v_q(\psi(f))\log q,\\
\log L(f)&=\sum_{q}v_q(L(f))\log q,
\end{align*}
where the sums run over all primes $q$.
Since $f/\log_2 f\ll \psi(f)\ll f\log_2 f$ for all $f\ge3$, it suffices to show 
\begin{equation}\label{eq:psi-L}
 \log\psi(f)-\log L(f)=\frac{1}{2}y\log y +\frac{1}{2}C_Ky +O\left(\frac{y(\log_2y)^3}{\log y}\right)
\end{equation}
for all but $o(x/\sqrt{\log x})$ split-free $f\le x$, where $y=\log_2x$. As in \cite{EPS}, we divide the primes $q$ into the following four ranges:
\begin{alignat*}{2}
&I_1\colon q\le y/\log y, \quad\quad\quad &&I_2\colon y/\log y<q\le y\log y,\\
&I_3\colon y\log y<q\le y^2,\quad\quad\quad &&I_4\colon q>y^2.  
\end{alignat*}
We estimate the contribution to the left-hand side of \eqref{eq:psi-L} from the primes in each of these four intervals separately. In the proof, we shall frequently resort to \eqref{eq:sfcount} without further notice. Note that $q^2\nmid f$ for any $q\notin I_1$ for all but $o(S_{\alpha}(x))=o(x/\sqrt{\log x})$ split-free $f\le x$. This follows from the estimate
\begin{align*}
S_{\alpha}(x)^{-1}\sum_{q\notin I_1}\sum_{\substack{f\le x\\ q^2\mid f}}\alpha(f)&=S_{\alpha}(x)^{-1}\sum_{q>y/\log y}\alpha(q^2)\sum_{f\le x/q^2}\alpha(f)\\
&\ll\sum_{y/\log y<q\le\sqrt{x}}\frac{1}{q^2}\left(1-\frac{\log q^2}{\log 3x}\right)^{-1/2}\\
&\ll\sum_{y/\log y<q\le x^{1/3}}\frac{1}{q^2}+(\log x)^{1/2}\sum_{x^{1/3}<q\le x}\frac{1}{q^2}\ll \frac{1}{y}.
\end{align*}
We will make use of this fact when examining the contributions from $q\in I_3\cup I_4$. Finally, we set $\Rr_i\colonequals\{a\in\Z/\Delta\Z\colon\chi(a)=i\}$ and $\Pp_i\colonequals\{p\text{~prime}\colon\chi(p)=i\}$ for $i\in\{0,\pm1\}$.

\subsection{A cutoff of $\log\psi$}\label{subsec:logpsi}
Before estimating the contribution from primes in $I_1\cup I_2$, we define a cutoff of $\log\psi$ by
\[h(f)\colonequals\sum_{q\le y\log y}v_q(\psi(f))\log q.\]
We wish to determine the typical size of the additive function $h$ by applying Theorem \ref{thm:WTK} to $h$ and the multiplicative weight $\alpha$. It boils down to estimating
\begin{align*}
A_{\alpha,h}(x)&=\sum_{p^k\le x}\alpha(p^k)\frac{h(p^k)}{p^{k}}\left(1-\frac{1}{p}\right),\\
B_{\alpha,h}(x)&=\sum_{p^k\le x}\alpha(p^k)\frac{h(p^k)^2}{p^{k}}\left(1-\frac{\log p^k}{\log 3x}\right)^{-1/2}.
\end{align*}
Using the trivial bound $h(p^k)\le\log\psi(p^k)\ll \log p^k$, we find that
\begin{align*}
\sum_{p^k\le x}\alpha(p^k)\frac{h(p^k)}{p^{k+1}}\ll\sum_{p^k\le x}\frac{\log p^k}{p^{k+1}}\ll\sum_{p\le x}\frac{\log p}{p^2}\ll1
\end{align*}
and that
\begin{align*}
\sum_{\substack{p^k\le x\\k\ge2}}\alpha(p^k)\frac{h(p^k)^j}{p^{k}}\left(1-\frac{\log p^k}{\log 3x}\right)^{-1/2}&\ll\sum_{\substack{p^k\le \sqrt{x}\\k\ge2}}\frac{(\log p^k)^j}{p^{k}}+(\log x)^{j+1/2}\sum_{\substack{\sqrt{x}<p^k\le x\\k\ge2}}\frac{1}{p^{k}}\\
&\ll1+(\log x)^{j+1/2}\sum_{\substack{\sqrt{x}<p^k\le x\\p\le x^{1/4}}}\frac{1}{p^{k}}+(\log x)^{j+1/2}\sum_{\substack{p>x^{1/4}\\k\ge2}}\frac{1}{p^{k}}\\
&\ll 1+\frac{(\log x)^{j+1/2}}{x^{1/2}}\sum_{p\le x^{1/4}}1+(\log x)^{j+1/2}\sum_{p>x^{1/4}}\frac{1}{p^2}\ll1,
\end{align*}
where $j\in\{1,2\}$. 
It follows that
\begin{align*}
A_{\alpha,h}(x)&=\sum_{p\le x}\alpha(p)\frac{h(p)}{p}+O(1),\\
B_{\alpha,h}(x)&=\sum_{p\le x}\alpha(p)\frac{h(p)^2}{p}\left(1-\frac{\log p}{\log 3x}\right)^{-1/2}+O(1).
\end{align*}
Let us estimate $A_{\alpha,h}(x)$ first. By definition of $h$, we have
\begin{align*}
 A_{\alpha,h}(x)&=\sum_{q\le y\log y}\log q\sum_{\substack{p\le x\\p\notin\Pp_1}}\frac{v_q(p-\chi(p))}{p}+O(1)\\
 &=\sum_{q\le y\log y}\log q\sum_{i\ge1}\sum_{\substack{p\le x\\p\notin\Pp_1\\p\equiv\chi(p)\pmod*{q^i}}}\frac{1}{p}+O(1)\\
 &=\sum_{q\le y\log y}\log q\sum_{i\ge1}\sum_{\substack{p\le x\\p\in\Pp_{-1}\\p\equiv-1\pmod*{q^i}}}\frac{1}{p}+O(1).
\end{align*}
For each $n\in\N$, let $\beta(n)\colonequals\#\{a\in\Rr_{-1}\colon a\equiv-1\pmod*{(n,\Delta)}\}$. The contribution to $A_{\alpha,h}(x)$ from the primes $q\mid\Delta$ is
\begin{align*}
\sum_{q\mid \Delta}\log q\sum_{i\ge1}\sum_{\substack{p\le x\\p\in\Pp_{-1}\\p\equiv-1\pmod*{q^i}}}\frac{1}{p}&= \sum_{q\mid \Delta}\log q\sum_{i\ge1}\sum_{\substack{a\in\Rr_{-1}\\a\equiv-1\pmod*{(q^i,\Delta)}}}\sum_{\substack{p\le x\\p\equiv a\pmod*{|\Delta|}\\p\equiv -1\pmod*{q^i}}}\frac{1}{p}.
\end{align*}
To estimate the innermost sum, we appeal to the following estimate due independently to Norton \cite{Nor76} and Pomerance \cite{Pom77}:
\begin{equation}\label{eq:NP}
\sum_{\substack{m<p\le x\\ p\equiv a\pmod*{m}}}\frac{1}{p}=\frac{\log_2x}{\phi(m)}+O\left(\frac{\log3m}{\phi(m)}\right)
\end{equation}
uniformly for all $x\ge3$, all $m\in\N$, and all $a\in\Z$ coprime to $m$. If $q$ is odd, then the Chinese remainder theorem and \eqref{eq:NP} imply that 
\[\sum_{\substack{p\le x\\p\equiv a\pmod*{|\Delta|}\\p\equiv -1\pmod*{q^i}}}\frac{1}{p}=\frac{y}{\phi(|\Delta|q^{i-1})}+O\left(\frac{\log q^i}{\varphi(q^i)}\right)=\frac{y+O(\log q^i)}{\phi(|\Delta|)q^{i-1}},\]
where one observes that the contribution from the only possible prime $<|\Delta|q^{i-1}$ appearing in the sum gets subsumed into the error.
If $2^j\parallel\Delta$ for some $1\le j\le3$, then the Chinese remainder theorem and \eqref{eq:NP} yield
\[\sum_{\substack{p\le x\\p\equiv a\pmod*{|\Delta|}\\p\equiv -1\pmod*{2^i}}}\frac{1}{p}=\frac{y+O(1)}{\phi(|\Delta|)},\]
for $1\le i\le j$, and
\[\sum_{\substack{p\le x\\p\equiv a\pmod*{|\Delta|}\\p\equiv -1\pmod*{2^i}}}\frac{1}{p}=\frac{y}{\phi(|\Delta|2^{i-j})}+O\left(\frac{\log 2^i}{\varphi(2^i)}\right)=\frac{y+O(\log 2^i)}{\phi(|\Delta|)2^{i-j}}\]
for $i>j$, where $2^i-1\le a''<|\Delta|2^{i-j}$ is some positive integer. Hence, the contribution to $A_{\alpha,h}(x)$ from the primes $q\mid\Delta$ is
\[\frac{1}{\phi(|\Delta|)}\left(\sum_{\substack{q\mid \Delta\\q>2}}\beta(q)\log q\sum_{i\ge1}\frac{1}{q^{i-1}}+1_{2\mid\Delta}\log 2\left(\sum_{i=1}^{j}\beta(2^i)+\sum_{i>j}\frac{\beta(2^j)}{2^{i-j}}\right)\right)y+O(1)=c_1y+O(1),\]
where
\[c_1\colonequals\frac{1}{\phi(|\Delta|)}\left(\sum_{q\mid\Delta}\frac{\beta(q)q\log q}{q-1}+\left(1_{4\mid\Delta}(2\beta(4)-\beta(2))+1_{8\mid\Delta}(2\beta(8)-\beta(4))\right)\log 2\right).\]
Using the orthogonality relations between Dirichlet characters and the fact that $\chi$ (mod $|\Delta|$) is primitive, one can show that
\[\beta(d)=\frac{\phi(|\Delta|)(1-\chi(-1)1_{d=|\Delta|})}{2\phi(d)}\]
for any positive integer $d\mid\Delta$. As a consequence, since $\chi(-1)=\sgn(\Delta)$, we have
\[c_1=\frac{1}{2}\left(\sum_{q\mid\Delta}\frac{q\log q}{(q-1)^2}-\frac{\sgn(\Delta)1_{\Dd}(\Delta)|\Delta|\log\Rad(|\Delta|)}{\phi(|\Delta|)^2}\right),\]
where $\Dd$ is defined by \eqref{eq:Dd}. On the other hand, the contribution to $A_{\alpha,h}(x)$ from the primes $q\nmid\Delta$ is
\begin{align*}
\sum_{\substack{q\le y\log y\\q\nmid \Delta}}\log q\sum_{i\ge1}\sum_{\substack{p\le x\\p\in\Pp_{-1}\\p\equiv-1\pmod*{q^i}}}\frac{1}{p}&= \sum_{\substack{q\le y\log y\\q\nmid \Delta}}\log q\sum_{i\ge1}\sum_{a\in\Rr_{-1}}\sum_{\substack{p\le x\\p\equiv a\pmod*{|\Delta|}\\p\equiv -1\pmod*{q^i}}}\frac{1}{p}\\
&=\frac{1}{2}\sum_{\substack{q\le y\log y\\q\nmid \Delta}}\log q\sum_{i\ge1}\frac{y+O(\log q^i)}{\phi(q^i)}\\
&=\frac{y}{2}\sum_{q\le y\log y}\frac{q\log q}{(q-1)^2}-\frac{y}{2}\sum_{q\mid \Delta}\frac{q\log q}{(q-1)^2}+O\left((\log y)^2\sum_{q\le y\log y}\frac{1}{q}\right)\\
&=\frac{y}{2}\sum_{q\le y\log y}\frac{q\log q}{(q-1)^2}-\frac{y}{2}\sum_{q\mid \Delta}\frac{q\log q}{(q-1)^2}+O\left((\log y)^2\log_2y\right),
\end{align*}
by \eqref{eq:NP} and the fact that $\#\Rr_{-1}=\phi(|\Delta|)/2$. It is shown in \cite{EPS} that
\[\sum_{q\le y\log y}\frac{q\log q}{(q-1)^2}=\log y+\log_2y+\left(-\gamma+\sum_{p}\frac{\log p}{(p-1)^2}\right)+O\left(e^{-\sqrt{\log y}}\right),\]
where $\gamma=0.57721...$ is the Euler--Mascheroni constant. Hence, the contribution to $A_{\alpha,h}(x)$ from the primes $q\nmid\Delta$ is
\[\frac{1}{2}y\log y+\frac{1}{2}y\log_2y +\frac{1}{2}\left(-\gamma+\sum_{p}\frac{\log p}{(p-1)^2}-\sum_{q\mid \Delta}\frac{q\log q}{(q-1)^2}\right)y+O\left(ye^{-\sqrt{\log y}}\right).\]
We conclude that
\begin{equation}\label{eq:A_h}
A_{\alpha,h}(x)=\frac{1}{2}y\log y+\frac{1}{2}y\log_2y +\frac{c_2}{2}y+O\left(ye^{-\sqrt{\log y}}\right),
\end{equation}
where
\[c_2\colonequals\sum_{p}\frac{\log p}{(p-1)^2}-\frac{\sgn(\Delta)1_{\Dd}(\Delta)|\Delta|\log\Rad(|\Delta|)}{\phi(|\Delta|)^2}-\gamma.\]
Next, we estimate $B_{\alpha,h}(x)$. Since the primes $p\in\Pp_{0}$ contribute $O(1)$, we have
\begin{align*}
B_{\alpha,h}(x)&=\sum_{\substack{p\le x\\p\in\Pp_{-1}}}\frac{h(p)^2}{p}\left(1-\frac{\log p}{\log 3x}\right)^{-1/2}+O(1)\\
&=\sum_{q_1,q_2\le y\log y}\log q_1\log q_2\sum_{\substack{p\le x\\p\in\Pp_{-1}}}\frac{v_{q_1}(p-\chi(p))v_{q_2}(p-\chi(p))}{p}\left(1-\frac{\log p}{\log 3x}\right)^{-1/2}+O(1)\\
&=\sum_{q_1,q_2\le y\log y}\log q_1\log q_2\sum_{i,j\ge1}\sum_{\substack{p\le x\\p\in\Pp_{-1}\\p\equiv\chi(p)\pmod*{[q_1^i,q_2^j]}}}\frac{1}{p}\left(1-\frac{\log p}{\log 3x}\right)^{-1/2}+O(1)\\
&\le\sum_{q_1,q_2\le y\log y}\log q_1\log q_2\sum_{i,j\ge1}\sum_{\substack{p\le x\\p\equiv-1\pmod*{[q_1^i,q_2^j]}}}\frac{1}{p}\left(1-\frac{\log p}{\log 3x}\right)^{-1/2}+O(1),
\end{align*}
where $[q_1^i,q_2^j]$ denotes the least common multiple of $q_1^i$ and $q_2^j$. Since $(1-\log p/\log 3x)^{-1/2}\ll1$ when $p\le\sqrt{x}$, the estimates for $H_1$ and $H_2$ from \cite[p.\ 367]{EPS} lead to
\[\sum_{q_1,q_2\le y\log y}\log q_1\log q_2\sum_{i,j\ge1}\sum_{\substack{p\le \sqrt{x}\\p\equiv-1\pmod*{[q_1^i,q_2^j]}}}\frac{1}{p}\left(1-\frac{\log p}{\log 3x}\right)^{-1/2}\ll y(\log y)^2.\]
It remains to bound the tail
\[\sum_{q_1,q_2\le y\log y}\log q_1\log q_2\sum_{i,j\ge1}\sum_{\substack{\sqrt{x}<p\le x\\p\equiv-1\pmod*{[q_1^i,q_2^j]}}}\frac{1}{p}\left(1-\frac{\log p}{\log 3x}\right)^{-1/2}.\]
When $[q_1^i,q_2^j]\le x^{1/3}$, we have by Brun--Tichmarsh and partial summation that
\[\sum_{\substack{\sqrt{x}<p\le x\\p\equiv-1\pmod*{[q_1^i,q_2^j]}}}\frac{1}{p}\left(1-\frac{\log p}{\log 3x}\right)^{-1/2}\ll\frac{1}{\phi([q_1^i,q_2^j])}.\]
Borrowing the relevant estimates from \cite[p.\ 367]{EPS}, we see that the contribution to the tail from $q_1^i$ and $q_2^j$ with $[q_1^i,q_2^j]\le x^{1/3}$ is
\[\ll\sum_{q_1,q_2\le y\log y}\log q_1\log q_2\sum_{i,j\ge1}\frac{1}{\phi([q_1^i,q_2^j])}\ll y(\log y)^2.\]
On the other hand, the contribution to the tail from $q_1^i$ and $q_2^j$ with $[q_1^i,q_2^j]>x^{1/3}$ is trivially
\[(\log x)^{1/2}\sum_{q_1,q_2\le y\log y}\log q_1\log q_2\sum_{\substack{i,j\ge1\\ x^{1/3}<[q_1^i,q_2^j]\le x+1}}\sum_{\substack{\sqrt{x}<p\le x\\p\equiv-1\pmod*{[q_1^i,q_2^j]}}}\frac{1}{p}.\]
Note that
\[\sum_{\substack{\sqrt{x}<p\le x\\p\equiv-1\pmod*{[q_1^i,q_2^j]}}}\frac{1}{p}\le\sum_{\substack{n\le x\\n\equiv-1\pmod*{[q_1^i,q_2^j]}}}\frac{1}{n}\ll\frac{\log x}{[q_1^i,q_2^j]}.\]
Thus, the contribution to the tail from $q_1^i$ and $q_2^j$ with $[q_1^i,q_2^j]>x^{1/3}$ is
\[\ll(\log x)^{3/2}\sum_{q_1,q_2\le y\log y}\log q_1\log q_2\sum_{\substack{i,j\ge1\\ x^{1/3}<[q_1^i,q_2^j]\le x+1}}\frac{1}{[q_1^i,q_2^j]}.\]
The contribution to the last line from the diagonal terms with $q_1=q_2$ is
\begin{align*}
&\le2(\log x)^{3/2}\sum_{q\le y\log y}(\log q)^2\sum_{\substack{i\ge j\ge1\\ x^{1/3}<q^i\le x+1}}\frac{1}{q^i}\\
&\le2(\log x)^{3/2}\sum_{q\le y\log y}(\log q)^2\sum_{\substack{i\ge1\\ q^i>x^{1/3}}}\frac{i}{q^i}\\
&\ll\frac{(\log x)^{5/2}}{x^{1/3}}\sum_{q\le y\log y}\log q\\
&\ll\frac{(\log x)^{5/2}y\log y}{x^{1/3}},
\end{align*}
while the contribution from the off-diagonal terms with $q_1\ne q_2$ is
\begin{align*}
&\le2(\log x)^{3/2}\sum_{q_1<q_2\le y\log y}\log q_1\log q_2\sum_{\substack{i,j\ge1\\ x^{1/3}<q_1^iq_2^j\le x+1}}\frac{1}{q_1^iq_2^j}\\
&\le 2(\log x)^{3/2}\sum_{q_1<q_2\le y\log y}\log q_1\log q_2\sum_{q_1^i\le x}\frac{1}{q_1^i}\sum_{q_2^j>x^{1/3}/q_1^i}\frac{1}{q_2^j}\\
&\ll\frac{(\log x)^{3/2}}{x^{1/3}}\sum_{q_1<q_2\le y\log y}\log q_1\log q_2\sum_{q_1^i\le x}1\\
&\le\frac{(\log x)^{5/2}}{x^{1/3}}\sum_{q_1,q_2\le y\log y}\log q_2\\
&\ll\frac{(\log x)^{5/2}y^2\log y}{x^{1/3}}.
\end{align*}
Collecting all the contributions to the tail shows that the tail is $\ll y(\log y)^2$, whence $B_{\alpha,h}(x)\ll y(\log y)^2$. Combining this estimate with \eqref{eq:A_h} and invoking Theorem \ref{thm:WTK}, we conclude that 
\begin{equation}\label{eq:normalh}
\left|h(f)-\frac{1}{2}y\log y-\frac{1}{2}y\log_2y -\frac{c_2}{2}y\right|<\frac{y}{\log y}
\end{equation}
holds for all but $o(S_{\alpha}(x))$ split-free positive integers $f\le x$. 

\subsection{The contribution from $I_1$}\label{subsec:I_1}
With \eqref{eq:normalh} in hand, we proceed to estimate the contribution to $\log L(f)$ from the primes in $I_1$. To this end, we show that 
\begin{equation}\label{eq:I_1a}
\sum_{\substack{q^i>y^2/\log^2y\\ i>1,\,q^i\parallel L(f)}}\log q^i<(\log y)^2  
\end{equation}
for all but $o(S_{\alpha}(x))$ split-free $f\le x$. Once we have this, we can conclude that the contribution to $\log L(f)$ from the primes $q\in I_1$ is at most
\begin{equation}\label{eq:conI_1}
\sum_{\substack{q\in I_1,\,q^i\parallel L(f)\\q^i\le y^2/\log^2y}}\log q^i+(\log y)^2\ll\sum_{q\in I_1}\log y+(\log y)^2\ll\frac{y}{\log y}    \end{equation}
for all but $o(S_{\alpha}(x))$ split-free $f\le x$. Now we prove \eqref{eq:I_1a} by averaging the left-hand side over split-free $f\le x$. More precisely, we show
\begin{equation}\label{eq:I_1b}
S_{\alpha}(x)^{-1}\sum_{f\le x}\alpha(f)\sum_{\substack{q^i>y^2/\log^2y\\ i>1,\,q^i\parallel L(f)}}\log q^i\ll \log y,
\end{equation}
from which \eqref{eq:I_1a} follows at once. Note that if $q^i\mid L(f)$, then $q^i\mid f$ with $q$ ramified, or $q^{i+1}\mid f$ with $q$ unramified, or $f$ has a prime factor $p$ with $p\equiv\chi(p)\pmod{q^i}$. The first case contributes at most
\begin{align*}
S_{\alpha}(x)^{-1}\sum_{\substack{q^i>y^2/\log^2y\\q\mid \Delta,\,i>1}}\log q^i\sum_{\substack{f\le x\\ q^{i}\mid f}}\alpha(f)&\ll\sum_{\substack{y^2/\log^2y<q^i\le x\\ q\mid \Delta,\,i>1}}\frac{\log q^i}{q^{i}}\left(1-\frac{\log q^{i}}{\log 3x}\right)^{-1/2}\\
&\ll\sum_{\substack{y^2/\log^2y<q^i\le \sqrt{x}\\q\mid \Delta,\,i>1}}\frac{\log q^i}{q^{i}}+(\log x)^{1/2}\sum_{\substack{q^i>\sqrt{x}\\ q\mid \Delta}}\frac{\log q^i}{q^{i}}\\
&\ll\sum_{q\mid\Delta}\sum_{i>1}\frac{\log q^i}{q^{i}}+\frac{(\log x)^{3/2}}{\sqrt{x}}\ll1,
\end{align*}
Similarly, the contribution to the left-hand side of \eqref{eq:I_1b} from the second case is
\begin{align*}
&\le S_{\alpha}(x)^{-1}\sum_{\substack{q^i>y^2/\log^2y\\ i>1}}\log q^i\sum_{\substack{f\le x\\ q^{i+1}\mid f}}\alpha(f)\\
&\ll\sum_{\substack{y^2/\log^2y<q^i\le x/q\\ i>1}}\frac{\log q^i}{q^{i+1}}\left(1-\frac{\log q^{i+1}}{\log 3x}\right)^{-1/2}\\
&\ll\sum_{\substack{y^2/\log^2y<q^i\le \sqrt{x}/q\\ i>1}}\frac{\log q^i}{q^{i+1}}+(\log x)^{1/2}\sum_{\substack{\sqrt{x}<q^{i+1}\le x\\ i>1}}\frac{\log q^i}{q^{i+1}}\\
&\ll\sum_{q}\sum_{i>1}\frac{\log q^i}{q^{i+1}}+\frac{(\log x)^{3/2}}{\sqrt{x}}\sum_{q\le x^{1/3}}1\ll1.
\end{align*}
Finally, the contribution to the left-hand side of \eqref{eq:I_1b} from the last case is
\begin{align*}
&\le S_{\alpha}(x)^{-1}\sum_{\substack{q^i>y^2/\log^2y\\ i>1}}\log q^i\sum_{\substack{p\le x\\ p\equiv\chi(p)\pmod*{q^i}\\p\notin\Pp_{1}}}\sum_{\substack{f\le x\\ p\mid f}}\alpha(f)\\
&= S_{\alpha}(x)^{-1}\sum_{\substack{y^2/\log^2y<q^i\le x+1\\ i>1}}\log q^i\sum_{\substack{p\le x\\p\equiv\chi(p)\pmod*{q^i}}}\alpha(p)\sum_{f\le x/p}\alpha(f)\\
&\ll\sum_{\substack{q^i>y^2/\log^2y\\ i>1}}\log q^i\sum_{\substack{p\le x\\p\equiv\chi(p)\pmod*{q^i}}}\frac{\alpha(p)}{p}\left(1-\frac{\log p}{\log 3x}\right)^{-1/2}\\
&\le\sum_{\substack{y^2/\log^2y<q^i\le x+1\\ i>1}}\log q^i\sum_{\substack{p\le x\\p\equiv-1\pmod*{q^i}}}\frac{1}{p}\left(1-\frac{\log p}{\log 3x}\right)^{-1/2}\\
&\ll\sum_{\substack{y^2/\log^2y<q^i\le x+1\\ i>1}}\log q^i\sum_{\substack{p\le \sqrt{x}\\p\equiv-1\pmod*{q^i}}}\frac{1}{p}+\sum_{\substack{y^2/\log^2y<q^i\le x+1\\ i>1}}\log q^i\sum_{\substack{\sqrt{x}<p\le x\\p\equiv-1\pmod*{q^i}}}\frac{1}{p}\left(1-\frac{\log p}{\log 3x}\right)^{-1/2}.
\end{align*} 
The first double sum on the last line follows easily from Brun--Titchmarsh and partial summation. Indeed, the first double sum is
\begin{align*}
\ll y\sum_{\substack{q^i>y^2/\log^2y\\ i>1}}\frac{\log q^i}{\phi(q^i)}&\le y\sum_{\substack{q^i>y^2/\log^2y\\ q\le y/\log y}}\frac{\log q^i}{\phi(q^i)}+y\sum_{\substack{q> y/\log y\\i>1}}\frac{\log q^i}{\phi(q^i)}\\
&\ll\frac{y\log y}{y^2/\log^2y}\sum_{q\le y/\log y}1+y\sum_{q> y/\log y}\frac{\log q}{q^2}\\
&\ll \log y.
\end{align*}
We split the second double sum into the two following subsums
\begin{align*}
&\sum_{\substack{y^2/\log^2y<q^i\le x^{1/3}\\ i>1}}\log q^i\sum_{\substack{\sqrt{x}<p\le x\\p\equiv-1\pmod*{q^i}}}\frac{1}{p}\left(1-\frac{\log p}{\log 3x}\right)^{-1/2},\\
&\sum_{\substack{x^{1/3}<q^i\le x+1\\ i>1}}\log q^i\sum_{\substack{\sqrt{x}<p\le x\\p\equiv-1\pmod*{q^i}}}\frac{1}{p}\left(1-\frac{\log p}{\log 3x}\right)^{-1/2}.
\end{align*}
By Brun--Titchmarsh and partial summation, we have
\[\sum_{\substack{\sqrt{x}<p\le x\\p\equiv-1\pmod*{q^i}}}\frac{1}{p}\left(1-\frac{\log p}{\log 3x}\right)^{-1/2}\ll\frac{1}{\phi(q^i)}\]
whenever $q^i\le x^{1/3}$, the first subsum above is
\[\ll\sum_{\substack{q^i>y^2/\log^2y\\ i>1}}\frac{\log q^i}{\phi(q^i)}\ll\frac{\log y}{y}.\]
Similarly, the second subsum is 
\begin{align*}
&\ll(\log x)^{1/2}\sum_{\substack{x^{1/3}<q^i\le x+1\\ i>1}}\log q^i\sum_{\substack{\sqrt{x}<p\le x\\p\equiv-1\pmod*{q^i}}}\frac{1}{p}\\
&\ll y(\log x)^{1/2}\sum_{\substack{q^i>x^{1/3}\\ i>1}}\frac{\log q^i}{\phi(q^i)}\\
&\le y(\log x)^{1/2}\sum_{\substack{q\le x^{1/6}\\q^i>x^{1/3}}}\frac{\log q^i}{\phi(q^i)}+y(\log x)^{1/2}\sum_{\substack{q>x^{1/6}\\i>1}}\frac{\log q^i}{\phi(q^i)}\\
&\ll\frac{y(\log x)^{3/2}}{x^{1/3}}\sum_{q\le x^{1/6}}1+y(\log x)^{1/2}\sum_{q>x^{1/6}}\frac{\log q}{q^2}\\
&\ll\frac{y(\log x)^{1/2}}{x^{1/6}}.
\end{align*}
Gathering the estimates above completes the proof of \eqref{eq:I_1b}. 

\subsection{The contribution from $I_2$}\label{subsec:I_2}
We now turn to the most elaborate part of the proof: estimating the contribution to $\log L(f)$ from the primes in $I_2$. Following the argument in \cite{EPS}, we define analogously
\begin{align*}
\Qq(q)&\colonequals\{p\le x\colon p\notin\Pp_1\text{~and~}p\equiv\chi(p)\pmod*{q}\},\\
\Qq_1(q)&\colonequals\{p\in \Qq(q)\cap[2,x^{1/y}]\colon p\not\equiv\chi(p)\pmod*{qq'}\text{~for any~}q'\in I_2\},\\
\Qq_2(q)&\colonequals\{p\in \Qq(q)\colon p\equiv\chi(p)\pmod*{qq'}\text{~for some~}q'\in I_2\},\\
\Qq_3(q)&\colonequals\{p\in \Qq(q)\cap(x^{1/y},x]\colon p\not\equiv\chi(p)\pmod*{qq'}\text{~for any~}q'\in I_2\},
\end{align*}
for each prime $q\in I_2$, so that $\Qq(q)=\Qq_1(q)\cup\Qq_2(q)\cup\Qq_3(q)$. Moreover, for every $p\in\Qq(q)$ we have $p\equiv-1\pmod*{q}$, and in particular, $p>q-1>y/\log y-1$. In view of \eqref{eq:I_1a}, we know that the contribution to $\log L(f)$ from the primes $q\in I_2$ is
\begin{equation}\label{eq:conI_2a}
\sum_{\substack{q\in I_2\\(f,\Qq_1(q))>1}}\log q+O\left(\sum_{q\in I_2}\sum_{\substack{p\mid f\\ p\in\Qq_2(q)}}\log q+\sum_{q\in I_2}\sum_{\substack{p\mid f\\ p\in\Qq_3(q)}}\log q\right)+O\left((\log y)^2\right)
\end{equation}
for all but $o(S_{\alpha}(x))$ split-free $f\le x$. 

We first show that the second big $O$ term is negligible by averaging. We start by proving
\begin{equation}\label{eq:conQ_2}
\sum_{q\in I_2}\sum_{\substack{p\mid f\\ p\in\Qq_2(q)}}\log q<\frac{y(\log_2y)^3}{\log y}
\end{equation}
for all but $o(S_{\alpha}(x))$ split-free $f\le x$. The proof is straightforward. We have
\begin{align*}
&\hspace*{5.5mm}S_{\alpha}(x)^{-1}\sum_{f\le x}\alpha(f)\sum_{q\in I_2}\sum_{\substack{p\mid f\\ p\in\Qq_2(q)}}\log q\\
&=S_{\alpha}(x)^{-1}\sum_{q\in I_2}\log q\sum_{p\in\Qq_2(q)}\sum_{\substack{f\le x\\ p\mid f}}\alpha(f)\\
&\ll\sum_{q,q'\in I_2}\log q\sum_{\substack{p\le x\\p\notin\Pp_1\\p\equiv\chi(p)\pmod*{qq'}}}\frac{1}{p}\left(1-\frac{\log p}{\log 3x}\right)^{-1/2}\\
&\le\sum_{q,q'\in I_2}\log q\sum_{\substack{p\le x\\p\equiv-1\pmod*{qq'}}}\frac{1}{p}\left(1-\frac{\log p}{\log 3x}\right)^{-1/2}\\
&\ll\sum_{q,q'\in I_2}\log q\sum_{\substack{p\le \sqrt{x}\\p\equiv-1\pmod*{qq'}}}\frac{1}{p}+\sum_{q,q'\in I_2}\log q\sum_{\substack{\sqrt{x}<p\le x\\p\equiv-1\pmod*{qq'}}}\frac{1}{p}\left(1-\frac{\log p}{\log 3x}\right)^{-1/2}.
\end{align*}
It follows from Brun--Tichmarsh and partial summation that the first double sum is 
\[\ll y\sum_{q,q'\in I_2}\frac{\log q}{\phi(qq')}\ll y\log y\left(\sum_{q\in I_2}\frac{1}{q}\right)^2\ll\frac{y(\log_2y)^2}{\log y},\]
and the second double sum is 
\[\ll \sum_{q,q'\in I_2}\frac{\log q}{\phi(qq')}\ll\frac{(\log_2y)^2}{\log y}.\]
Hence, we obtain
\[S_{\alpha}(x)^{-1}\sum_{f\le x}\alpha(n)\sum_{q\in I_2}\sum_{\substack{p\mid f\\ p\in\Qq_2(q)}}\log q\ll\frac{y(\log_2y)^2}{\log y},\]
of which \eqref{eq:conQ_2} is an immediate consequence.

In a similar fashion, we show that
\begin{equation}\label{eq:conQ_3}
\sum_{q\in I_2}\sum_{\substack{p\mid f\\ p\in\Qq_3(q)}}\log q<(\log y)^2.
\end{equation}
for all but $o(S_{\alpha}(x))$ split-free $f\le x$. Indeed, we have
\begin{align*}
&\hspace*{5.5mm}S_{\alpha}(x)^{-1}\sum_{f\le x}\alpha(f)\sum_{q\in I_2}\sum_{\substack{p\mid f\\ p\in\Qq_3(q)}}\log q\\
&=S_{\alpha}(x)^{-1}\sum_{q\in I_2}\log q\sum_{p\in\Qq_3(q)}\sum_{\substack{f\le x\\ p\mid f}}\alpha(f)\\
&\ll\sum_{q\in I_2}\log q\sum_{\substack{x^{1/y}<p\le x\\p\equiv-1\pmod*{q}}}\frac{1}{p}\left(1-\frac{\log p}{\log 3x}\right)^{-1/2}\\
&\ll\sum_{q\in I_2}\log q\sum_{\substack{x^{1/y}<p\le \sqrt{x}\\p\equiv-1\pmod*{q}}}\frac{1}{p}+\sum_{q\in I_2}\log q\sum_{\substack{\sqrt{x}<p\le x\\p\equiv-1\pmod*{q}}}\frac{1}{p}\left(1-\frac{\log p}{\log 3x}\right)^{-1/2}\\
&\ll\log y\sum_{q\in I_2}\frac{\log q}{\phi(q)}+\sum_{q\in I_2}\frac{\log q}{\phi(q)}\\
&\ll(\log y)\log_2y,
\end{align*}
which implies \eqref{eq:conQ_3}.

Inserting \eqref{eq:conQ_2} and \eqref{eq:conQ_3} into \eqref{eq:conI_2a}, we see that the contribution to $\log L(f)$ from the primes $q\in I_2$ is
\begin{equation}\label{eq:conI_2b}
\sum_{\substack{q\in I_2\\(f,\Qq_1(q))>1}}\log q+O\left(\frac{y(\log_2y)^3}{\log y}\right)
\end{equation}
for all but $o(S_{\alpha}(x))$ split-free $f\le x$. Now we estimate the sum in \eqref{eq:conI_2b}. Let
\[g(f)\colonequals\sum_{\substack{q\in I_2\\(f,\Qq_1(q))=1}}1.\]
Although $g$ is not additive, we can still determine its typical size by computing its first and second moments weighted by $\alpha$. We begin with its first moment 
\[S_{\alpha}(x)^{-1}\sum_{f\le x}\alpha(f)g(f)=S_{\alpha}(x)^{-1}\sum_{q\in I_2}\sum_{\substack{f\le x\\(f,\Qq_1(q))=1}}\alpha(f).\]
The fundamental lemma of sieve theory \cite[Theorem 19.1]{Kou19} ensures that there exist upper and lower bound sieve weights $\lambda^{\pm}$ (not to be confused with Carmichael's $\lambda$-function) of level $D=x^{u/y}$ with $u=\log y$, which are 1-bounded and supported on positive integers $n\le D$ composed entirely of primes factors from $\Qq_1(q)$, such that $1\ast\lambda^-\le 1_{(n,\Qq_1(q))=1}\le1\ast\lambda^+$ and 
\begin{equation}\label{eq:FL}
\sum_{d\le D}\frac{\lambda^{\pm}(d)\alpha(d)}{d}=\left(1+O\left(u^{-u/2}\right)\right)\prod_{p\in \Qq_1(q)}\left(1-\frac{1}{p}\right),    
\end{equation}
where $\ast$ denotes the Dirichlet convolution.
 Then we have
\begin{equation}\label{eq:lambda+}
\sum_{\substack{f\le x\\(f,\Qq_1(q))=1}}\alpha(f)=\sum_{f\le x}\alpha(f)1_{(f,\Qq_1(q))=1}\le\sum_{f\le x}\alpha(f)\sum_{d\mid f}\lambda^+(d)=\sum_{d\le D}\lambda^+(d)\sum_{\substack{f\le x\\ d\mid f}}\alpha(f).    
\end{equation}
Note that
\begin{align*}
\sum_{d\le D}\lambda^+(d)\sum_{\substack{f\le x\\ d\mid f}}\alpha(f)&=\sum_{d\le D}\lambda^+(d)\alpha(d)c_{\alpha}\frac{x}{d}\left(\log\frac{x}{d}\right)^{-1/2}\left(1+O\left(\left(\log\frac{x}{d}\right)^{-1}\right)\right)\\
&=c_{\alpha}x(\log x)^{-1/2}\sum_{d\le D}\frac{\lambda^+(d)\alpha(d)}{d}\left(1+O\left(\frac{\log D}{\log x}\right)\right)\left(1+O\left(\frac{1}{\log x}\right)\right)\\
&=c_{\alpha}x(\log x)^{-1/2}\sum_{d\le D}\frac{\lambda^+(d)\alpha(d)}{d}\left(1+O\left(\frac{\log y}{y}\right)\right)\\
&=c_{\alpha}x(\log x)^{-1/2}\left(\sum_{d\le D}\frac{\lambda^+(d)\alpha(d)}{d}+O\left(\frac{\log y}{y}\sum_{d\le D}\frac{|\lambda^+(d)|\alpha(d)}{d}\right)\right).
\end{align*}
For each $q\in I_2$, we have by \eqref{eq:NP} that
\[\sum_{p\in \Qq_1(q)}\frac{1}{p}\le\sum_{p\in \Qq(q)}\frac{1}{p}\le\sum_{\substack{p\le x\\p\in\Pp_{-1}\\p\equiv-1\pmod*{q}}}\frac{1}{p}=\frac{\phi(|\Delta|)}{2}\cdot\frac{y+O(\log q)}{\phi(|\Delta|)\phi(q)}=\frac{y+O(\log y)}{2\phi(q)}=\frac{y+O(\log y)}{2q},\]
where the last equality results from the inequality $y/(q\phi(q))<2(\log y)/q$. It follows that
\[\sum_{d\le D}\frac{|\lambda^+(d)|\alpha(d)}{d}\le  \sum_{\substack{d\ge1\\p\mid d\Rightarrow p\in\Qq_1(q)}}\frac{1}{d}\le\prod_{p\in \Qq_1(q)}\left(1-\frac{1}{p}\right)^{-1}\ll\exp\left(\sum_{p\in \Qq_1(q)}\frac{1}{p}\right)\ll e^{\frac{y}{2q}},\]
whence
\[\sum_{d\le D}\lambda^+(d)\sum_{\substack{f\le x\\ d\mid f}}\alpha(f)=c_{\alpha}x(\log x)^{-1/2}\left(\sum_{d\le D}\frac{\lambda^+(d)\alpha(d)}{d}+O\left(\frac{\log y}{y}e^{\frac{y}{2q}}\right)\right).\]
Inserting this into \eqref{eq:lambda+} and applying \eqref{eq:FL}, we have
\[\sum_{\substack{f\le x\\(f,\Qq_1(q))=1}}\alpha(f)\le S_{\alpha}(x)\left(1+O\left(u^{-u/2}\right)\right)\prod_{p\in \Qq_1(q)}\left(1-\frac{1}{p}\right)+O\left(\frac{S_{\alpha}(x)\log y}{y}e^{\frac{y}{2q}}\right).\]
Arguing with $\lambda^-$ in place of $\lambda^+$ leads to the reversed inequality. Therefore, we obtain
\[\sum_{\substack{f\le x\\(f,\Qq_1(q))=1}}\alpha(f)=S_{\alpha}(x)\left(1+O\left(u^{-u/2}\right)\right)\prod_{p\in \Qq_1(q)}\left(1-\frac{1}{p}\right)+O\left(\frac{S_{\alpha}(x)\log y}{y}e^{\frac{y}{2q}}\right),\]
Summing on $q\in I_2$, we conclude that
\begin{equation}\label{eq:g1a}
S_{\alpha}(x)^{-1}\sum_{f\le x}\alpha(f)g(f)=\left(1+O\left(u^{-u/2}\right)\right)\sum_{q\in I_2}\prod_{p\in \Qq_1(q)}\left(1-\frac{1}{p}\right)+O\left(\frac{\sqrt{y}}{\log y}\right),
\end{equation}
since
\begin{equation}\label{eq:sumI_2e^{y/2q}}
\sum_{q\in I_2}e^{\frac{y}{2q}}\le\sum_{y/\log y<q\le 2y/\log y}e^{\frac{\log y}{2}}+\sum_{2y/\log y<q\le y\log y}e^{\frac{\log y}{4}}\ll\frac{y^{3/2}}{(\log y)^2}.   
\end{equation}
To complete the estimation of the first weighted moment of $g$, it suffices to estimate
\[\sum_{q\in I_2}\prod_{p\in \Qq_1(q)}\left(1-\frac{1}{p}\right).\]
A similar sum is handled in \cite{EPS}, and the method there can be easily adapted to treat our sum above. So we will be brief. We start by observing that
\begin{align*}
\sum_{p\in \Qq_1(q)}\frac{1}{p}&=\sum_{\substack{p\le x^{1/y}\\p\notin\Pp_1\\p\equiv\chi(p)\pmod*{q}}}\frac{1}{p}-\sum_{\substack{p\le x^{1/y}\\p\in \Qq_2(q)}}\frac{1}{p}\\
&=\sum_{\substack{p\le x^{1/y}\\p\in\Pp_{-1}\\p\equiv-1\pmod*{q}}}\frac{1}{p}+O\left(\sum_{q'\in I_2}\sum_{\substack{p\le x\\p\equiv-1\pmod*{qq'}}}\frac{1}{p}\right)\\
&=\frac{\phi(|\Delta|)}{2}\cdot\frac{\log\log x^{1/y}+O(\log q)}{\phi(|\Delta|)\phi(q)}+O\left(y\sum_{q'\in I_2}\frac{1}{\phi(qq')}\right)\\
&=\frac{y}{2q}+O\left(\frac{y\log_2y}{q\log y}\right).
\end{align*}
Thus, following \cite[pp.\ 369--370]{EPS}, we have
\begin{align}\label{eq:g1b}
\sum_{q\in I_2}\prod_{p\in \Qq_1(q)}\left(1-\frac{1}{p}\right)&=\sum_{q\in I_2}\exp\left(-\frac{y}{2q}+O\left(\frac{y\log_2y}{q\log y}\right)\right)\nonumber\\  
&=\sum_{q\in I_2}e^{-\frac{y}{2q}}\left(1+O\left(\frac{y\log_2y}{q\log y}\right)\right)+O\left(\sum_{q\le y/4\log_2y}\frac{1}{(\log y)^2}\right)\nonumber\\
&=\sum_{q\in I_2}e^{-\frac{y}{2q}}+O\left(\frac{y(\log_2y)^2}{(\log y)^2}\right).
\end{align}
By partial summation, we see that
\begin{equation}\label{eq:sume^{-y/2t}}
\sum_{q\in I_2}e^{-\frac{y}{2q}}=\left(\pi(y\log y)-\pi\left(\frac{y}{\log y}\right)\right)e^{-\frac{1}{2\log y}}-\int_{y/\log y}^{y\log y}e^{-\frac{y}{2t}}\frac{y}{2t^2}\left(\pi(t)-\pi\left(\frac{y}{\log y}\right)\right)\,dt.    
\end{equation}
By the Prime Number Theorem, we have
\[\left(\pi(y\log y)-\pi\left(\frac{y}{\log y}\right)\right)e^{-\frac{1}{2\log y}}=y-\frac{y\log_2y}{\log y}+\frac{y}{2\log y}+O\left(\frac{y(\log_2y)^2}{(\log y)^2}\right).\]
Following \cite[p 371]{EPS}, the integral in \eqref{eq:sume^{-y/2t}} equals
\begin{align*}
&\hspace*{5.5mm}\int_{y/\log y}^{y\log y}e^{-\frac{y}{2t}}\frac{y}{2t^2}\left(\frac{t}{\log t}+O\left(\frac{t}{(\log t)^2}\right)\right)\,dt-\pi\left(\frac{y}{\log y}\right)\left(e^{-\frac{1}{2\log y}}-e^{-\frac{\log y}{2}}\right)\\
&=\int_{y/\log y}^{y\log y}e^{-\frac{y}{2t}}\frac{y}{2t}\left(\frac{1}{\log y}+O\left(\frac{\log_2y}{(\log y)^2}\right)\right)\,dt+O\left(\frac{y}{(\log y)^2}\right)\\
&=\frac{y}{2\log y}\int_{2/\log y}^{2\log y}e^{-\frac{1}{v}}\frac{dv}{v}+O\left(\frac{y(\log_2y)^2}{(\log y)^2}\right)\\
&=\frac{y}{2\log y}\left(e^{-\frac{1}{2\log y}}\log(2\log y)+e^{-\frac{\log y}{2}}\log\frac{\log y}{2}-\int_{2/\log y}^{2\log y}e^{-\frac{1}{v}}\frac{\log v}{v^2}\,dv\right)+O\left(\frac{y(\log_2y)^2}{(\log y)^2}\right)\\
&=\frac{y}{2\log y}\left(\log_2y+\log2-\int_{0}^{\infty}e^{-\frac{1}{v}}\frac{\log v}{v^2}\,dv\right)+O\left(\frac{y(\log_2y)^2}{(\log y)^2}\right)\\
&=\frac{y\log_2y}{2\log y}+\frac{1}{2}(\log 2-\gamma)\frac{y}{\log y}+O\left(\frac{y(\log_2y)^2}{(\log y)^2}\right),
\end{align*}
by partial integration and the Prime Number Theorem, where we have exploited the identity
\[\gamma=-\Gamma'(1)=-\int_{0}^{\infty}e^{-t}\log t\,dt=\int_{0}^{\infty}e^{-\frac{1}{v}}\frac{\log v}{v^2}\,dv,\]
with $\Gamma$ being the Gamma function. Combining these estimates with \eqref{eq:g1b} and \eqref{eq:sume^{-y/2t}} yields
\[\sum_{q\in I_2}\prod_{p\in \Qq_1(q)}\left(1-\frac{1}{p}\right)=y-\frac{3y\log_2y}{2\log y}+\frac{1}{2}(1-\log 2+\gamma)\frac{y}{\log y}+O\left(\frac{y(\log_2y)^2}{(\log y)^2}\right).\]
Carrying this back into \eqref{eq:g1a}, we conclude that the first weighted moment of $g$ is
\begin{equation}\label{eq:g1c}
S_{\alpha}(x)^{-1}\sum_{f\le x}\alpha(f)g(f)=y-\frac{3y\log_2y}{2\log y}+\frac{1}{2}(1-\log 2+\gamma)\frac{y}{\log y}+O\left(\frac{y(\log_2y)^2}{(\log y)^2}\right).    
\end{equation}

Next, we estimate the second weighted moment of $g$ defined by
\[S_{\alpha}(x)^{-1}\sum_{f\le x}\alpha(f)g(f)^2=S_{\alpha}(x)^{-1}\sum_{q_1,q_2\in I_2}\sum_{\substack{f\le x\\(f,\Qq_1(q_1)\cup\Qq_1(q_2))=1}}\alpha(f).\]
The contribution from the diagonal terms with $q_1=q_2$ is exactly the first weighted moment of $g$ whose estimate is provided by \eqref{eq:g1c}. On the other hand, a similar application of the fundamental lemma of sieve theory to $\Qq_1(q_1)\cup\Qq_1(q_2)$, but with $D=x^{u/y}$ and $u=\log_2y$ this time, implies that the contribution from the off-diagonal terms with $q_1\ne q_2$ is
\[\left(1+O\left(u^{-u/2}\right)\right)\sum_{\substack{q_1,q_2\in I_2\\q_1\ne q_2}}\prod_{p\in \Qq_1(q_1)\cup\Qq_1(q_2)}\left(1-\frac{1}{p}\right)+O\left(\frac{\log_2 y}{y}\sum_{\substack{q_1,q_2\in I_2\\q_1\ne q_2}}e^{\frac{y}{2q_1}+\frac{y}{2q_2}}\right),\]
where the second error term is evidently $O(y^2\log_2y/(\log y)^4)$ by \eqref{eq:sumI_2e^{y/2q}}. Since $\Qq_1(q_1)\cap\Qq_1(q_2)=\emptyset$ for $q_1\ne q_2$, we have
\begin{align*}
&\hspace*{5mm}\left(1+O\left(u^{-u/2}\right)\right)\sum_{\substack{q_1,q_2\in I_2\\q_1\ne q_2}}\prod_{p\in \Qq_1(q_1)\cup\Qq_1(q_2)}\left(1-\frac{1}{p}\right)\\
&=\left(1+O\left(u^{-u/2}\right)\right)\left(\sum_{q\in I_2}\prod_{p\in \Qq_1(q)}\left(1-\frac{1}{p}\right)\right)^2-\left(1+O\left(u^{-u/2}\right)\right)\sum_{q\in I_2}\prod_{p\in \Qq_1(q)}\left(1-\frac{1}{p}\right)^2\\
&=\left(1+O\left(u^{-u/2}\right)\right)\left(S_{\alpha}(x)^{-1}\sum_{f\le x}\alpha(f)g(f)+O\left(\frac{\sqrt{y}}{\log y}\right)\right)^2+O(\pi(y\log y))\\
&=\left(S_{\alpha}(x)^{-1}\sum_{f\le x}\alpha(f)g(f)\right)^2+O\left(y^2u^{-u/2}\right),
\end{align*}
by \eqref{eq:g1a} and \eqref{eq:g1c}. Thus, the contribution from the off-diagonal terms with $q_1\ne q_2$ is
\[\left(S_{\alpha}(x)^{-1}\sum_{f\le x}\alpha(f)g(f)\right)^2+O\left(\frac{y^2\log_2y}{(\log y)^4}\right).\]
It follows that the second weighted moment of $g$ is
\[S_{\alpha}(x)^{-1}\sum_{f\le x}\alpha(f)g(f)^2=\left(S_{\alpha}(x)^{-1}\sum_{f\le x}\alpha(f)g(f)\right)^2+O\left(\frac{y^2\log_2y}{(\log y)^4}\right).\]
Therefore,
\begin{align*}
&\hspace*{5.5mm}S_{\alpha}(x)^{-1}\sum_{f\le x}\alpha(f)\left(g(f)-S_{\alpha}(x)^{-1}\sum_{f\le x}\alpha(f)g(f)\right)^2\\
&=S_{\alpha}(x)^{-1}\sum_{f\le x}\alpha(f)g(f)^2-\left(S_{\alpha}(x)^{-1}\sum_{f\le x}\alpha(f)g(f)\right)^2\ll \frac{y^2\log_2y}{(\log y)^4},
\end{align*}
which in conjunction with \eqref{eq:g1c} allows us to deduce that
\begin{equation}\label{eq:normalg}
\left|g(f)-\left(y-\frac{3y\log_2y}{2\log y}+\frac{1}{2}(1-\log 2+\gamma)\frac{y}{\log y}\right)\right|<\frac{y(\log_2y)^3}{(\log y)^2}
\end{equation}
holds for all but $o(S_{\alpha}(x))$ split-free $f\le x$. 

Combining \eqref{eq:conI_2b} with \eqref{eq:normalg}, we find, as in \cite{EPS}, that the contribution to $\log L(f)$ from the primes $q\in I_2$ is
\begin{align}\label{eq:conI_2c}
&\hspace*{5.5mm}(\log y+O(\log_2y))\sum_{\substack{q\in I_2\\(f,\Qq_1(q))>1}}1+O\left(\frac{y(\log_2y)^3}{\log y}\right)\nonumber\\
&=(\log y+O(\log_2y))\left(\pi(y\log y)-\pi\left(\frac{y}{\log y}\right)-g(f)\right)+O\left(\frac{y(\log_2y)^3}{\log y}\right)\nonumber\\
&=(\log y+O(\log_2y))\left(\frac{y\log_2y}{2\log y}+\frac{1}{2}(1+\log 2-\gamma)\frac{y}{\log y}+O\left(\frac{y(\log_2y)^3}{(\log y)^2}\right)\right)+O\left(\frac{y(\log_2y)^3}{\log y}\right)\nonumber\\
&=\frac{1}{2}y\log_2y+\frac{1}{2}(1+\log 2-\gamma)y+O\left(\frac{y(\log_2y)^3}{\log y}\right).
\end{align}

\subsection{Contribution from $I_3$}\label{subsec:I_3}
Now we estimate the contribution to the left-hand side of \eqref{eq:psi-L} from the primes $q\in I_3$. As remarked at the beginning of Section \ref{sec:Lforder}, we have that $q^2\nmid f$ for any $q\notin I_1$ for all but $o(S_{\alpha}(x))$ split-free $f\le x$. So we may assume $q^2\nmid f$ for any $q\in I_3$. We may assume further that $q^2\nmid L(f)$ for any $q\in I_3$, thanks to the estimate
\begin{align*}
&\hspace*{5.5mm}S_{\alpha}(x)^{-1}\sum_{q\in I_3}\sum_{\substack{p\le x\\p\equiv\chi(p)\pmod*{q^2}}}\sum_{\substack{f\le x\\p\mid f}}\alpha(f)\\
&\ll\sum_{q\in I_3}\sum_{\substack{p\le x\\p\equiv\chi(p)\pmod*{q^2}}}\frac{\alpha(p)}{p}\left(1-\frac{\log p}{\log 3x}\right)^{-1/2}\\
&\le\sum_{q\in I_3}\sum_{\substack{p\le x\\p\equiv-1\pmod*{q^2}}}\frac{1}{p}\left(1-\frac{\log p}{\log 3x}\right)^{-1/2}\\
&\ll\sum_{q\in I_3}\sum_{\substack{p\le \sqrt{x}\\p\equiv-1\pmod*{q^2}}}\frac{1}{p}+\sum_{q\in I_3}\sum_{\substack{\sqrt{x}<p\le x\\p\equiv-1\pmod*{q^2}}}\frac{1}{p}\left(1-\frac{\log p}{\log 3x}\right)^{-1/2}\\
&\ll y\sum_{q>y\log y}\frac{1}{\phi(q^2)}+\sum_{q>y\log y}\frac{1}{\phi(q^2)}\\
&\ll\frac{1}{(\log y)^2}.
\end{align*}
If $q^2\nmid f$ and $v_q(L(f))=1$, then $p\not\equiv\chi(p)\pmod*{q^2}$ for any $p\mid f$, and the number $\omega(f,\Qq(q))$ of primes $p\in\Qq(q)$ dividing $f$ is $v_q(\psi(f))$. 
Thus, setting
\[G(f)\colonequals\sum_{\substack{q\in I_3\\\omega(f,\Qq(q))>1}}\omega(f,\Qq(q))\log q,\]
we have
\[\sum_{q\in I_3}\left(v_q(\psi(f))-v_q(L(f))\right)\log q=\sum_{\substack{q\in I_3\\v_q(L(f))=1}}\left(v_q(\psi(f))-1\right)\log q\le G(f)\]
for those $f$ with $q^2\nmid f$ and $q^2\nmid L(f)$ for any $q\in I_3$. We show that the mean value of $G$ over split-free $f\le x$ is
\begin{equation}\label{eq:G(f)}
 S_{\alpha}(x)^{-1}\sum_{f\le x}\alpha(f)G(f)\ll\frac{y}{\log y}.   
\end{equation}
From this it follows that $G(f)<y\log_2y/\log y$ for all but $o(S_{\alpha}(x))$ split-free $f\le x$. Hence, the contribution to the left-hand side of \eqref{eq:psi-L} from the primes $q\in I_3$ is
\begin{equation}\label{eq:conI_3}
 \sum_{q\in I_3}\left(v_q(\psi(f))-v_q(L(f))\right)\log q<\frac{y\log_2y}{\log y}   
\end{equation}
for all but $o(S_{\alpha}(x))$ split-free $f\le x$. We now prove \eqref{eq:G(f)}. We have 
\begin{align*}
S_{\alpha}(x)^{-1}\sum_{f\le x}\alpha(f)G(f)&=S_{\alpha}(x)^{-1}\sum_{q\in I_3}\log q\sum_{k\ge2}k\sum_{\substack{f\le x\\\omega(f,\Qq(q))=k}}\alpha(f)\\
&\le S_{\alpha}(x)^{-1}\sum_{q\in I_3}\log q\sum_{k\ge2}k\sum_{p_1<\cdots<p_k\in\Qq(q)}\sum_{\substack{f\le x\\p_1\cdots p_k\mid f}}\alpha(f)\\
&\ll\sum_{q\in I_3}\log q\sum_{k\ge2}k\sum_{\substack{p_1\cdots p_k\le x\\p_1<\cdots<p_k\\\forall i,\,p_i\equiv-1\pmod*{q}}}\frac{1}{p_1\cdots p_k}\left(1-\frac{\log(p_1\cdots p_k)}{\log 3x}\right)^{-1/2}.
\end{align*}
By \eqref{eq:NP}, we have
\begin{align*}
&\hspace*{5.5mm}\sum_{q\in I_3}\log q\sum_{k\ge2}k\sum_{\substack{p_1\cdots p_k\le \sqrt{x}\\p_1<\cdots<p_k\\\forall i,\,p_i\equiv-1\pmod*{q}}}\frac{1}{p_1\cdots p_k}\left(1-\frac{\log(p_1\cdots p_k)}{\log 3x}\right)^{-1/2}\\
&\ll\sum_{q\in I_3}\log q\sum_{k\ge2}\frac{1}{(k-1)!}\left(\sum_{\substack{p\le x\\p\equiv-1\pmod*{q}}}\frac{1}{p}\right)^{k}\\
&=\sum_{q\in I_3}\log q\sum_{k\ge2}\frac{1}{(k-1)!}\left(\frac{y+O(\log q)}{\phi(q)}\right)^{k}\ll\frac{y}{\log y},
\end{align*}
where the last inequality follows from $y/\phi(q)\ll 1/\log y=o(1)$, as in \cite{EPS}.
We still have to estimate
\[\sum_{q\in I_3}\log q\sum_{k\ge2}k\sum_{\substack{\sqrt{x}<p_1\cdots p_k\le x\\p_1<\cdots<p_k\\\forall i,\,p_i\equiv-1\pmod*{q}}}\frac{1}{p_1\cdots p_k}\left(1-\frac{\log(p_1\cdots p_k)}{\log 3x}\right)^{-1/2}.\]
Put $m=\lceil y/\log y\rceil$. Stirling's formula implies that 
\[m!\gg m^{m+1/2}e^{-m}=e^{y+O(y/\log y)}\ge(\log x)^{1/2}.\]
Since $y/\phi(q)=o(1)$, we have
\begin{align*}
&\hspace*{5.5mm}\sum_{q\in I_3}\log q\sum_{k\ge m+1}k\sum_{\substack{\sqrt{x}<p_1\cdots p_k\le x\\p_1<\cdots<p_k\\\forall i,\,p_i\equiv-1\pmod*{q}}}\frac{1}{p_1\cdots p_k}\left(1-\frac{\log(p_1\cdots p_k)}{\log 3x}\right)^{-1/2}\\
&\ll(\log x)^{1/2}\sum_{q\in I_3}\log q\sum_{k\ge m+1}k\sum_{\substack{\sqrt{x}<p_1\cdots p_k\le x\\p_1<\cdots<p_k\\\forall i,\,p_i\equiv-1\pmod*{q}}}\frac{1}{p_1\cdots p_k}\\
&\le(\log x)^{1/2}\sum_{q\in I_3}\log q\sum_{k\ge m+1}\frac{1}{(k-1)!}\left(\sum_{\substack{p\le x\\p\equiv-1\pmod*{q}}}\frac{1}{p}\right)^{k}\\
&\ll\sum_{q\in I_3}\log q\sum_{k\ge m+1}\frac{m!}{(k-1)!}\left(\frac{y+O(\log y)}{\phi(q)}\right)^{k}\\
&\le e\sum_{q\in I_3}\log q\left(\frac{y+O(\log y)}{\phi(q)}\right)^2\\
&\ll y^2\sum_{q>y\log y}\frac{\log q}{q^2}\ll\frac{y}{\log y}.
\end{align*}
Thus, it remains to estimate
\begin{equation}\label{eq:I_3m}
\sum_{q\in I_3}\log q\sum_{k=2}^{m}k\sum_{\substack{\sqrt{x}<p_1\cdots p_k\le x\\p_1<\cdots<p_k\\\forall i,\,p_i\equiv-1\pmod*{q}}}\frac{1}{p_1\cdots p_k}\left(1-\frac{\log(p_1\cdots p_k)}{\log 3x}\right)^{-1/2}   
\end{equation}
Thanks to the constraint $k\le m$, we have $x^{1/2k}/q=x^{(1/2-o(1))/k}$. Besides, the constraints $p_1\cdots p_k>\sqrt{x}$ and $p_1<\cdots<p_k$ imply that $p_k>(p_1\cdots p_k)^{1/k}>x^{1/2k}$. Let $m_k\colonequals p_1\cdots p_{k-1}$, $a_k\colonequals \max\left(x^{1/2k},\sqrt{x}/m_k\right)$ and $b_k\colonequals x/m_k\in(x^{1/2k},x]$. Since
\[-\frac{d}{dt}\left(\frac{1}{t}\left(1-\frac{\log(m_kt)}{\log 3x}\right)^{-1/2}\right)\le \frac{1}{t^2}\left(1-\frac{\log(m_kt)}{\log 3x}\right)^{-1/2}\]

for any $t\in(0,b_k]$, it follows from Brun--Titchmarsh and partial summation that
\begin{align*}
&\hspace*{5.5mm}\sum_{\substack{a_k<p_k\le b_k\\p_k\equiv-1\pmod*{q}}}\frac{1}{p_k}\left(1-\frac{\log(p_1\cdots p_k)}{\log 3x}\right)^{-1/2}\\
&\ll\frac{(\log x)^{1/2}\pi(b_k;q,-1)}{b_k}-\int_{a_k}^{b_k}\pi(t;q,-1)\frac{d}{dt}\left(\frac{1}{t}\left(1-\frac{\log(m_kt)}{\log 3x}\right)^{-1/2}\right)\,dt\\
&\le\frac{(\log x)^{1/2}\pi(b_k;q,-1)}{b_k}+\int_{a_k}^{b_k}\frac{\pi(t;q,-1)}{t^2}\left(1-\frac{\log(m_kt)}{\log 3x}\right)^{-1/2}\,dt\\
&\ll\frac{k}{\phi(q)\sqrt{\log x}}+\frac{1}{\phi(q)}\int_{a_k}^{b_k}\frac{1}{t\log t}\left(1-\frac{\log(m_kt)}{\log 3x}\right)^{-1/2}\,dt\\
&\ll\frac{k}{\phi(q)\sqrt{\log x}}+\frac{k}{\phi(q)\log x}\int_{\sqrt{x}}^{x}\frac{1}{t}\left(1-\frac{\log t}{\log 3x}\right)^{-1/2}\,dt\\
&\ll\frac{k}{\phi(q)\sqrt{\log x}}+\frac{k}{\phi(q)}\int_{\frac{\log \sqrt{x}}{\log 3x}}^{\frac{\log x}{\log 3x}}(1-t)^{-1/2}\,dt\\
&\ll\frac{k}{\phi(q)}.
\end{align*}
Hence, we have
\begin{align*}
&\hspace*{5.5mm}\sum_{q\in I_3}\log q\sum_{k=2}^{m}k\sum_{\substack{\sqrt{x}<p_1\cdots p_k\le x\\p_1<\cdots<p_k\\\forall i,\,p_i\equiv-1\pmod*{q}}}\frac{1}{p_1\cdots p_k}\left(1-\frac{\log(p_1\cdots p_k)}{\log 3x}\right)^{-1/2}\\
&\ll\sum_{q\in I_3}\frac{\log q}{\phi(q)}\sum_{k\ge2}k^2\sum_{\substack{p_1<\cdots<p_{k-1}\le x\\\forall i,\,p_i\equiv-1\pmod*{q}}}\frac{1}{p_1\cdots p_{k-1}}\\
&\ll\sum_{q\in I_3}\frac{\log q}{\phi(q)}\sum_{k\ge2}\frac{k}{(k-2)!}\left(\frac{y+O(\log y)}{\phi(q)}\right)^{k-1}\\
&\ll \sum_{q\in I_3}\frac{\log q}{\phi(q)}\left(\frac{y+O(\log y)}{\phi(q)}\right)\ll\frac{1}{\log y}.
\end{align*}
Collecting the estimates above yields \eqref{eq:G(f)}.
\subsection{Contribution from $I_4$ and completion of the proof}\label{subsec:I_4}
Finally, we estimate the contribution to the left-hand side of \eqref{eq:psi-L} from the primes $q\in I_4$. We show that 
\begin{equation}\label{eq:conI_4}
 \sum_{q\in I_4}\left(v_q(\psi(f))-v_q(L(f))\right)\log q=0 
\end{equation}
for all but $o(S_{\alpha}(x))$ split-free $f\le x$. Since we may assume $q^2\nmid f$ for any $q\in I_4$, it suffices to show that all but $o(S_{\alpha}(x))$ split-free $f\le x$ are divisible by at most one prime in $\Qq(q)$ for any $q\in I_4$. This follows quickly from an averaging argument. We have
\[S_{\alpha}(x)^{-1}\sum_{q\in I_4}\sum_{p_1<p_2\in\Qq(q)}\sum_{\substack{f\le x\\p_1p_2\mid f}}\alpha(f)\ll\sum_{y^2<q\le\sqrt{x}+1}\sum_{\substack{p_1p_2\le x\\p_1<p_2\\p_1,\,p_2\equiv-1\pmod*{q}}}\frac{1}{p_1p_2}\left(1-\frac{\log p_1p_2}{\log 3x}\right)^{-1/2},\]
where the constraint $q\le\sqrt{x}+1$ holds because $q$ divides $p_1+1<\sqrt{x}+1$. It is easy to see that
\begin{align*}
\sum_{y^2<q\le\sqrt{x}+1}\sum_{\substack{p_1p_2\le x^{2/3}\\p_1<p_2\\p_1,\,p_2\equiv-1\pmod*{q}}}\frac{1}{p_1p_2}\left(1-\frac{\log p_1p_2}{\log 3x}\right)^{-1/2}&\ll\sum_{y^2<q\le\sqrt{x}+1}\sum_{\substack{p_1,\,p_2\le x\\p_1,\,p_2\equiv-1\pmod*{q}}}\frac{1}{p_1p_2}\\
&\ll y^2\sum_{q>y^2}\frac{1}{\phi(q)^2}\ll\frac{1}{\log y}
\end{align*}
and 
\begin{align*}
\sum_{\log x<q\le\sqrt{x}+1}\sum_{\substack{p_1p_2\le x\\p_1<p_2\\p_1,\,p_2\equiv-1\pmod*{q}}}\frac{1}{p_1p_2}\left(1-\frac{\log p_1p_2}{\log 3x}\right)^{-1/2}&\ll\sum_{\log x<q\le\sqrt{x}+1}\sum_{\substack{p_1,\,p_2\le x\\p_1,\,p_2\equiv-1\pmod*{q}}}\frac{(\log x)^{1/2}}{p_1p_2}\\
&\le(\log x)^{1/2}y^2\sum_{q>\log x}\frac{1}{\phi(q)^2}\\
&\ll \frac{y}{\sqrt{\log x}}.
\end{align*}
Hence, it remains to estimate the sum
\[\sum_{y^2<q\le\log x}\sum_{\substack{x^{2/3}<p_1p_2\le x\\p_1<p_2\\p_1,\,p_2\equiv-1\pmod*{q}}}\frac{1}{p_1p_2}\left(1-\frac{\log p_1p_2}{\log 3x}\right)^{-1/2}.\]
But the argument used to estimate \eqref{eq:I_3m} shows that this sum is 
\[\ll\sum_{y^2<q\le\log x}\frac{1}{\phi(q)}\sum_{\substack{p\le x\\p\equiv-1\pmod*{q}}}\frac{1}{p}\ll y\sum_{q>y^2}\frac{1}{\phi(q)^2}\ll\frac{1}{y\log y}.\]
Therefore, we have
\[S_{\alpha}(x)^{-1}\sum_{q\in I_4}\sum_{p_1<p_2\in\Qq(q)}\sum_{\substack{f\le x\\p_1p_2\mid f}}\alpha(f)\ll\frac{1}{\log y},\]
from which we verify that all but $o(S_{\alpha}(x))$ split-free $f\le x$ are divisible by at most one prime in $\Qq(q)$ for any $q\in I_4$.

Combining \eqref{eq:normalh}, \eqref{eq:conI_1}, \eqref{eq:conI_2c}, \eqref{eq:conI_3} and \eqref{eq:conI_4} completes the proof of \eqref{eq:psi-L} and hence that of Proposition \ref{prop:Lforder} for $L(f)$.

\section{Interlude: The radical of $\psi(f)/L(f)$}\label{sec:radicalpsiL} 
We devote this section to establishing the following result on the typical size of the radical of $\psi(f)/L(f)$ with $f$ split-free, which will be needed in the proof of Proposition \ref{prop:keythm2}.
\begin{prop}\label{prop:radpsiL} For almost all split-free $f$,
\[ \Rad \frac{\psi(f)}{L(f)} = (\log{f})^{\frac{1}{2} + O(1/\log_3 f)}. \]
\end{prop}
Our proof of Proposition \ref{prop:radpsiL} adapts that of \cite[Theorem 1]{pollack21} on the typical size of $\omega(\phi(n)/\lambda(n))$. As in Section \ref{sec:Lforder}, the complication comes mainly from estimation of weighted sums. Our setup is similar to that in \cite{pollack21}. Let 
\[W(f)\colonequals\log\Rad\frac{\psi(f)}{L(f)}=\sum_{p\mid\frac{\psi(f)}{L(f)}}\log p.\] 
Let $y=\log_2 x$, and define $\I\colonequals(y/\log y,y\log y]$, which is exactly same as the interval $I_2$ introduced in Section \ref{sec:Lforder}, and
\[\J_f\colonequals\{p\in \I\colon\exists\text{~distinct primes~}q_1,q_2\mid f\text{~with~}q_i\le x^{1/y}\text{~and~}q_i\equiv\chi(q_i)\pmod*{p}\text{~for~}i=1,2\}\]
for split-free $f\in\N$. Finally, we set
\[\widetilde{W}(f)\colonequals\sum_{p\in \J_f}\log p,\]
which we think of as an approximation to $W(f)$. It is not hard to see that Proposition \ref{prop:radpsiL} follows from the following estimates combined:
\begin{enumerate}[label=(\Roman*)]
\item $W(f)-\widetilde{W}(f)$ on average:
\[S_{\alpha}(x)^{-1}\sum_{f\le x}\alpha(f)\left(W(f)-\widetilde{W}(f)\right)\ll\frac{y}{\log y}.\]
\item The first moment of $\widetilde{W}$:
\[S_{\alpha}(x)^{-1}\sum_{f\le x}\alpha(f)\widetilde{W}(f)=\left(\frac{1}{2}+O\left(\frac{1}{\log y}\right)\right)y.\]
\item The second moment of $\widetilde{W}$:
\[S_{\alpha}(x)^{-1}\sum_{f\le x}\alpha(f)\widetilde{W}(f)^2=\left(\frac{1}{4}+O\left(\frac{1}{\log y}\right)\right)y^2.\] 
\end{enumerate}
We shall prove (I) and (II) but only sketch the proof of (III) while leaving the full details to the interested reader. Before embarking on the proof, we point out that the primes $p$ included in the definition of $W(f)$ satisfy either (i) $p^2\mid f$, or (ii) $f$ is divisible by two distinct primes $q_1,q_2$ with $q_i\equiv\chi(q_i)\pmod*{p}$ for $i=1,2$. The primes in set $\J_f$ hence belong to the second category. 
\subsection{$W(f)-\widetilde{W}(f)$ on average}\label{subsection:(I)}
We start by proving (I). Write $W(f)-\widetilde{W}(f)=T_0(f)+T_1(f)$, where
\begin{align*}
 T_0(f)&=\sum_{\substack{p\mid\frac{\psi(f)}{L(f)}\\p\le y/\log y}}\log p,\\ 
 T_1(f)&=\sum_{\substack{p\mid\frac{\psi(f)}{L(f)},\,p\notin \J_f\\p> y/\log y}}\log p.
\end{align*}
The Prime Number Theorem implies that $T_0(f)<2y/\log y$, so that
\[S_{\alpha}(x)^{-1}\sum_{f\le x}\alpha(f)T_0(f)<\frac{2y}{\log y}.\]
To estimate the mean value of $T_1$, we observe that the primes $p$ included in definition of $T_1(f)$ satisfy (i) $p^2\mid f$, or (ii) $p\in\I$ and $q\mid f$ for some prime $q\equiv\chi(q)\pmod*{p}$ with $q>x^{1/y}$, or (iii) $p> y\log y$ and there exist distinct primes $q_1,q_2\mid f$ with $q_i\equiv\chi(q_i)\pmod*{p}$ for $i=1,2$. The contribution from those primes $p$ satisfying (i) is
\[S_{\alpha}(x)^{-1}\sum_{f\le x}\alpha(f)\sum_{\substack{p>y/\log y\\p^2\mid f}}\log p=S_{\alpha}(x)^{-1}\sum_{p>y/\log y}\alpha(p^2)\log p\sum_{f\le x/p^2}\alpha(f)\ll\frac{\log y}{y},\]
by the same argument at the beginning of Section \ref{sec:Lforder}. Next, the contribution from those primes $p$ satisfying (ii) is
\begin{align*}
S_{\alpha}(x)^{-1}\sum_{f\le x}\alpha(f)\sum_{p\in\I}\log p\sum_{\substack{x^{1/y}<q\le x\\q\mid f\\q\equiv\chi(q)\pmod*{p}}}1&=
S_{\alpha}(x)^{-1}\sum_{p\in\I}\log p\sum_{\substack{x^{1/y}<q\le x\\q\equiv\chi(q)\pmod*{p}}}\sum_{\substack{f\le x\\q\mid f}}\alpha(f)\\
&\ll\sum_{p\in\I}\log p\sum_{\substack{x^{1/y}<q\le x\\q\equiv-1\pmod*{p}}}\frac{1}{q}\left(1-\frac{\log q}{\log 3x}\right)^{-1/2}.
\end{align*}
Since Brun--Tichmarsh and partial summation imply 
\begin{align*}
\sum_{\substack{x^{1/y}<q\le \sqrt{x}\\q\equiv-1\pmod*{p}}}\frac{1}{q}\left(1-\frac{\log q}{\log 3x}\right)^{-1/2}&\ll\sum_{\substack{x^{1/y}<q\le \sqrt{x}\\q\equiv-1\pmod*{p}}}\frac{1}{q}\ll\frac{\log y}{\phi(p)},\\
\sum_{\substack{\sqrt{x}<q\le x\\q\equiv-1\pmod*{p}}}\frac{1}{q}\left(1-\frac{\log q}{\log 3x}\right)^{-1/2}&\ll\frac{1}{\phi(p)},
\end{align*}
we find that the contribution from those primes $p$ satisfying (ii) is
\[\ll\log y\sum_{p\in\I}\frac{\log p}{\phi(p)}\ll(\log y)\log_2y.\]
Finally, the contribution from those primes $p$ satisfying (iii) is
\begin{align*}
&\hspace*{5.5mm}S_{\alpha}(x)^{-1}\sum_{f\le x}\alpha(f)\sum_{p>y\log y}\log p\sum_{\substack{q_1q_2\le x\\q_1q_2\mid f,\,q_1\ne q_2\\\forall i,\,q_i\equiv\chi(q_i)\pmod*{p}}}1&\\
&\ll\sum_{p>y\log y}\log p\sum_{\substack{q_1q_2\le x\\q_1<q_2\\q_1,q_2\equiv-1\pmod*{p}}}\frac{1}{q_1q_2}\left(1-\frac{\log q_1q_2}{\log 3x}\right)^{-1/2}.
\end{align*}
If $p>\log x$, then we have
\[\sum_{\substack{q_1q_2\le x\\q_1<q_2\\q_1,q_2\equiv-1\pmod*{p}}}\frac{1}{q_1q_2}\left(1-\frac{\log q_1q_2}{\log 3x}\right)^{-1/2}\ll(\log x)^{1/2}\left(\sum_{\substack{q\le x\\q\equiv-1\pmod*{p}}}\frac{1}{q}\right)^2\ll\frac{(\log x)^{1/2}y^2}{\phi(p)^2};\]
if $y\log y<p\le \log x$, then we have
\begin{align*}
\sum_{\substack{q_1q_2\le x^{2/3}\\q_1<q_2\\q_1,q_2\equiv-1\pmod*{p}}}\frac{1}{q_1q_2}\left(1-\frac{\log q_1q_2}{\log 3x}\right)^{-1/2}&\ll\left(\sum_{\substack{q\le x\\q\equiv-1\pmod*{p}}}\frac{1}{q}\right)^2\ll\frac{y^2}{\phi(p)^2},\\
\sum_{\substack{x^{2/3}<q_1q_2\le x\\q_1<q_2\\q_1,q_2\equiv-1\pmod*{p}}}\frac{1}{q_1q_2}\left(1-\frac{\log q_1q_2}{\log 3x}\right)^{-1/2}&\le\sum_{\substack{q_1\le\sqrt{x}\\q_1\equiv-1\pmod*{p}}}\frac{1}{q_1}\sum_{\substack{x^{2/3}/q_1<q_2\le x/q_1\\q_2\equiv-1\pmod*{p}}}\frac{1}{q_2}\left(1-\frac{\log q_1q_2}{\log 3x}\right)^{-1/2}\\
&\ll\frac{1}{\phi(p)}\sum_{\substack{q_1\le\sqrt{x}\\q_1\equiv-1\pmod*{p}}}\frac{1}{q_1}\ll\frac{y}{\phi(p)^2},
\end{align*}
where the estimate for the inner sum over $x^{2/3}/q_1<q_2\le x/q_1$ follows from the proof of \eqref{eq:I_3m} in Subsection \ref{subsec:I_3}. Hence, the contribution from those primes $p$ satisfying (iii) is
\[\ll y^2(\log x)^{1/2}\sum_{p>\log x}\frac{\log p}{\phi(p)^2}+y^2\sum_{p>y\log y}\frac{\log p}{\phi(p)^2}\ll\frac{y}{\log y}.\]
Gathering all the contributions to the mean value of $T_1$ above yields
\[S_{\alpha}(x)^{-1}\sum_{f\le x}\alpha(f)T_1(f)\ll \frac{y}{\log y}.\]
Combining this with the mean value estimate for $T_0$ proves (I).

\subsection{The first moment of $\widetilde{W}(f)$}\label{subsection:(II)}
We now turn to the proof of (II). Following \cite{pollack21}, let us write
\begin{equation}\label{eq:1stW}
\sum_{f\le x}\alpha(f)\widetilde{W}(f)=\sum_{p\in\I}\log p\sum_{\substack{f\le x\\p\in\J_f}}\alpha(f)=\sum_{p\in\I}(S_{\alpha}(x)-N_0(p)-N_1(p))  \log{p}
\end{equation}
where $N_i(p)$ is defined to be
\[\#\{f\le x\text{~split-free}\colon\exists\text{~exactly \emph{i} distinct primes~}q\mid f\text{~with~}q\le x^{1/y}\text{~and~}q\equiv\chi(q)\pmod*{p}\}\]
for $i=0,1$. By the Prime Number Theorem, we have 
\begin{equation}\label{eq:N_ast}
S_{\alpha}(x)^{-1}\sum_{p\in\I}S_{\alpha}(x)\log p=\int_{\I}1\,dt+O\left(\frac{y}{(\log y)^A}\right)    
\end{equation}
for any fixed $A>0$. 

To estimate the mean value of $N_0(p)\log p$, we apply the fundamental lemma of sieve theory to the set of primes  $q\le x^{1/y}$ with $q\notin\Pp_1$ and $q\equiv\chi(q)\pmod*{p}$, where $D=x^{u/y}$ and $u=\log y$, in exactly the same way as in Subsection \ref{subsec:I_2}, to obtain
\[S_{\alpha}(x)^{-1}N_0(p)=\left(1+O\left(u^{-u/2}\right)\right)\prod_{\substack{q\le x^{1/y}\\q\in\Pp_{-1}\\q\equiv-1\pmod*{p}}}\left(1-\frac{1}{q}\right)+O\left(\frac{\log y}{y}e^{\frac{y}{2p}}\right).\]
Since 
\[\sum_{\substack{q\le x^{1/y}\\q\in\Pp_{-1}\\q\equiv-1\pmod*{p}}}\frac{1}{q}=\frac{\phi(|\Delta|)}{2}\cdot\frac{\log\log x^{1/y}+O(\log p)}{\phi(|\Delta|)\phi(p)}=\frac{y}{2p}+O\left(\frac{(\log y)^2}{y}\right),\]
we have
\begin{align*}
S_{\alpha}(x)^{-1}N_0(p)&=\left(1+O\left(u^{-u/2}\right)\right)\exp\left(-\frac{y}{2p}+O\left(\sum_{q\equiv-1\pmod*{p}}\frac{1}{q^2}+\frac{(\log y)^2}{y}\right)\right)+O\left(\frac{\log y}{y}e^{\frac{y}{2p}}\right)\\
&=\left(1+O\left(u^{-u/2}\right)\right)\exp\left(-\frac{y}{2p}+O\left(\frac{1}{p\log p}+\frac{(\log y)^2}{y}\right)\right)+O\left(\frac{\log y}{y}e^{\frac{y}{2p}}\right)\\
&=\left(1+O\left(\frac{(\log y)^2}{y}\right)\right)e^{-\frac{y}{2p}}+O\left(\frac{\log y}{y}e^{\frac{y}{2p}}\right).
\end{align*}
It follows from \eqref{eq:sumI_2e^{y/2q}} that
\[S_{\alpha}(x)^{-1}\sum_{p\in\I}N_0(p)\log p=\left(1+O\left(\frac{(\log y)^2}{y}\right)\right)\sum_{p\in\I}e^{-\frac{y}{2p}}\log p+O\left(\sqrt{y}\right).\]
By partial summation and the Prime Number Theorem, we have
\[\sum_{p\in\I}e^{-\frac{y}{2p}}\log p=\int_{\I}e^{-\frac{y}{2t}}\,dt+O\left(\frac{y}{(\log y)^A}\right)\]
for any fixed $A>0$. Inserting this into the last equation above, we find that
\begin{equation}\label{eq:N_0}
S_{\alpha}(x)^{-1}\sum_{p\in\I}N_0(p)\log p=\int_{\I}e^{-\frac{y}{2t}}\,dt+O\left(\frac{y}{(\log y)^A}\right).
\end{equation}

Finally, we estimate the mean value of  $N_1(p)$. Note that $N_1(p)$ can be approximated by
\[\sum_{\substack{q\le x^{1/y}\\q\notin\Pp_1\\q\equiv\chi(q)\pmod*{p}}}\#\left\{m\le x/q\text{~split-free}\colon m\text{~not divisible by any~}q'\le x^{1/y}\text{~with~}q'\equiv\chi(q')\pmod*{p}\right\},\]
with an error being
\[\ll\#\left\{f\le x\text{~split-free}\colon f\text{~divisible by any~}q^2\text{~for some~}q>y/\log y\right\}.\]
This error is evidently
\[\ll\sum_{q>y/\log y}\sum_{\substack{f\le x\\ q^2\mid f}}\alpha(f)\ll\frac{S_{\alpha}(x)}{y},\]
as we have seen at the beginning of Section \ref{sec:Lforder}. For each $q$ included in the approximation of $N_1(p)$ above, our application of the fundamental lemma of sieve theory in Subsection \ref{subsec:I_2}, again with $D=x^{u/y}$ and $u=\log y$, implies that the corresponding term is 
\begin{align*}
&=\left(1+O\left(u^{-u/2}\right)\right)\frac{S_{\alpha}(x)}{q}\prod_{\substack{q'\le x^{1/y}\\q'\in\Pp_{-1}\\q'\equiv-1\pmod*{p}}}\left(1-\frac{1}{q}\right)+O\left(\frac{S_{\alpha}(x)\log y}{yq}e^{\frac{y}{2p}}\right)\\
&=\left(1+O\left(\frac{(\log y)^2}{y}\right)\right)\frac{S_{\alpha}(x)}{q}e^{-\frac{y}{2p}}+O\left(\frac{S_{\alpha}(x)\log y}{yq}e^{\frac{y}{2p}}\right),
\end{align*}
as shown in the preceding subsection. Summing over $q$, we see that $N_1(p)$ is
\begin{align*}
&=S_{\alpha}(x)\left(1+O\left(\frac{(\log y)^2}{y}\right)\right)\left(\frac{y}{2p}+O\left(\frac{(\log y)^2}{y}\right)\right)e^{-\frac{y}{2p}}+O\left(S_{\alpha}(x)\frac{e^{\frac{y}{2p}}}{p}\log y\right),\\
&=S_{\alpha}(x)\left(1+O\left(\frac{(\log y)^2}{y}\right)\right)\frac{y}{2p}e^{-\frac{y}{2p}}+O\left(S_{\alpha}(x)\frac{e^{\frac{y}{2p}}}{p}\log y\right).
\end{align*}
By the Prime Number Theorem and partial summation, we have
\begin{align*}
\sum_{p\in I}e^{\frac{y}{2p}}\frac{\log p}{p}&=\int_{\I}e^{\frac{y}{2t}}\,d\left(\log t+O\left(\frac{1}{\log t}\right)\right)\\
&=\int_{\I}e^{\frac{y}{2t}}\frac{dt}{t}+O\left(\frac{\sqrt{y}}{\log y}\right)\\
&=\int_{\frac{1}{2\log y}}^{\frac{\log y}{2}}\frac{e^{v}}{v}\,dv+O\left(\frac{\sqrt{y}}{\log y}\right)\\
&\ll x^{1/3}\int_{\frac{1}{2\log y}}^{\frac{\log y}{3}}\frac{dv}{v}+\frac{1}{\log y}\int_{\frac{\log y}{3}}^{\frac{\log y}{2}}e^{v}\,dv+O\left(\frac{\sqrt{y}}{\log y}\right)\\
&\ll\frac{\sqrt{y}}{\log y}.
\end{align*}
Thus, it follows, upon summing over $p\in\I$ and applying partial summation and the Prime Number Theorem, that
\begin{align}\label{eq:N_1}
S_{\alpha}(x)^{-1}\sum_{p\in\I}N_1(p)\log p&=\left(1+O\left(\frac{(\log y)^2}{y}\right)\right)\int_{\I}\frac{y}{2t}e^{-\frac{y}{2t}}\,dt+O\left(\frac{y}{(\log y)^A}\right)\nonumber\\
&=\int_{\I}\frac{y}{2t}e^{-\frac{y}{2t}}\,dt+O\left(\frac{y}{(\log y)^A}\right)
\end{align}
for any fixed $A>0$, since the integral is $\asymp y\log_2y$.

Inserting \eqref{eq:N_ast}, \eqref{eq:N_0} and \eqref{eq:N_1} into \eqref{eq:1stW} and observing that
\[\frac{d}{dt}\left(t-te^{-\frac{y}{2t}}\right)=1-e^{-\frac{y}{2t}}-\frac{y}{2t}e^{-\frac{y}{2t}},\]
we conclude, by taking $A=1$, that
\begin{align*}
S_{\alpha}(x)^{-1}\sum_{f\le x}\alpha(f)\widetilde{W}(f)&=\int_{\I}\left(1-e^{-\frac{y}{2t}}-\frac{y}{2t}e^{-\frac{y}{2t}}\right)\,dt+O\left(\frac{y}{\log y}\right)\\
&=y\log y\left(1-e^{-\frac{1}{2\log y}}\right)-\frac{y}{\log y}\left(1-\frac{1}{\sqrt{y}}\right)+O\left(\frac{y}{\log y}\right)\\
&=\left(\frac{1}{2}+O\left(\frac{1}{\log y}\right)\right)y,
\end{align*}
which is (II). It may be worth noting that one can obtain an asymptotic expansion of any length for the first weighted moment of $\widetilde{W}$ by taking $A$ suitably large but fixed.
\subsection{The second moment of $\widetilde{W}(f)$: a sketch}\label{subsection:(III)}
Finally, we outline the proof (III) and invite the reader to fill in the necessary details. We start by writing
\begin{equation}\label{eq:2ndWa}
S_{\alpha}(x)^{-1}\sum_{f\le x}\alpha(f)\widetilde{W}(f)^2=S_{\alpha}(x)^{-1}\sum_{p_1,p_2\in\I}\log p_1\log p_2\sum_{\substack{f\le x\\p_1,p_2\in\J_f}}\alpha(f).    
\end{equation}
For each pair $(i_1,i_2)\in\{0,1\}^2$, we denote by $N_{i_1,i_2}(p_1,p_2)$ the number of split-free $f\le x$ such that for every $k\in\{1,2\}$, there are exactly $i_k$ distinct primes $q_k\mid f$ with $q_k\le x^{1/y}$ and $q_k\equiv\chi(q_k)\pmod*{p_k}$. Furthermore, let
\[M_{i_1,i_2}\colonequals\int_{\I}\left(\frac{y}{2t}\right)^{i_1}e^{-i_2\frac{y}{2t}}\,dt.\]
By the inclusion-exclusion principle, we have
\begin{equation}\label{eq:2ndWb}
\sum_{\substack{f\le x\\p_1,p_2\in\J_f}}\alpha(f)=S_{\alpha}(x)-\sum_{k\in\{1,2\}}\sum_{i_k\in\{0,1\}}N_{i_k}(p_k)+\sum_{(i_1,i_2)\in\{0,1\}^2}N_{i_1,i_2}(p_1,p_2).   
\end{equation}
The contribution to \eqref{eq:2ndWa} from the term $S_{\alpha}(x)$ in \eqref{eq:2ndWb} is obviously
\[\left(\sum_{p\in\I}\log p\right)^2=M_{0,0}^2+O\left(\frac{y^2}{(\log y)^A}\right)\]
for any fixed $A>0$. In addition, the contribution to \eqref{eq:2ndWa} from the double sum in \eqref{eq:2ndWb} is 
\[2\left(\sum_{p\in\I}\log p\right)\left(S_{\alpha}(x)^{-1}\sum_{p\in\I}N_0(p)\log p+S_{\alpha}(x)^{-1}\sum_{p\in\I}N_1(p)\log p\right),\]
which, according to \eqref{eq:N_0} and \eqref{eq:N_1}, is equal to
\[2M_{0,0}\left(M_{0,1}+M_{1,1}\right)+O\left(\frac{y^2}{(\log y)^A}\right)\]
for any fixed $A>0$. Hence, it remains to estimate the contribution to \eqref{eq:2ndWa} from the last sum in \eqref{eq:2ndWb}. Note first that the contribution to \eqref{eq:2ndWa} from the diagonal terms with $p_1=p_2$ in the last sum in \eqref{eq:2ndWb} is
\[S_{\alpha}(x)^{-1}\sum_{p\in\I}N_0(p)(\log p)^2+S_{\alpha}(x)^{-1}\sum_{p\in\I}N_1(p)(\log p)^2\ll\frac{y\log_2y}{(\log y)^A}\]
for any fixed $A>0$, which is negligible. Thus, it is sufficient to estimate the contribution to \eqref{eq:2ndWa} from the off-diagonal terms with $p_1\ne p_2$. We claim that
\begin{equation}\label{eq:N_{i_1,i_2}}
S_{\alpha}(x)^{-1}\sum_{p_1\ne p_2\in\I}N_{i_1,i_2}(p_1,p_2)\log p_1\log p_2=M_{i_1,1}M_{i_2,1}+O\left(\frac{y^2\log_2y}{(\log y)^2}\right)   
\end{equation}
for each pair $(i_1,i_2)\in\{0,1\}^2$. Inserting all the estimates above into \eqref{eq:2ndWa} yields
\[S_{\alpha}(x)^{-1}\sum_{f\le x}\alpha(f)\widetilde{W}(f)^2=\left(M_{0,0}-M_{0,1}-M_{1,1}\right)^2+O\left(\frac{y^2\log_2y}{(\log y)^2}\right)=\left(\frac{1}{4}+O\left(\frac{1}{\log y}\right)\right)y^2,\]
which is (III). Again, it is possible to make $O(y^2/\log y)$ explicit and have $O\left(y^2\log_2y/(\log y)^2\right)$ as the error instead.

Taking the case $(i_1,i_2)=(1,1)$ for example, we now illustrate briefly how \eqref{eq:N_{i_1,i_2}} can be derived by adapting the argument in \cite{pollack21}. For any $p_1\ne p_2\in\I$, the quantity $S_{\alpha}(x)^{-1}N_{1,1}(p_1,p_2)$ includes particularly the contribution from those split-free $f\le x$ of the form $f=q_1q_2m$, where $q_1,q_2\le x^{1/y}$ are distinct primes not in $\Pp_1$, satisfying $q_k\equiv\chi(q_k)\pmod*{p_k}$ for all $k=1,2$ and $q_k\not\equiv\chi(q_k)\pmod*{p_l}$ for all $(k,l)=(1,2),(2,1)$, and where $m\le x/q_1q_2$ is split-free and free of prime factors $q\le x^{1/y}$ with $q\equiv\chi(q)\pmod*{p_k}$ for all $k=1,2$. The rest of the split-free $f\le x$ that contribute to $S_{\alpha}(x)^{-1}N_{1,1}(p_1,p_2)$ satisfy either $q^2\mid f$ for some prime $q\ge y/\log y-1>y/2\log y$ or $q\mid f$ for some prime $q$ with $q\equiv\chi(q)\pmod*{p_1p_2}$. Thus, the contribution to $S_{\alpha}(x)^{-1}N_{1,1}(p_1,p_2)$ from these residual split-free $f\le x$ does not exceed
\[S_{\alpha}(x)^{-1}\sum_{y/2\log y< q\le \sqrt{x}}\sum_{\substack{f\le x\\q^2\mid f}}\alpha(f)+S_{\alpha}(x)^{-1}\sum_{\substack{q\le x\\q\equiv\chi(q)\pmod*{p_1p_2}}}\sum_{\substack{f\le x\\q\mid f}}\alpha(f)\ll\frac{1}{y}+\frac{y}{p_1p_2},\]
where the sums over $q$ can be easily estimated by dividing the ranges of $q$ as in Section \ref{sec:Lforder}. Summing on $p_1,p_2\in\I$, we see that the contribution to the left-hand side of \eqref{eq:N_{i_1,i_2}} from these residual split-free $f\le x$ is $O(y(\log y)^2)$, which is negligible compared to the error term in \eqref{eq:N_{i_1,i_2}}. Now we turn to those split-free $f\le x$ of the form $f=q_1q_2m$ with the properties described above. Given $p_1,p_2,q_1,q_2$, the count of these $f$, when divided by $S_{\alpha}(x)$, is
\[\left(1+O\left(u^{-u/2}\right)\right)\frac{1}{q_1q_2}\prod_{\substack{q\le x^{1/y}\\q\notin\Pp_{1}\\\exists k=1,2,\,q\equiv\chi(q)\pmod*{p_k}}}\left(1-\frac{1}{p}\right)+O\left(\frac{\log_2 y}{yq_1q_2}e^{\frac{y}{2p_1}+\frac{y}{2p_2}}\right),\]
by the fundamental lemma of sieve theory with $D=x^{u/y}$ and $u=\log_2y$, as in Subsection \ref{subsec:I_2}. The contribution to the left-hand side of \eqref{eq:N_{i_1,i_2}} from the second error term above is 
\begin{align*}
&\ll\frac{\log_2y}{y}\sum_{p_1,p_2\in\I}e^{\frac{y}{2p_1}+\frac{y}{2p_2}}\log p_1\log p_2\sum_{\substack{q_1,q_2\le x^{1/y}\\\forall k=1,2,\,q_k\equiv-1\pmod*{p_k}}}\frac{1}{q_1q_2}\\
&\ll y\log_2y\sum_{p_1,p_2\in\I}e^{\frac{y}{2p_1}+\frac{y}{2p_2}}\frac{\log p_1\log p_2}{p_1p_2}\\
&\ll y\log_2y\left(\frac{\sqrt{y}}{\log y}\right)^2\ll\frac{y^2\log_2y}{(\log y)^2},
\end{align*}
which matches the error in \eqref{eq:N_{i_1,i_2}}. The main contribution to the left-hand side of \eqref{eq:N_{i_1,i_2}} from those split-free $f\le x$ of the form $f=q_1q_2m$ with the properties described above is 
\[\left(1+O\left(u^{-u/2}\right)\right)\sum_{p_1,p_2\in\I}\log p_1\log p_2\sum_{\substack{q_1,q_2\le x^{1/y}\\q_1,q_2\notin\Pp_{1}\\\forall k\ne l\in\{1,2\},\,q_k\equiv\chi(q_k)\pmod*{p_k}\\q_k\not\equiv\chi(q_k)\pmod*{p_l}}}\frac{1}{q_1q_2}\prod_{\substack{q\le x^{1/y}\\q\notin\Pp_{1}\\\exists k=1,2,\,q\equiv\chi(q)\pmod*{p_k}}}\left(1-\frac{1}{p}\right).\]
Now an argument analogous to the one displayed on \cite[p.\ 214]{pollack21} shows that
\[\sum_{\substack{q_1,q_2\le x^{1/y}\\q_1,q_2\notin\Pp_{1}\\\forall k\ne l\in\{1,2\},\,q_k\equiv\chi(q_k)\pmod*{p_k}\\q_k\not\equiv\chi(q_k)\pmod*{p_l}}}\frac{1}{q_1q_2}\prod_{\substack{q\le x^{1/y}\\q\notin\Pp_{1}\\\exists k=1,2,\,q\equiv\chi(q)\pmod*{p_k}}}\left(1-\frac{1}{p}\right)=\left(1+O\left(\frac{(\log y)^2}{y}\right)\right)\frac{y^2}{4p_1p_2}e^{-\frac{y}{2p_1}-\frac{y}{2p_2}}.\]
Summing over $p_1\ne p_2\in\I$ with the weights $\log p_1\log p_2$ attached, we find that the main contribution to the left-hand side of \eqref{eq:N_{i_1,i_2}} from those split-free $f\le x$ of the form $f=q_1q_2m$ is 
\begin{align*}
&=\left(1+O\left(u^{-u/2}\right)\right)\left(y\sum_{p\in\I}e^{-\frac{y}{2p}}\frac{\log p}{2p}\right)^2+O\left(y^2\sum_{p\in \I}\frac{(\log p)^2}{p^2}e^{-\frac{y}{p}}\right)\\
&=\left(1+O\left(u^{-u/2}\right)\right)\left(\int_{\I}\frac{y}{2t}e^{-\frac{y}{2t}}\,dt+O\left(\frac{y}{(\log y)^{A+1}}\right)\right)^2+O\left(y(\log y)^2\right)\\
&=M_{1,1}^2+O\left(\frac{y^2}{(\log y)^{A}}\right)
\end{align*}
for any fixed $A>0$. Gathering the contributions above verifies \eqref{eq:N_{i_1,i_2}} in the case $(i_1,i_2)=(1,1)$.

\section{Application of the algebro-analytic machine: GRH and the unit index corresponding to an inert prime}\label{sec:algebraicnonsense} For the rest of the paper, we assume that $K= \Q(\sqrt{D})$, where $D>1$ is squarefree. We write $\epsilon$ for the fundamental unit of $\Oo_K$ and $\sigma_0$ for the nontrivial element of $\Gal(K/\Q)$. All number fields appearing below are viewed as subfields of $\C$, and odd order roots of real numbers are understood as taking their real values.

Recall from \S\ref{sec:bigpicture} that $\ell(f)$ denotes the order of (the image of) $\epsilon$ in the group $\PreCl(\Oo_f) = (\Oo_K/f\Oo_K)^{\times}/\langle \text{images of integers prime to $f$}\rangle$. It will be important for the proof of Proposition \ref{prop:keythm2} to understand the distribution of the numbers $\ell(p)$, as $p$ varies over primes inert in $K$. Certainly $\ell(p)$ divides $\#\PreCl(\Oo_p) = \psi(p) = p+1$ for each such $p$. The main result of this section is a GRH-conditional estimate for how often $\ell(p) \mid \frac{p+1}{q}$ for a given odd prime $q$ (Proposition \ref{prop:hooleychen}).

The method of proof is essentially that introduced by Hooley \cite{hooley67} to study Artin's primitive root conjecture under the assumption of GRH. Of course, our setting is a bit different than Hooley's, as the arithmetic is taking place over a real quadratic field instead of $\Q$. Happily, we can quickly deduce what we need from work of Chen \cite{chen02}. (Closely related arguments can be found in papers of Roskam \cite{roskam00}, Kataoka \cite{kataoka03}, and Pollack \cite[\S2]{pollackm}.)

Let $q$ be an odd prime. We will relate the condition that $\ell(p)\mid \frac{p+1}{q}$ to the splitting behavior of $p$ in the number field
\[ E_q := \Q(\zeta_q, \sqrt[q]{\epsilon}). \] (As usual, $\zeta_m:=\exp(2\pi \mathrm{i}/m)$.) Note that $E_q$ contains $\Q(\epsilon) = K$, so that $E_q = K(\zeta_q,\sqrt[q]{\epsilon})$.

To bring our setup into alignment with Chen's, we need to write $E_q$ in a slightly different way. Let $s=1$ if $\Nm_{K/\Q}(\epsilon)=1$ and let $s=2$ if $\Nm_{K/\Q}(\epsilon) = -1$. Put $\eta = \epsilon^s$. Then $\Nm_{K/\Q}(\eta)=1$ (in fact, $\eta$ generates the group of norm $1$ units). Since $q$ is an odd prime while $s=1$ or $s=2$, we have that $E_q = K(\zeta_q, \sqrt[q]{\epsilon}) = K(\zeta_q, \sqrt[q]{\eta})$.

The next three lemmas are special cases of Chen's results in \cite{chen02}. After all three results have been stated, we say a few words about how they may be deduced from \cite{chen02}.

\begin{lem}\label{lem:Eqdeg} Let $q$ be an odd prime. If $q\nmid D$, then $[E_q:\Q] = 2q\phi(q)$,
\end{lem}

Set $\tilde{\eta} = \sigma_0(\eta)$. Then $\eta \tilde{\eta} = 1$ and $E_q$ contains a $q$th root of $\tilde{\eta}$, namely $1/\sqrt[q]{\eta}$. It follows that 
\[ E_q= K(\zeta_q, \sqrt[q]{\eta}) = K(\zeta_q, \sqrt[q]{\eta}, \sqrt[q]{\tilde{\eta}}) =  \Q(\zeta_q, \sqrt[q]{\eta}, \sqrt[q]{\tilde{\eta}}) ,\]
rendering apparent that $E_q$ is the splitting field over $\Q$ of 
\begin{align}\notag F_q(X):&= (X^q-\eta)(X^q-\tilde{\eta})\\
&= X^{2q} - \mathrm{Tr}_{K/\Q}(\eta) X^q + 1 \in \Z[X].\label{eq:Fqdef} \end{align}
Hence, $E_q$ is a Galois extension of $\Q$.

Let $\tau$ represent complex conjugation. Define $\Cc_q^{-} \subset \Gal(E_q/\Q)$ by  
\[ \Cc_q^{-} = \{\sigma \in \Gal(E_q/\Q): \sigma|_{K} = \sigma_0, \sigma|_{\Q(\zeta_q)} = \tau|_{\Q(\zeta_q)}, \sigma^2=\textrm{id}\}.\]
It is straightforward to check that $\Cc_q^{-}$ is a conjugation-stable subset of $\Gal(E_q/\Q)$. 

\begin{lem}\label{lem:chensingleton} Let $q$ be an odd prime not dividing $D$. Then $\#\Cc_q^{-}=1$.  
\end{lem} 

Write $\eta = \frac{1}{2}(u + v\sqrt{D})$ with integers $u$ and $v$.

\begin{lem}\label{lem:chenfrobenius} Suppose that $p$ is an odd prime not dividing $v$. Let $q$ be an odd prime. Then 
\[ p\text{ is inert in $K$},~p \equiv -1\pmod*{q},\text{ and } \eta^{\frac{p+1}{q}}\equiv 1\pmod*{p\Oo_K} \Longleftrightarrow \Frob_{E_q/\Q,p} \in \Cc_q^{-}. \]
\end{lem}

\begin{proof}[Deduction of Lemmas \ref{lem:Eqdeg}--\ref{lem:chenfrobenius} from Chen's work] Chen's setup in \cite{chen02} is as follows: $K_0$ is a quadratic field (not necessarily real) with nontrivial automorphism $\sigma_0$. The element $\alpha\in K_0^{\times}$ has norm $1$ but is not a root of unity. For each prime $q$, the field $E_q$ is defined as $K_0(\zeta_q, \sqrt[q]{\alpha})$, and when $n$ is squarefree, $E_n$ is defined as the compositum of the $E_q$ for primes $q$ dividing $n$. 

Under these assumptions, Chen determines the degrees of the extensions $E_n/\Q$ in Lemma 1.6 of \cite{chen02}, while Lemma 1.7 of \cite{chen02} describes the sizes of the sets
\[ \Cc_n^{-}:= \{\sigma \in \Gal(E_n/\Q): \sigma|_{K} = \sigma_0, \sigma|_{\Q(\zeta_n)} = \tau|_{\Q(\zeta_n)}, \sigma^2=\textrm{id},\text{ and } \sigma|_{K_0(\sqrt{\alpha}+1/\sqrt{\alpha})}=\textrm{id}\text{ if } 2\mid n\}.\] 
For brevity, we quote only the relevant cases: Suppose $n$ is odd and squarefree and that $K_0$ is real. Let $s_0$ be the largest positive integer for which $\alpha \in (K_0^{\times})^{s_0}$, and let $n_1 = n/\gcd(n,s_0)$. Then Chen's Lemma 1.6 asserts that
\[ [E_n:\Q] = \frac{2 n_1  \phi(n)}{[K_0 \cap \Q(\zeta_n):\Q]}, \] while her Lemma 1.7 claims that $\#\Cc_n^{-} = 1$ unless $K \subset \Q(\zeta_n)$.

Our Lemmas \ref{lem:Eqdeg} and \ref{lem:chensingleton} follow from these two results (respectively) upon taking $K_0=K$, $\alpha=\eta$, and $n=q$. Here we note that Chen's $s_0$ coincides with our integer $s$, and that
\[ \Q \subset K\cap \Q(\zeta_q) \subset \Q(\zeta_{4|D|}) \cap \Q(\zeta_q) = \Q(\zeta_{\gcd(4|D|,q)}) = \Q,\]
since we are assuming $q$ is an odd prime not dividing $D$.

What about Lemma \ref{lem:chenfrobenius}? Write $K_0 = \Q(\sqrt{D_0})$ with $D_0$ squarefree, and express $\alpha$ as $u_0+v_0\sqrt{D_0}$ with $u_0, v_0 \in \Q$. In Lemma 1.4(a), Chen assumes that $p$ is an odd prime inert in $K$ for which $v_p(v_0) = 0$. She then shows that for each odd prime $q$,
\begin{equation}\label{eq:chenimplication} p \equiv -1\pmod*{q} \text{ and } \eta^{\frac{p+1}{q}}\equiv 1\pmod*{p\Oo_K} \Longleftrightarrow \Frob_{E_q/\Q,p} \in \Cc_q^{-}. \end{equation}
Specializing to $K=K_0$ (so that $D=D_0$) and $\alpha=\eta$ (so that $u_0=\frac{1}{2}u, v_0 = \frac{1}{2}v$), the forward implication in \eqref{eq:chenimplication} yields the forward implication of Lemma \ref{lem:chenfrobenius}. One can also deduce the backward implication of Lemma \ref{lem:chenfrobenius} from the backward implication in \eqref{eq:chenimplication}; it suffices to observe that $\Frob_{E_q/\Q,p} \in \Cc_q^{-}$ implies that $p$ is inert in $K_0$. Indeed, if $\Frob_{E_q/\Q,p} \in \Cc_q^{-}$, then $\Frob_{K_0/\Q,p} = \Frob_{E_q/\Q,p}|_{K_0} = \sigma_0$.
\end{proof}

The following conditional version of the Chebotarev density theorem has been extracted from Serre's paper \cite{serre81} (take $K=\Q$ in Serre's eq.\ ($20_{{\rm R}}$)).

\begin{thmcheb} Let $L/\Q$ be a Galois extension, and let $\Cc$ be a subset of $\Gal(L/\Q)$ stable under conjugation. For all $y\ge 2$,
\[ \#\{\text{primes }p\le y: \Frob_{L/\Q,\,p} \subset \Cc\} = \frac{|\Cc|}{[L:\Q]} \Li(y) + O\left((\#\Cc)y^{1/2} \log\bigg([L:\Q]y\prod_{\ell \mid \Delta_L}\ell\bigg)\right).  \]
Here the implied constant is absolute.
\end{thmcheb}

\begin{lem}[assuming GRH]\label{lem:chebapp} Let $q$ be an odd prime not dividing $D$. For every real number $y\ge 2$, the number of primes $p\le y$ inert in $K$ for which $p\equiv -1\pmod{q}$ and $\eta^{(p+1)/q}\equiv 1\pmod{p\Oo_K}$ is
\[ \frac{\Li(y)}{2q\phi(q)} + O\left(y^{1/2}\log(qy)\right). \]
Here and below, implied constants may depend on $K$.
\end{lem}

\begin{proof} By adjusting the $O$-constant, we may ignore any set of primes $p$ of cardinality $O_K(1)$. Bearing in mind the results of Lemmas \ref{lem:chensingleton} and \ref{lem:chenfrobenius}, Lemma \ref{lem:chebapp}  follows from the GRH-dependent Chebotarev density theorem upon taking $L= E_q$ and $\Cc=\Cc_q^{-}$, once we show that the error term of the theorem is subsumed by our $O$-expression. 

To that end, we argue that a rational prime $p$ is unramified in $E_q$ whenever $p\nmid 2  q D v_0$. Observe that $p$ is unramified in $E_q$ as long as $F_q(X)$, as defined in \eqref{eq:Fqdef}, has no multiple roots mod $p$. We test this condition by checking for common roots of $F_q$ and $F_q'$. If $p\nmid 2q$, then each root $\rho \in \overline{\mathbb{F}}_p$ of $F_q'(X)$ has $2\rho^q = \Tr_{K/\Q}(\eta)$. But then 
\begin{align*} 4 F_q(\rho) &= 4 \rho^{2q} - 4  \Tr_{K/\Q}(\eta) \rho^q + 4 \\
&= 4 \Nm_{K/\Q}(\eta) - \Tr_{K/\Q}(\eta)^2 \\&= -D v^2,
\end{align*}
which is nonvanishing mod $p$ under our assumption that $p\nmid 2qDv$.

The primes that ramify in $E_q$ are precisely those dividing $\Delta_{E_q}$. Hence,
\[ \prod_{\ell \mid \Delta_{E_q}} \ell \le 2q D v,\quad\text{and}\quad [E_q:\Q] y \prod_{\ell \mid \Delta_{E_q}}\ell \le 4q^3 Dv y,    \]
using Lemma \ref{lem:Eqdeg} for the last inequality. Therefore, $$\log\big([E_q:\Q] y \prod_{\ell \mid \Delta_{E_q}}\ell\big)  \ll \log(qy),$$
and $(\#\Cc) y^{1/2} \log\big([E_q:\Q] y \prod_{\ell \mid \Delta_{E_q}}\ell\big) \ll y^{1/2} \log(qy)$.
\end{proof}

We now come to the main result of this section.

\begin{prop}[assuming GRH]\label{prop:hooleychen} Let $q$ be an odd prime, and let $y\ge 2$. The number of primes $p \le y$ for which $p$ is inert in $K$, $p\equiv -1\pmod{q}$, and $\ell(p) \mid \frac{p+1}{q}$ is
\[ \ll \frac{\Li(y)}{q^2} + y^{1/2}\log(qy). \]
\end{prop}

\begin{proof} For the finitely many odd primes $q$ dividing $D$, the result is trivial if we choose a sufficiently large implied constant (depending on $K$). So we may suppose that $q\nmid D$. 

We now argue that if $p$ is inert in $K$ and $p\equiv -1\pmod{q}$, then 
\[ \ell(p) \mid \frac{p+1}{q} \Longleftrightarrow \eta^{\frac{p+1}{q}}\equiv 1\pmod{p\Oo_K}. \]
Once this is proved, Proposition \ref{prop:hooleychen} follows from Lemma \ref{lem:chebapp}.

If $\eta^{(p+1)/q}\equiv 1\pmod{p\Oo_K}$, then $\epsilon^{(p+1)/q} \equiv \pm 1\pmod{p\Oo_K}$. Thus, $\epsilon^{(p+1)/q} \in \Oo_p$, and $\ell(p) \mid \frac{p+1}{q}$.  Conversely, if $\ell(p) \mid \frac{p+1}{q}$, then $\epsilon^{(p+1)/q} \equiv h\pmod{p\Oo_K}$ for some rational integer $h$ coprime to $p$. As $p$ is inert in $K$, the $p$th power map modulo $p\Oo_K$ is induced by $\sigma_0$, and 
\[ (\eta^{(p+1)/q})^q \equiv \eta \cdot \eta^{p} \equiv \Nm_{K/\Q}(\eta) \equiv 1\pmod{p\Oo_K}. \]
Hence, if $s$ denotes the order of $\eta^{(p+1)/q}$ in $(\Oo_K/p\Oo_K)^{\times}$, then $s \mid q$. On the other hand, $s$ coincides with the order of $h$ in $(\Z/p\Z)^{\times}$, and so $s\mid p-1$. Using that $q\mid p+1$, we conclude that $s$ divides
\[ \gcd(q,p-1) = \gcd(q,(p+1)-2) = \gcd(q,-2)=1. \]
Thus, $s=1$ and $\eta^{(p+1)/q} \equiv  1\pmod{p\Oo_K}$. 
\end{proof}

We conclude this section with a bound on the number of small values of $\ell(p)$.

\begin{lem}\label{lem:smallorderbound} For each $y\ge 1$, the number of inert primes $p$ for which $\ell(p) \le y$ is $O(y^2)$.
\end{lem}

As a consequence of Lemma \ref{lem:smallorderbound}, the number of inert primes $p\le t$ with $\ell(p) \le 2p^{1/2}/\log{p}$ is $O(t/(\log{t})^2)$, for all $t\ge 2$. This will be needed in the arguments of Section \ref{sec:completionkeythm2}.

\begin{proof} Let $p$ be a prime inert in $K$, and let $s$ be the order of $\epsilon^{\ell(p)}$ in $(\Oo_K/p\Oo_K)^{\times}$. Since $\epsilon^{\ell(p)}$ is congruent to a rational integer unit mod $p$, we have that $s\mid p-1$. Noting that $\epsilon^{p+1}\equiv \epsilon\cdot \epsilon^p\equiv \Nm_{K/\Q}(\epsilon) \equiv \pm 1 \pmod{p\Oo_K}$, we also have that $s\mid 2(p+1)$. Thus, $s\mid \gcd(p-1,2(p+1))\mid 4$, and $\epsilon^{4\ell(p)}=1$ in $\Oo_K/p\Oo_K$.

Now assume that $\ell(p)\le y$. Then $\epsilon^{j}\equiv 1\pmod{p\Oo_K}$ for some positive integer $j\le 4y$, and 
\[ p\text{ divides } \prod_{1 \le j \le 4y} \Nm(\epsilon^j-1). \]
As $\Nm(\epsilon^j-1) = \exp(O(j))$, the product on $j$ has size $\exp(O(y^2))$ and so is divisible by $O(y^2)$ distinct primes.
\end{proof}

\section{Proof of Proposition \ref{prop:keythm2}}\label{sec:completionkeythm2}

We now turn to the proof of Proposition \ref{prop:keythm2}. Our arguments borrow heavily from those of Erd\H{o}s--Pomerance--Schmutz and Li--Pomerance \cite{LP03}. Let
\[ z_f := \frac{\log_2{f}}{\log_4 f}. \]

The next lemma is proved by an averaging argument analogous to one appearing on p.\ 368 of \cite{EPS}.

\begin{lem}\label{lem:eps} For almost all split-free $f$, the $z_f$-smooth part of $L(f)$ is of size $(\log{f})^{O(1/\log_4 f)}$.
\end{lem}

\begin{proof} The proof bears a great resemblance to that of \eqref{eq:I_1b} in Subsection \ref{subsec:I_1}. For any $n\in\N$ and $y,w\ge1$, denote by $S(n,y)$ the $y$-smooth part of $n$ and put
\[S(n,y,w)\colonequals\prod_{\substack{p^k\parallel n\\p\le y\\p^{k}>w}}p^k.\]
With $z=\log_2 x/\log_4 x$ and $w=(\log_2 x)\log_4x$, we have
\begin{align*}
S_{\alpha}(x)^{-1}\sum_{f\le x}\alpha(f)\log S(L(f),z_f,w)&=S_{\alpha}(x)^{-1}\sum_{f\le x}\alpha(f)\sum_{\substack{p^k\parallel L(f)\\ p\le z_f,\,p^k>w}}\log p^k\\ 
&\le S_{\alpha}(x)^{-1}\sum_{\substack{p\le z\\p^k>w}}\log p^k\sum_{\substack{f\le x\\p^k\mid L(f)}}\alpha(f).
\end{align*}
Observe that if $p^k\mid L(f)$, then $p^k\mid f$ with $p$ ramified, or $p^{k+1}\mid f$ with $p$ unramified, or $f$ has a prime factor $q$ with $q\equiv\chi(q)\pmod{p^k}$. The contribution from the first case is
\begin{align*}
&\le S_{\alpha}(x)^{-1}\sum_{\substack{p\mid\Delta\\w<p^{k}\le x}}\log p^k\sum_{f\le x/p^{k}}\alpha(f)\\
&\ll\sum_{\substack{p\mid\Delta\\w<p^{k}\le x}}\frac{\log p^k}{p^{k}}\left(1-\frac{\log p^{k}}{\log 3x}\right)^{-1/2}\\
&\ll\sum_{\substack{p\mid\Delta\\w<p^{k}\le \sqrt{x}}}\frac{\log p^k}{p^{k}}+(\log x)^{1/2}\sum_{\substack{p\mid\Delta\\p^{k}>\sqrt{x}}}\frac{\log p^k}{p^{k}}\\
&\ll\frac{\log w}{w}+\frac{(\log x)^{3/2}}{\sqrt{x}}\ll\frac{\log w}{w},
\end{align*}
Analogously, the second case contributes an amount 
\begin{align*}
&\le S_{\alpha}(x)^{-1}\sum_{\substack{p\le z\\w<p^{k}\le x/p}}\log p^k\sum_{f\le x/p^{k+1}}\alpha(f)\\
&\ll\sum_{\substack{p\le z\\w<p^{k}\le x/p}}\frac{\log p^k}{p^{k+1}}\left(1-\frac{\log p^{k+1}}{\log 3x}\right)^{-1/2}\\
&\ll\sum_{\substack{p\le z\\w<p^{k}\le \sqrt{x}/p}}\frac{\log p^k}{p^{k+1}}+(\log x)^{1/2}\sum_{\substack{p\le z\\p^{k+1}>\sqrt{x}}}\frac{\log p^k}{p^{k+1}}\\
&\ll\frac{\log w}{w}\sum_{p\le z}\frac{1}{p}+\frac{(\log x)^{3/2}}{\sqrt{x}}\sum_{p\le z}1\ll\frac{(\log w)\log_2z}{w}.
\end{align*}
Finally, the last case contributes an amount 
\begin{align*}
&\le S_{\alpha}(x)^{-1}\sum_{\substack{p\le z\\p^k>w}}\log p^k\sum_{\substack{q\le x\\q\equiv\chi(q)\pmod*{p^k}}}\alpha(q)\sum_{f\le x/q}\alpha(f)\\
&\ll\sum_{\substack{p\le z\\p^k>w}}\log p^k\sum_{\substack{q\le x\\q\equiv\chi(q)\pmod*{p^k}}}\frac{\alpha(q)}{q}\left(1-\frac{\log q}{\log 3x}\right)^{-1/2}\\
&\le\sum_{\substack{p\le z\\p^k>w}}\log p^k\sum_{\substack{q\le x\\q\equiv-1\pmod*{p^k}}}\frac{1}{q}\left(1-\frac{\log q}{\log 3x}\right)^{-1/2}\\
&\ll\sum_{\substack{p\le z\\p^k>w}}\log p^k\sum_{\substack{q\le \sqrt{x}\\q\equiv-1\pmod*{p^k}}}\frac{1}{q}+\sum_{\substack{p\le z\\p^k>w}}\log p^k\sum_{\substack{\sqrt{x}<q\le x\\q\equiv-1\pmod*{p^k}}}\frac{1}{q}\left(1-\frac{\log q}{\log 3x}\right)^{-1/2}.
\end{align*}
The first double sum on the last line is clearly
\[\ll\sum_{\substack{p\le z\\p^k>w}}\frac{\log p^k}{\phi(p^k)}\log_2 x\ll\frac{(\log w)\log_2x}{w}\sum_{p\le z}1\ll\frac{z(\log w)\log_2x}{w\log z}\]
by Brun--Titchmarsh and partial summation. To estimate the second double sum, we split it into the two following subsums
\begin{align*}
& \sum_{\substack{p\le z\\w<p^k\le x^{1/3}}}\log p^k\sum_{\substack{\sqrt{x}<q\le x\\q\equiv-1\pmod*{p^k}}}\frac{1}{q}\left(1-\frac{\log q}{\log 3x}\right)^{-1/2},\\
&\sum_{\substack{p\le z\\p^k>x^{1/3}}}\log p^k\sum_{\substack{\sqrt{x}<q\le x\\q\equiv-1\pmod*{p^k}}}\frac{1}{q}\left(1-\frac{\log q}{\log 3x}\right)^{-1/2}.
\end{align*}
Since Brun--Titchmarsh and partial summation imply that
\[\sum_{\substack{\sqrt{x}<q\le x\\q\equiv-1\pmod*{p^k}}}\frac{1}{q}\left(1-\frac{\log q}{\log 3x}\right)^{-1/2}\ll\frac{1}{\phi(p^k)}\]
whenever $p^k\le x^{1/3}$, the first subsum above is
\[\ll\sum_{\substack{p\le z\\p^k>w}}\frac{\log p^k}{\phi(p^k)}\ll\frac{z(\log w)}{w\log z}.\]
On the other hand, the second subsum is 
\begin{align*}
&\ll(\log x)^{1/2}\sum_{\substack{p\le z\\p^k>x^{1/3}}}\log p^k\sum_{\substack{q\le x\\q\equiv-1\pmod*{p^k}}}\frac{1}{q} \\
&\ll(\log x)^{1/2}\log_2x\sum_{\substack{p\le z\\p^k>x^{1/3}}}\frac{\log p^k}{\phi(p^k)}\\
&\ll\frac{z(\log x)^{3/2}\log_2x}{x^{1/3}\log z}.
\end{align*}
Collecting the estimates above, we find
\[S_{\alpha}(x)^{-1}\sum_{f\le x}\alpha(f)\log S(L(f),z_f,w)\ll\frac{z(\log w)\log_2x}{w\log z}\ll\frac{z}{\log_4x},\]
from which it follows that for all but $o(S_{\alpha}(x))$ split-free $f\in\N\cap[1,x]$, we have $\log S(L(f),z_f,w)\le z$. Since 
\[\sum_{\substack{p\le z\\ p^k\le w}}\log p^k\ll z\left(\frac{\log w}{\log z}\right)^2\ll z,\]
we have 
\[\log S(L(f),z_f)\le\sum_{\substack{p\le z\\ p^k\le w}}\log p^k+\log S(L(f),z_f,w)\ll z\]
for all but $o(S_{\alpha}(x))$ split-free $f\in\N\cap[1,x]$. If we restrict to $f \in (\sqrt{x},x]$, then $z$ and $z_f$ are of the same order of magnitude, and the final  ``$\ll z$'' can be replaced with ``$\ll z_f$.'' This completes the proof of Lemma \ref{lem:eps}.
\end{proof}

A similar argument shows that $L(f)$ is rarely divisible by the square of a prime exceeding $z_f$.

\begin{lem}\label{lem:nolargesquares} For almost all split-free $f$, there is no prime $p > z_f$ for which $p^2 \mid L(f)$.
\end{lem}
\begin{proof} Observe that
\begin{align*}
S_{\alpha}(x)^{-1}\sum_{\substack{f\le x\\p^2\mid L(f)\text{~for some~}p>z_f}}\alpha(f)&\le S_{\alpha}(x)^{-1}\sum_{p>z/2}\sum_{\substack{\sqrt{x}<f\le x\\p^2\mid L(f)}}\alpha(f)+O\left(x^{-1/2}\right)\\
&\le S_{\alpha}(x)^{-1}\sum_{p>z/2}\sum_{\substack{f\le x\\p^2\mid L(f)}}\alpha(f)+O\left(x^{-1/2}\right),
\end{align*}
where $z=\log_2 x/\log_4 x$. The double sum on the last line is bounded above by
\[S_{\alpha}(x)^{-1}\sum_{p>z/2}\sum_{\substack{f\le x\\p^3\mid f}}\alpha(f)+S_{\alpha}(x)^{-1}\sum_{p>z/2}\sum_{\substack{q\le x\\q\equiv\chi(q)\pmod*{p^2}}}\sum_{\substack{f\le x\\ q\mid f}}\alpha(f).\]
Note that
\begin{align*}
S_{\alpha}(x)^{-1}\sum_{p>z/2}\sum_{\substack{f\le x\\p^3\mid f}}\alpha(f)&\ll\sum_{z/2<p\le x^{1/3}}\frac{1}{p^3}\left(1-\frac{\log p^3}{\log3x}\right)^{-1/2}\\
&\ll\sum_{z/2<p\le x^{1/4}}\frac{1}{p^3}+(\log x)^{1/2}\sum_{x^{1/4}<p\le x^{1/3}}\frac{1}{p^3}\\
&\ll \frac{1}{z^2\log z}
\end{align*}
and that
\begin{align*}
&\hspace{5.5mm}S_{\alpha}(x)^{-1}\sum_{p>z/2}\sum_{\substack{q\le x\\q\equiv\chi(q)\pmod*{p^2}}}\sum_{\substack{f\le x\\ q\mid f}}\alpha(f)\\
&\ll\sum_{z/2<p\le \sqrt{x}}\sum_{\substack{q\le x\\q\equiv\chi(q)\pmod*{p^2}}}\frac{\alpha(q)}{q}\left(1-\frac{\log q}{\log3x}\right)^{-1/2}\\
&\le \sum_{z/2<p\le \sqrt{x}}\sum_{\substack{q\le x\\q\equiv-1\pmod*{p^2}}}\frac{1}{q}\left(1-\frac{\log q}{\log3x}\right)^{-1/2}\\
&\ll\sum_{z/2<p\le \sqrt{x}}\sum_{\substack{q\le \sqrt{x}\\q\equiv-1\pmod*{p^2}}}\frac{1}{q}+\sum_{z/2<p\le \sqrt{x}}\sum_{\substack{\sqrt{x}<q\le x\\q\equiv-1\pmod*{p^2}}}\frac{1}{q}\left(1-\frac{\log q}{\log3x}\right)^{-1/2}\\
&\ll\sum_{z/2<p\le \sqrt{x}}\frac{\log_2x}{\phi(p^2)}+\sum_{z/2<p\le x^{1/6}}\frac{1}{\phi(p^2)}+(\log x)^{1/2}\sum_{x^{1/6}<p\le\sqrt{x}}\frac{\log_2x}{\phi(p^2)}\\
&\ll\frac{\log_2x}{z\log z}.
\end{align*}
So
\[S_{\alpha}(x)^{-1}\sum_{\substack{f\le x\\p^2\mid L(f)\text{~for some~}p>z_f}}\alpha(f)\ll\frac{\log_2x}{z\log z}\ll\frac{\log_4x}{\log_3x}. \qedhere\]
\end{proof}

The following lemma is a variant of \cite[Proposition 1]{LP03}.
\begin{lem}[under GRH]\label{lem:LiP} For almost all split-free $f$, the ratio $\frac{L(f)}{\ell(f)}$ is $z_f$-smooth.
\end{lem}

\begin{proof} Put
\[ z := \frac{1}{2} \frac{\log_2{x}}{\log_4{x}}. \] 
Since $z < z_f$ whenever $\sqrt{x} < f \le x$, it is enough to show all but $o(S_{\alpha}(x))$ split-free $f\le x$ have $L(f)/\ell(f)$ being $z$-smooth. 

We first show that the number of split-free $f\in\N\cap[1,x]$ divisible by an inert prime $p>z-1$ with $\ell(p)<2p^{1/2}/\log p$ is $o(S_{\alpha}(x))$. It is clear that the count is at most
\[\sum_{\substack{p\in(z-1,x]\cap\Pp_{-1}\\\ell(p)<2p^{1/2}/\log p}}\sum_{\substack{f\le x\\ p\mid f}}\alpha(f)\ll S_{\alpha}(x)\sum_{\substack{p\in(z-1,x]\cap\Pp_{-1}\\\ell(p)<2p^{1/2}/\log p}}\frac{1}{p}\left(1-\frac{\log p}{\log 3x}\right)^{-1/2}.\]
By Lemma \ref{lem:smallorderbound}, the number of inert primes $p\le t$ with $\ell(p)<2p^{1/2}/\log p$ is $O(t/\log^2 t)$, for all $t\ge 2$. By partial summation, we obtain
\[\sum_{\substack{p\in[2,t]\cap\Pp_{-1}\\\ell(p)<2p^{1/2}/\log p}}\frac{1}{p}=c+O\left(\frac{1}{\log t}\right),\]
where $c\ge0$ is constant. From this it follows by partial summation again that
\[\sum_{\substack{p\in(\sqrt{x},x]\cap\Pp_{-1}\\\ell(p)<2p^{1/2}/\log p}}\frac{1}{p}\left(1-\frac{\log p}{\log 3x}\right)^{-1/2}\ll\frac{1}{\sqrt{\log x}}.\]
Thus, the number of split-free $f\in\N\cap[1,x]$ divisible by an inert prime $p>z-1$ with $\ell(p)<2p^{1/2}/\log p$ is 
\[\ll S_{\alpha}(x)\sum_{\substack{p\in(z-1,\sqrt{x}]\cap\Pp_{-1}\\\ell(p)<2p^{1/2}/\log p}}\frac{1}{p}+S_{\alpha}(x)\sum_{\substack{p\in(\sqrt{x},x]\cap\Pp_{-1}\\\ell(p)<2p^{1/2}/\log p}}\frac{1}{p}\left(1-\frac{\log p}{\log 3x}\right)^{-1/2}\ll\frac{S_{\alpha}(x)}{\log z}.\]

Next, we claim that for each odd prime $q$, the number of split-free $f\in\N\cap[1,x]$ divisible by a prime $p\equiv -1\pmod{q}$ with $q^2/(4\log^2q)<p\le q^2\log^4q$ is $O\left(S_{\alpha}(x)\log_2q/\sqrt{\log q}\right)$. Evidently, the count is bounded above by
\[\sum_{\substack{q^2/(4\log^2q)<p\le\min(q^2\log^4q,x)\\p\equiv -1\pmod*{q}}}\sum_{\substack{f\le x\\ p\mid f}}\alpha(f)\ll S_{\alpha}(x)\sum_{\substack{q^2/(4\log^2q)<p\le\min(q^2\log^4q,x)\\p\equiv -1\pmod*{q}}}\frac{1}{p}\left(1-\frac{\log p}{\log 3x}\right)^{-1/2}.\]
By Brun--Titchmarsh and partial summation, we find that if $q\le x^{1/3}$, then the count is
\[\ll S_{\alpha}(x)\sum_{\substack{q^2/(4\log^2q)<p\le q^2\log^4q\\p\equiv -1\pmod*{q}}}\frac{1}{p}\ll\frac{S_{\alpha}(x)\log_2q}{q\log q};\]
if $q>x^{1/3}$, then the count is
\[\ll S_{\alpha}(x)(\log x)^{1/2}\sum_{\substack{q^2/(4\log^2q)<p\le q^2\log^4q\\p\equiv -1\pmod*{q}}}\frac{1}{p}\ll\frac{S_{\alpha}(x)\log_2q}{q\sqrt{\log q}}.\]
This proves our claim.

Moving on, we assert that for each odd prime $q$, the number of split-free $f\in\N\cap[1,x]$ divisible by an inert prime $p\equiv -1\pmod{q}$ with $p>q^2\log^4q$ and $q\mid \frac{p+1}{\ell(p)}$ is 
\[\ll S_{\alpha}(x)\left(\frac{1}{q\sqrt{\log q}}+\frac{\log_2x}{q^2}\right).\]
The count here is at most
\[\sum_{\substack{q^2\log^4q<p\le x\\p\in A_q}}\sum_{\substack{f\le x\\p\mid f}}\alpha(f)\ll S_{\alpha}(x)\sum_{\substack{q^2\log^4q<p\le x\\p\in A_q}}\frac{1}{p}\left(1-\frac{\log p}{\log 3x}\right)^{-1/2}.\]
where $A_q$ denotes the set of inert primes $p\equiv -1\pmod{q}$ with $q\mid \frac{p+1}{\ell(p)}$. By Proposition \ref{prop:hooleychen},
\[ \#(A_q\cap[2,y])\ll \frac{\Li(y)}{q^2} + y^{1/2}\log(qy) \]
for any $y\ge2$. In particular, we have $\#(A_q\cap[2,y])\ll \sqrt{y}\log q$ for $y\in[2,q^4\log^4q]$ and $\#(A_q\cap[2,y])\ll y/(q^2\log y)$ for $q^4\log^4q<y\le x$. Thus, if $q\le x^{1/8}$, then the count is
\[\ll S_{\alpha}(x)\left(\sum_{\substack{q^2\log^4q<p\le q^4\log^4q\\p\in A_q}}\frac{1}{p}+\sum_{\substack{q^4\log^4q<p\le x^{2/3}\\p\in A_q}}\frac{1}{p}+\sum_{\substack{x^{2/3}<p\le x\\p\in A_q}}\frac{1}{p}\left(1-\frac{\log p}{\log 3x}\right)^{-1/2}\right).\]
Since partial summation yields
\begin{align*}
 \sum_{\substack{q^2\log^4q<p\le q^4\log^4q\\p\in A_q}}\frac{1}{p}&\ll \frac{1}{q\log q},\\
 \sum_{\substack{q^4\log^4q<p\le x^{2/3}\\p\in A_q}}\frac{1}{p}&\ll\frac{\log_2x}{q^2},\\
 \sum_{\substack{x^{2/3}<p\le x\\p\in A_q}}\frac{1}{p}\left(1-\frac{\log p}{\log 3x}\right)^{-1/2}&\ll\frac{1}{q^2},
\end{align*}
it follows that the count is 
\[\ll S_{\alpha}(x)\left(\frac{1}{q\log q}+\frac{\log_2x}{q^2}\right).\]
If $q>x^{1/8}$, the count is 
\begin{align*}
&\ll S_{\alpha}(x)(\log x)^{1/2}\left(\sum_{\substack{q^2\log^4q<p\le q^4\log^4q\\p\in A_q}}\frac{1}{p}+\sum_{\substack{q^4\log^4q<p\le x\\p\in A_q}}\frac{1}{p}\right)\\
&\ll S_{\alpha}(x)(\log q)^{1/2}\left(\frac{1}{q\log q}+\frac{1}{q^2}\right)\\
&\ll\frac{S_{\alpha}(x)}{q\sqrt{\log q}}.
\end{align*}
This confirms our assertion.

We can now complete the proof of Lemma \ref{lem:LiP}. Suppose $q_f:=P^{+}(L(f)/\ell(f))> z$. Since $x$ can be assumed large, $q_f$ is unramified in $K$. As 
\[ L(f) = \lcm\{\psi(p^k): p^k \parallel f\}\quad\text{while}\quad \ell(f) = \lcm\{\ell(p^k): p^k\parallel f\}, \]
there is a prime power $p^k\parallel f$ which \[ v_{q_f}(\psi(p^k)) > v_{q_f}(\ell(p^k)).\] 
Then either $p=q_f$ and $k\ge 2$, implying $q_f^2\mid f$, or $p \in A_{q_f}$. The number of split-free $f\in\N\cap[1,x]$ with $q^2\mid f$ for some prime $q>z$ does not exceed
\begin{align*}
\sum_{\substack{q>z}}\sum_{\substack{f\le x\\q^2\mid f}}\alpha(f)&\ll S_{\alpha}(x)\sum_{z<q\le\sqrt{x}}\frac{1}{q^2}\left(1-\frac{\log q^2}{\log 3x}\right)^{-1/2}\\
&\ll S_{\alpha}(x)\sum_{z<q\le x^{1/3}}\frac{1}{q^2}+S_{\alpha}(x)(\log x)^{1/2}\sum_{x^{1/3}<q\le\sqrt{x}}\frac{1}{q^2}\\
&\ll\frac{S_{\alpha}(x)}{z\log z},
\end{align*}
which is negligible. It remains to estimate the number of split-free $f\in\N\cap[1,x]$  divisible by some prime $p\in A_{q_f}$, where $q_f>z$. The assumption $p\in A_{q_f}$ implies trivially that $p\ge q_f-1>z-1$. We have handled the case when $\ell(p)<2p^{1/2}/\log p$, so we may assume $\ell(p)\ge 2p^{1/2}\log p$. Since $\ell(p)\le(p+1)/q_f<2p/q_f$, we have $p>q_f^2/(4\log^2q_f)$. But the number of split-free $f\in\N\cap[1,x]$ divisible by a prime $p\in A_q$ with $p>q^2/(4\log^2q)$ for some prime $q>z$ is 
\[\ll S_{\alpha}(x)\sum_{q>z}\left(\frac{\log_2q}{q\sqrt{\log q}}+\frac{\log_2x}{q^2}\right)\ll S_{\alpha}(x)\left(\frac{\log_2z}{\sqrt{\log z}}+\frac{\log_2x}{z\log z}\right)\ll S_{\alpha}(x)\frac{\log_4x}{\sqrt{\log_3x}},\]
which is acceptable. This completes the proof of Lemma \ref{lem:LiP}.
\end{proof}

With Lemmas \ref{lem:eps}, \ref{lem:nolargesquares}, and \ref{lem:LiP} in hand, we can make short work of Proposition \ref{prop:keythm2}. Since
\[ \#\PrinCl(\Oo_f) = \frac{\psi(f)}{\ell(f)} = \frac{\psi(f)}{L(f)} \frac{L(f)}{\ell(f)}, \]
we have both
\[ \RadExp \PrinCl(\Oo_f) \ge \Rad \frac{\psi(f)}{L(f)}, \]
and
\[ \RadExp \PrinCl(\Oo_f) \le \left(\Rad \frac{\psi(f)}{L(f)}\right)  \left(\Rad \frac{L(f)}{\ell(f)}\right). \]
By Lemma \ref{lem:LiP},
\[ \Rad \frac{L(f)}{\ell(f)} \le \prod_{p \le z_f} p \le (\log{f})^{O(1/\log_4 f)} \]
for almost all split-free $f$. Proposition \ref{prop:keythm2}(a) follows by combining the last three displays and inserting the estimate of Proposition \ref{prop:radpsiL} for $\Rad\frac{\psi(f)}{L(f)}$.

Turning to (b), let $R = \frac{\Exp \PrinCl(\Oo_f)}{\RadExp \PrinCl(\Oo_f)}$. Then $\Rad(R) \mid \RadExp \PrinCl(\Oo_f)$, and 
\[ R \cdot \Rad(R) \mid \Exp \PrinCl(\Oo_f) \mid \Exp \PreCl(\Oo_f) \mid L(f). \]
Hence, $R$ is a divisor of $L(f)$ and $R$ is supported on primes $p$ for which $p^2 \mid L(f)$. By Lemmas \ref{lem:eps} and \ref{lem:nolargesquares}, $R\le (\log{f})^{O(1/\log_4 f)}$ for almost all split-free $f$. 

Finally, Proposition \ref{prop:keythm2}(c) is immediate from Lemmas \ref{lem:eps} and \ref{lem:LiP}.

\bibliographystyle{amsplain}
\bibliography{HFD}

\providecommand{\bysame}{\leavevmode\hbox to3em{\hrulefill}\thinspace}
\providecommand{\MR}{\relax\ifhmode\unskip\space\fi MR }
\providecommand{\MRhref}[2]{%
  \href{http://www.ams.org/mathscinet-getitem?mr=#1}{#2}
}
\providecommand{\href}[2]{#2}
\begin{thebibliography}{10}

\bibitem{AGP94}
{W.\,R.} Alford, A.~Granville, and C.~Pomerance, \emph{There are infinitely many {C}armichael numbers}, Ann. of Math. (2) \textbf{139} (1994), 703--722.

\bibitem{anderson97}
{D.\,F.} Anderson, \emph{Elasticity of factorizations in integral domains: a survey}, Factorization in integral domains ({I}owa {C}ity, {IA}, 1996), Lecture Notes in Pure and Appl. Math., vol. 189, Dekker, New York, 1997, pp.~1--29.

\bibitem{carlitz60}
L.~Carlitz, \emph{A characterization of algebraic number fields with class number two}, Proc. Amer. Math. Soc. \textbf{11} (1960), 391--392.

\bibitem{CC00}
{S.\,T.} Chapman and J.~Coykendall, \emph{Half-factorial domains, a survey}, Non-{N}oetherian commutative ring theory, Math. Appl., vol. 520, Kluwer Acad. Publ., Dordrecht, 2000, pp.~97--115.

\bibitem{chen02}
Y.-M.~J. Chen, \emph{On primitive roots of one-dimensional tori}, J. Number Theory \textbf{93} (2002), 23--33.

\bibitem{coykendall01}
J.~Coykendall, \emph{Half-factorial domains in quadratic fields}, J. Algebra \textbf{235} (2001), 417--430.

\bibitem{EBK69}
{P. van} Emde~Boas and D.~Kruyswijk, \emph{A combinatorial problem on finite {A}belian groups {III}}, Report ZW 1969-008, Amsterdam (1969), 32~p.

\bibitem{EPS}
P.~Erd\H{o}s, C.~Pomerance, and E.~Schmutz, \emph{Carmichael's lambda function}, Acta Arith. \textbf{58} (1991), 363--385.

\bibitem{SFTh23}
{K.\,(S.)} Fan, \emph{Rough numbers and variations on the {E}rdős--{K}ac theorem}, Dartmouth College Ph.D Dissertations. \textbf{156} (2023).

\bibitem{SF23}
\bysame, \emph{Weighted {E}rd{\H o}s--{K}ac theorems via computing moments}, Acta Arith. \textbf{217} (2025), 99--158.

\bibitem{FPextremal}
{K.\,(S.)} Fan and P.~Pollack, \emph{Extremal elasticity of quadratic orders}, The ideal theory and arithmetic of rings, monoids, and semigroups (Palermo, 2024) ({S.\,T.} Chapman, ed.), Contemp. Math., Amer. Math. Soc., Providence, RI, to appear; arXiv: \url{https://arxiv.org/abs/2503.07801}.

\bibitem{GHK06}
A.~Geroldinger and F.~Halter-Koch, \emph{Non-unique factorizations}, Pure and Applied Mathematics (Boca Raton), vol. 278, Chapman \& Hall/CRC, Boca Raton, FL, 2006.

\bibitem{HT88}
{R.\,R.} Hall and G.~Tenenbaum, \emph{Divisors}, Cambridge Tracts in Mathematics, vol.~90, Cambridge University Press, Cambridge, 1988.

\bibitem{HK72}
F.~Halter-Koch, \emph{Einseinheitengruppen und prime {R}estklassengruppen in quadratischen {Z}ahlk\"{o}rpern}, J. Number Theory \textbf{4} (1972), 70--77.

\bibitem{HK83}
\bysame, \emph{Factorization of algebraic integers}, Grazer Math. Berichte \textbf{191} (1983).

\bibitem{HK95}
\bysame, \emph{Elasticity of factorizations in atomic monoids and integral domains}, J. Th\'{e}or. Nombres Bordeaux \textbf{7} (1995), 367--385.

\bibitem{hooley67}
C.~Hooley, \emph{On {A}rtin's conjecture}, J. Reine Angew. Math. \textbf{225} (1967), 209--220.

\bibitem{kataoka03}
N.~Kataoka, \emph{The distribution of prime ideals in a real quadratic field with units having a given index in the residue class field}, J. Number Theory \textbf{101} (2003), 349--375.

\bibitem{Kou19}
D.~Koukoulopoulos, \emph{The distribution of prime numbers}, Graduate Studies in Mathematics, American Mathematical Society, Providence, RI, 2019.

\bibitem{LP03}
S.~Li and C.~Pomerance, \emph{On generalizing {A}rtin's conjecture on primitive roots to composite moduli}, J. Reine Angew. Math. \textbf{556} (2003), 205--224.

\bibitem{LP07}
Florian Luca and Carl Pomerance, \emph{Irreducible radical extensions and {E}uler-function chains}, Combinatorial number theory, de Gruyter, Berlin, 2007, pp.~351--361.

\bibitem{MC99}
P.~Moree and J.~Cazaran, \emph{On a claim of {R}amanujan in his first letter to {H}ardy}, Exposition. Math. \textbf{17} (1999), no.~4, 289--311.

\bibitem{narkiewicz95}
W.~Narkiewicz, \emph{A note on elasticity of factorizations}, J. Number Theory \textbf{51} (1995), 46--47.

\bibitem{Nor76}
K.~Norton, K.\, \emph{On the number of restricted prime factors of an integer. {I}}, Illinois J. Math. \textbf{20} (1976), 681--705.

\bibitem{pollack21}
P.~Pollack, \emph{The number of non-cyclic {S}ylow subgroups of the multiplicative group modulo {$n$}}, Canad. Math. Bull. \textbf{64} (2021), 204--215.

\bibitem{PP23}
\bysame, \emph{Half-factorial real quadratic orders}, Arch. Math. (Basel) \textbf{122} (2024), 491--500.

\bibitem{pollackm}
\bysame, \emph{Maximally elastic quadratic fields}, J. Number Theory \textbf{267} (2025), 80--100.

\bibitem{Pom77}
C.~Pomerance, \emph{On the distribution of amicable numbers}, J. Reine Angew. Math. \textbf{293(294)} (1977), 217--222.

\bibitem{Rib01}
P.~Ribenboim, \emph{Classical theory of algebraic numbers}, Universitext, Springer-Verlag, New York, 2001.

\bibitem{roskam00}
H.~Roskam, \emph{A quadratic analogue of {A}rtin's conjecture on primitive roots}, J. Number Theory \textbf{81} (2000), 93--109, erratum in \textbf{85} (2000), 108.

\bibitem{serre81}
J.-P. Serre, \emph{Quelques applications du th\'{e}or\`eme de densit\'{e} de {C}hebotarev}, Inst. Hautes \'{E}tudes Sci. Publ. Math. (1981), no.~54, 323--401.

\bibitem{steffan86}
J.-L. Steffan, \emph{Longueurs des d\'{e}compositions en produits d'\'{e}l\'{e}ments irr\'{e}ductibles dans un anneau de {D}edekind}, J. Algebra \textbf{102} (1986), 229--236.

\bibitem{Ten15}
G.~Tenenbaum, \emph{Introduction to analytic and probabilistic number theory}, 3rd ed., Graduate Studies in Mathematics, vol. 163, American Mathematical Society, Providence, RI, 2015.

\bibitem{valenza90}
{R.\,J.} Valenza, \emph{Elasticity of factorization in number fields}, J. Number Theory \textbf{36} (1990), 212--218.

\bibitem{weber82}
H.~Weber, \emph{{B}eweis des {S}atzes, dass jede eigentlich primitive quadratische {F}orm unendlich viele {P}rimzahlen darzustellen f\"{a}hig ist}, Math. Ann. \textbf{20} (1882), 301--329.

\bibitem{zaks76}
A.~Zaks, \emph{Half factorial domains}, Bull. Amer. Math. Soc. \textbf{82} (1976), 721--723, corrigendum in \textbf{82} (1976), 965.

\bibitem{zaks80}
\bysame, \emph{Half-factorial-domains}, Israel J. Math. \textbf{37} (1980), 281--302.

\end{thebibliography}
\end{document}